\documentclass[11pt]{article}

\usepackage{color}
\usepackage{latexsym}
\usepackage{amssymb}
\usepackage{amsmath, amsfonts,amssymb,theorem,euscript,array,enumerate,amsfonts,mathrsfs}
\usepackage{graphicx}

\newtheorem{Theorem}{Theorem}[part]
\newtheorem{Definition}{Definition}[part]
\newtheorem{Proposition}{Proposition}[part]

\newtheorem{Lemma}{Lemma}[part]

\newtheorem{Remark}{Remark}[part]

\def\esssup_#1{\underset{#1}{\mathrm{ess\,sup\, }}}
\def\essinf_#1{\underset{#1}{\mathrm{ess\,inf\, }}}

\def \trans{^{\scriptscriptstyle{\intercal}}}

\def \Prod{\displaystyle\prod}

\def \trans{^{\scriptscriptstyle{\intercal }}}

\def \N{\mathbb{N}}
\def \R{\mathbb{R}}

\def \E{\mathbb{E}}
\def \F{\mathbb{F}}
\def \G{\mathbb{G}}
\def \P{\mathbb{P}}

\def \S{\mathbb{S}}

\def \Ac{{\cal A}}
\def \Bc{{\cal B}}

\def \Ec{{\cal E}}
\def \Fc{{\cal F}}
\def \Gc{{\cal G}}

\def \Lc{{\cal L}}
\def \Pc{{\cal P}}

\def \Xc{{\cal X}}

\def \ni{\noindent}

\def \eps{\varepsilon}

\def \ep{\hbox{ }\hfill$\Box$}

\def\reff#1{{\rm(\ref{#1})}}

\def\beqs{\begin{eqnarray*}}
\def\enqs{\end{eqnarray*}}
\def\beq{\begin{eqnarray}}
\def\enq{\end{eqnarray}}

\allowdisplaybreaks

\begin{document}

\title{Backward SDE Representation for Stochastic Control Problems with Non Dominated Controlled Intensity\thanks{The authors would like to thank Prof. Huy\^en Pham for helpful discussions and suggestions related to this work.}}

\author{S\'ebastien CHOUKROUN\thanks{Laboratoire de Probabilit\'es et Mod\`eles Al\'eatoires, CNRS, UMR 7599, Universit\'e Paris Diderot,
		\sf  choukroun at math.univ-paris-diderot.fr.}~~~
        Andrea COSSO\thanks{Laboratoire de Probabilit\'es et Mod\`eles Al\'eatoires, CNRS, UMR 7599, Universit\'e Paris Diderot,
		\sf  cosso at math.univ-paris-diderot.fr.}}

\maketitle

\date{}

\begin{abstract}
We are interested in stochastic control problems coming from mathematical finance and, in particular, related to model uncertainty, where the uncertainty affects both volatility and intensity. This kind of stochastic control problems is associated to a fully nonlinear integro-partial differential equation, which has the peculiarity that the measure $(\lambda(a,\cdot))_a$ characterizing the jump part is not fixed but depends on a parameter $a$ which lives in a compact set $A$ of some Euclidean space $\R^q$. We do not assume that the family $(\lambda(a,\cdot))_a$ is dominated. Moreover, the diffusive part can be degenerate. Our aim is to give a BSDE representation, known as nonlinear Feynman-Kac formula, for the value function associated to these control problems. For this reason, we introduce a class of backward stochastic differential equations with jumps and partially constrained diffusive part. We look for the minimal solution to this family of BSDEs, for which we prove uniqueness and existence by means of a penalization argument. We then show that the minimal solution to our BSDE provides the unique viscosity solution to our fully nonlinear integro-partial differential equation.
\end{abstract}

\vspace{5mm}


\vspace{5mm}

\noindent {\bf 2010 Math Subject Classification:}  60H10,  60H30, 60G55, 60G57, 93E20.

\newpage

\section{Introduction}

Recently, \cite{khapha12} introduced a new class of backward stochastic differential equations (BSDEs) with nonpositive jumps in order to provide a probabilistic representation formula, known as nonlinear Feynman-Kac formula, for fully nonlinear integro-partial differential equations (IPDEs) of the following type (we use the notation $x.y$ to denote the scalar product in $\R^d$):
\begin{align}
\label{E:HJBIntro0}
&\frac{\partial v}{\partial t} + \sup_{a\in A}\bigg[b(x,a).D_x v + \frac{1}{2}\text{tr}\big(\sigma\sigma\trans (x,a)D_x^2 v\big) + f(x,a) \\
&+ \int_E\big(v(t,x+\beta(x,a,e))-v(t,x)-\beta(x,a,e).D_x v(t,x)\big)\lambda(de) \bigg] \ = \ 0, \qquad \text{on }[0,T)\times\R^d\!, \notag \\
&\hspace{9.4cm}v(T,x) \ = \ g(x), \qquad\quad\!\!\! x\in\R^d\!, \notag
\end{align}
where $A$ is a compact subset of $\R^q$, $E$ is a Borelian subset of $\R^k\backslash\{0\}$, and $\lambda$ is a nonnegative $\sigma$-finite measure on $(E,\Bc(E))$ satisfying the integrability condition $\int_E(1\wedge|e|^2)\lambda(de)<\infty$. Notice that in \cite{khapha12} more general equations than \eqref{E:HJBIntro0} are considered, where the function $f=f(x,a,v,\sigma\trans(x,a)D_x v)$ depends also on $v$ and its gradient $D_x v$. However, the case $f=f(x,a)$ is particularly relevant, as \eqref{E:HJBIntro0} turns out to be the Hamilton-Jacobi-Bellman equation of a stochastic control problem where the state process is a jump-diffusion with drift $b$, diffusion coefficient $\sigma$ (possibly degenerate), and jump size $\beta$, which are all controlled. A special case of \eqref{E:HJBIntro0} is the Hamilton-Jacobi-Bellman equation associated to the uncertain volatility model in mathematical finance, which takes the following form:
\begin{equation}
\label{E:HJBIntroG}
\frac{\partial v}{\partial t} + G(D_x^2 v) \ = \ 0, \quad \text{on }[0,T)\times\R^d, \qquad\qquad v(T,x) \ = \ g(x), \quad x\in\R^d,
\end{equation}
where $G(M)=\frac{1}{2}\sup_{c\in C}[c M]$ and $C$ is a set of symmetric nonnegative matrices of order $d$. As described in \cite{peng06}, the unique viscosity solution to \eqref{E:HJBIntroG} is represented in terms of the so-called $G$-Brownian motion $B$ under the nonlinear expectation $\Ec(\cdot)$ as follows:
\[
v(t,x) \ = \ \Ec\big(g(x+B_t)\big).
\]
It is however not clear how to simulate $G$-Brownian motion. On the other hand, when $C$ can be identified with a compact subset $A$ of a Euclidean space $\R^q$, we have the probabilistic representation formula presented in \cite{khapha12}, which can be implemented numerically as shown in \cite{khalanpha13a} and \cite{khalanpha13b}. We recall that the results presented in \cite{khapha12} were generalized to the case of controller-and-stopper games in \cite{choukcossopham13} and to non-Markovian stochastic control problems in \cite{fuhrmanpham13}.

In the present paper, our aim is to generalize the results presented in \cite{khapha12} providing a probabilistic representation formula for the unique viscosity solution to the following fully nonlinear integro-PDE of Hamilton-Jacobi-Bellman type:
\begin{align}
\label{E:HJBIntro}
&\frac{\partial v}{\partial t} + \sup_{a\in A}\bigg[b(x,a).D_x v + \frac{1}{2}\text{tr}\big(\sigma\sigma\trans (x,a)D_x^2 v\big) + f(x,a) \\
&+ \int_E\big(v(t,x+\beta(x,a,e))-v(t,x)-\beta(x,a,e).D_x v(t,x)\big)\lambda(a,de) \bigg] = 0, \quad \text{on }[0,T)\times\R^d\!, \notag \\
&\hspace{9.85cm}v(T,x) = g(x), \qquad\quad\!\! x\in\R^d\!, \notag
\end{align}
where $\lambda$ is a transition kernel from $(\R^q,\Bc(\R^q))$ into $(E,\Bc(E))$, namely $\lambda(a,\cdot)$ is a nonnegative measure on $(E,\Bc(E))$ for every $a\in\R^q$ and $\lambda(\cdot,E')$ is a Borel measurable function for every $E'\in\Bc(E)$. We do not assume that the family of measures $(\lambda(a,\cdot))_{a\in\R^q}$ is dominated. Moreover, the diffusion coefficient $\sigma$ can be degenerate.

A motivation to the study of equation \eqref{E:HJBIntro} comes from mathematical finance and, in particular, from model uncertainty, when uncertainty affects both volatility and intensity. This topic was studied by means of second order BSDEs with jumps (2BSDEJs) in \cite{kaziposszhou12I} and \cite{kaziposszhou12II}, to which we refer also for the wellposedness of these kinds of backward equations. Model uncertainty is also strictly related to the theory of $G$-L\'evy processes and, more generally, of nonlinear L\'evy processes, see \cite{hu_peng09} and \cite{neufeld_nutz14}. In this case, the associated fully nonlinear integro-PDE, which naturally generalizes equation \eqref{E:HJBIntroG}, takes the following form:
\begin{align}
\label{E:HJBIntroLevy}
&\frac{\partial v}{\partial t} + \sup_{(b,c,F)\in\Theta}\bigg[b.D_x v + \frac{1}{2}\text{tr}\big(c D_x^2 v\big) \\
&+ \int_E\big(v(t,x+z)-v(t,x)-D_x v(t,x).z1_{\{|z|\leq1\}}\big)F(dz) \bigg] \ = \ 0, \qquad \text{on }[0,T)\times\R^d\!, \notag \\
&\hspace{8.1cm}v(T,x) \ = \ g(x), \qquad\quad\!\!\! x\in\R^d\!, \notag
\end{align}
where $\Theta$ denotes a set of L\'evy triplets $(b,c,F)$; here $b$ is a vector in $\R^d$, $c$ is a symmetric nonnegative matrix of order $d$, and $F$ is a L\'evy measure on $(\R^d,\Bc(\R^d))$. From  \cite{hu_peng09} and \cite{neufeld_nutz14}, we know that the unique viscosity solution to equation \eqref{E:HJBIntroLevy} is represented in terms of the so-called nonlinear L\'evy process $\Xc$ under the nonlinear expectation $\Ec(\cdot)$ as follows:
\[
v(t,x) \ = \ \Ec(g(x+\Xc_t)).
\]
If we are able to describe the set $\Theta$ by means of a parameter $a$ which lives in a compact set $A$ of an Euclidean space $\R^q$, then \eqref{E:HJBIntroLevy} can be written in the form \eqref{E:HJBIntro}. Therefore, $v$ is also given by our probabilistic representation formula, in which the forward process is possibly easier to simulate than a nonlinear L\'evy process.

More generally, we expect that the viscosity solution $v$ to equation \eqref{E:HJBIntro} should represent the value function of a stochastic control problem where, roughly speaking, the state process $X$ is a jump-diffusion process, which has the peculiarity that we may control the dynamics of $X$ changing its jump intensity, other than acting on the coefficients $b$, $\sigma$, and $\beta$ of the SDE solved by $X$. We refer to this problem as a stochastic optimal control problem with (non dominated) controlled intensity. Unfortunately, we did not find any reference in the literature for this kind of stochastic control problem. For this reason, and also because it will be useful to understand the general idea behind the derivation of our nonlinear Feynman-Kac formula, we describe it here, even if only formally. Let $(\bar\Omega,\bar\Fc,\bar\P)$ be a complete probability space satisfying the usual conditions on which a $d$-dimensional Brownian motion $\bar W$ $=$ $(\bar W_t)_{t\geq0}$ is defined. Let $\bar\F=(\bar\Fc_t)_{t\geq0}$ denote the usual completion of the natural filtration generated by $\bar W$ and $\bar\Ac$ the class of control processes $\alpha$, i.e., of $\bar\F$-predictable processes valued in $A$. Let also $\Omega'$ be the canonical space of the marked point process on $\R_+\times E$ (see Section \ref{S:Notation} below for a definition), with canonical right-continuous filtration $\F'$ and canonical random measure $\pi'$. Then, consider $(\Omega,\Fc,\F=(\Fc_t)_{t\geq0})$ defined as $\Omega$ $:=$ $\bar\Omega\times\Omega'$, $\Fc$ $:=$ $\bar\Fc\otimes\Fc_\infty'$, and  $\Fc_t$ $:=$ $\cap_{s>t}\bar\Fc_s\otimes\Fc_s'$. Moreover, we set $W(\omega)$ $:=$ $\bar W(\bar\omega)$, $\pi(\omega,\cdot)$ $:=$ $\pi'(\omega',\cdot)$, and $\Ac:=\{\alpha\colon\alpha(\omega)=\bar\alpha(\bar\omega),\,\forall\,\omega\in\Omega,\text{ for some }\bar\alpha\in\bar\Ac\}$. Suppose that for every $\alpha\in\Ac$ we are able to construct a measure $\P^\alpha$ on $(\Omega,\Fc)$ such that $W$ is a Brownian motion and $\pi$ is an integer-valued random measure with compensator $\lambda(\alpha_t,de)dt$ on $(\Omega,\Fc,\F,\P^\alpha)$. Then, consider the stochastic control problem with value function given by
($\E^\alpha$ denotes the expectation with respect to $\P^\alpha$)
\beq
\label{E:value_function}
v(t,x) &:=& \sup_{\alpha\in\Ac} \E^\alpha\bigg[\int_t^T f(X_s^{t,x,\alpha},\alpha_s)ds + g(X_T^{t,x,\alpha})\bigg],
\enq
where $X^{t,x,\alpha}$ has the controlled dynamics on $(\Omega,\Fc,\F,\P^\alpha)$
\beqs
dX_s^\alpha &=& b(X_s^\alpha,\alpha_s)ds + \sigma(X_s^\alpha,\alpha_s)dW_s + \int_E \beta(X_{s^-}^\alpha,\alpha_s,e)\tilde\pi(ds,de)
\enqs
starting from $x$ at time $t$, with $\tilde\pi(dt,de)$ $=$ $\pi(dt,de)$ $-$ $\lambda(\alpha_t,de)dt$ the compensated martingale measure of $\pi$. As mentioned above, even if we do not address this problem here, we expect that the above partial differential equation \eqref{E:HJBIntro} turns out to be the dynamic programming equation of the stochastic control problem with value function formally given by \eqref{E:value_function}. Having this in mind, we can now begin to describe the intuition, inspired by \cite{khaetal10} and \cite{khapha12}, behind the derivation of our Feynman-Kac representation formula for the HJB equation \eqref{E:HJBIntro} in terms of a forward backward stochastic differential equation (FBSDE).

The fundamental idea concerns the \emph{randomization} of the control, which is achieved introducing on $(\bar\Omega,\bar\Fc,\bar\P)$ a $q$-dimensional Brownian motion $\bar B=(\bar B_t)_{t\geq0}$, independent of $\bar W$. Now $\bar\F$ denotes the usual completion of the natural filtration generated by $\bar W$ and $\bar B$. We also set $B(\omega):=\bar B(\bar\omega)$, for all $\omega\in\Omega$, so that $B$ is defined on $\Omega$. Since the control lives in the compact set $A\subset\R^q$, we can not use directly $B$ to randomize the control, but we need to map $B$ on $A$. More precisely, we shall assume the existence of a surjection $h\colon\R^d\rightarrow A$ satisfying $h\in C^2(\R^d;A)$ (e.g., the existence of such a function $h$ is guaranteed when $A$ is a ball in $\R^q$). Then, for every $(t,x,a)\in[0,T]\times\R^d\times\R^q$, we consider the forward stochastic differential equation in $\R^d\times\R^q$:
\beq
X_s &=& x + \int_t^s b(X_r,I_r) dr + \int_t^s \sigma(X_r,I_r) dW_r + \int_t^s\int_E \beta(X_{r^-},I_r,e) \tilde\pi(dr,de), \label{FSDEX_Intro}\\
I_s &=& h(a + B_s - B_t), \label{FSDEI_Intro}
\enq
for all $t\leq s\leq T$, where $\tilde\pi(ds,de) = \pi(ds,de) - \lambda(I_s,de)ds$ is the compensated martingale measure of $\pi$, which is an integer-valued random measure with compensator $\lambda(I_s,de)ds$. Unlike \cite{khapha12}, we used a Brownian motion $B$ to randomize the control, instead of a Poisson random measure $\mu$ on $\R_+\times A$. From one hand, the Poisson random measure turns out to be more convenient to deal with a general compact set $A$, since $\mu$ is already supported by $\R_+\times A$, so that we do not have to impose the existence of a surjection $h$ from the entire space $\R^q$ onto $A$, as we did here. On the other hand, the choice of a Brownian motion $B$ is more convenient to derive a martingale representation theorem for our model. Indeed, in contrast with \cite{khapha12}, the intensity of the measure $\pi$ depends on the process $I$, therefore it is natural to expect a dependence between $\pi$ and the noise used to randomize the control. The advantage of $B$ with respect to $\mu$ is given by the fact that $B$ is \emph{orthogonal} to $\pi$, since $B$ is a continuous process (see the bottom of page 183 in \cite{jacshiryaev03} for a definition of orthogonality between a martingale and a random measure). Thanks to this orthogonality we are able to derive a martingale representation theorem in our context, which is essential for the derivation of our nonlinear Feynman-Kac representation formula.

Let us focus on the form of the stochastic differential equation \eqref{FSDEX_Intro}-\eqref{FSDEI_Intro}. We observe that the jump part of the driving factors in \eqref{FSDEX_Intro} is not given, but depends on the solution via its intensity. This makes the SDE \eqref{FSDEX_Intro}-\eqref{FSDEI_Intro} nonstandard. These kinds of equations were firstly studied in \cite{jacprotter82} and have also been used in the financial literature, see e.g. \cite{bechschw05}, \cite{crepey11}, \cite{crepeymat08}, \cite{cuchfilmayteich11}, \cite{filovschm11}. Notice that in \cite{bechschw05}, \cite{crepey11}, and \cite{crepeymat08}, $\lambda$ is absolutely continuous with respect to a given deterministic measure on $(E,\Bc(E))$, which allows to solve \eqref{FSDEX_Intro}-\eqref{FSDEI_Intro} bringing it back to a standard SDE, via a change of intensity ``\`{a} la Girsanov''. On the other hand, in the present paper, we shall tackle the above SDE solving firstly equation \eqref{FSDEI} for any $(t,a)\in[0,T]\times\R^q$, then constructing a probability measure $\P^{t,a}$ on $(\Omega,\Fc)$ such that the random measure $\pi(dt,de)$ admits $\lambda(I_s^{t,a},de)ds$ as compensator, and finally addressing \eqref{FSDEX}. In the appendix, we also prove additional properties of $\pi$ and $(X,I)$. More precisely, we present a characterization of $\pi$ in terms of Fourier and Laplace functionals, which shows that $\pi$ is a conditionally Poisson random measure (also known as doubly stochastic Poisson random measure or Cox random measure) relative to $\sigma(I_s^{t,a};s\geq0)$. Moreover, we study the Markov properties of the pair $(X,I)$.

Regarding the backward stochastic differential equation, as expected, it is driven by the Brownian motions $W$ and $B$, and by the random measure $\pi$, namely it is a BSDE with jumps with terminal condition $g(X_T^{t,x,a})$ and generator $f(X_\cdot^{t,x,a},I_\cdot^{t,a})$, as it is natural from the expression of the HJB equation \eqref{E:HJBIntro}. The backward equation is also characterized by a constraint on the diffusive part relative to $B$, which turns out to be crucial and entails the presence of an increasing process in the BSDE. In conclusion, for any $(t,x,a)\in[0,T]\times\R^d\times\R^q$, the backward stochastic differential equation has the following form:
\beq
\label{E:BSDEIntro}
Y_s & = & g(X_T^{t,x,a}) + \int_s^T f(X_r^{t,x,a},I_r^{t,a}) dr + K_T - K_s - \int_s^T Z_r dW_r \notag \\
& & - \int_s^T V_r dB_r - \int_s^T\int_E U_r(e) \tilde\pi(dr,de), \qquad   t \leq s \leq T,\,\P^{t,a}\,a.s.
\enq
and
\beq
\label{E:JumpConstrIntro}
|V_s| & = & 0 \qquad\qquad ds\otimes d\P^{t,a}\,a.e.
\enq
We refer to \eqref{E:BSDEIntro}-\eqref{E:JumpConstrIntro} as backward stochastic differential equation with jumps and partially constrained diffusive part. Notice that the presence of the increasing process $K$ in the backward equation does not guarantee the uniqueness of the solution. For this reason, we look only for the minimal solution $(Y,Z,V,U,K)$ to the above BSDE, in the sense that for any other solution $(\bar Y,\bar Z,\bar V,\bar U,\bar K)$ we must have $Y\leq\bar Y$. The existence of the minimal solution is based on a penalization approach as in \cite{khapha12}. We can now write down the nonlinear Feynman-Kac formula:
\[
v(t,x,a) \ := \ Y_t^{t,x,a}, \qquad (t,x,a)\in[0,T]\times\mathbb{R}^d\times\mathbb{R}^q.
\]
Observe that the function $v$ should not depend on $a$, but only on $(t,x)$. The function $v$ turns out to be independent of the variable $a$ as a consequence of the constraint \eqref{E:JumpConstrIntro}. Indeed, if $v$ were regular enough, then, for any $(t,x,a)\in[0,T]\times\R^d\times\R^q$, we would have
\[
V_s^{t,x,a} \ = \ D_hv(s,X_s^{t,x,a},I_s^{t,a})D_ah(a+B_s-B_t) \ = \ 0, \qquad ds\otimes d\P^{t,a}\,a.e.
\]
This would imply (see Subsection \ref{SubS:NonDep_a} below) that $v$ does not depend on its last argument. However, we do not know in general if the function $v$ is so regular in order to justify the previous passages. Therefore, the rigorous proof relies on viscosity solutions arguments. In the end, we prove that the function $v$ does not depend on the variable $a$ in the interior $\mathring{A}$ of $A$ and admits the following probabilistic representation formula:
\[
v(t,x) \ := \ Y_t^{t,x,a}, \qquad (t,x)\in[0,T]\times\mathbb{R}^d,
\]
for any $a\in\mathring{A}$. Moreover, $v$ is a viscosity solution to \eqref{E:HJBIntro}. Actually, $v$ is the unique viscosity solution to \eqref{E:HJBIntro}, as it follows from the comparison theorem proved in the Appendix. Notice that, due to the presence of the non dominated family of measures $(\lambda(a,\cdot))_{a\in A}$, we did not find in the literature a comparison theorem for viscosity solution to our equation \eqref{E:HJBIntro}. For this reason, we prove it in the Appendix, even though the main ideas are already contained in the paper \cite{barlesimbert08}, in particular the remarkable Jensen-Ishii's lemma for integro-partial differential equations.

The  rest of the paper is organized as follows. Section 2 introduces some notations and studies the construction of the solution to the forward equation \eqref{FSDEX_Intro}-\eqref{FSDEI_Intro}. Section 3 gives a detailed formulation of the BSDE with jumps and partially constrained diffusive part. In particular, Subsection 3.1 is devoted to the existence of the minimal solution to our BSDE  by a penalization approach. Section 4 makes the connection between the minimal solution to our BSDE and equation \eqref{E:HJBIntro}. In the Appendix, we prove a martingale representation theorem for our model, we collect some properties of the random measure $\pi$ and of the pair $(X,I)$, and we prove a comparison theorem for equation \eqref{E:HJBIntro}.

\section{Notations and preliminaries}
\label{S:Notation}

\setcounter{equation}{0} \setcounter{Assumption}{0}
\setcounter{Theorem}{0} \setcounter{Proposition}{0}
\setcounter{Corollary}{0} \setcounter{Lemma}{0}
\setcounter{Definition}{0} \setcounter{Remark}{0}

Let $(\bar\Omega,\bar\Fc,\bar\P)$ be a complete probability space satisfying the usual conditions on which are defined a $d$-dimensional Brownian motion $\bar W$ $=$ $(\bar W_t)_{t\geq0}$ and an independent $q$-dimensional Brownian motion $\bar B$ $=$ $(\bar B_t)_{t\geq0}$. We will always assume that $\bar\F=(\bar\Fc_t)_{t\geq0}$ is the usual completion of the natural filtration generated by $\bar W$ and $\bar B$. Let us introduce some additional notations.
\begin{enumerate}
\item[(i)] $\Omega'$ is the set of sequences $\omega'=(t_n,e_n)_{n\in\N}\subset(0,\infty]\times E_\Delta$, where $E_\Delta$ $=$ $E\cup\{\Delta\}$ and $\Delta$ is an external point of $E$. Moreover $t_n<\infty$ if and only if $e_n\in E$, and when $t_n<\infty$ then $t_n<t_{n+1}$. $\Omega'$ is equipped with the canonical marked point process $(T_n',\alpha_n')_{n\in\N}$, with associated canonical random measure $\pi'$, defined as
    \beqs
    T_n'(\omega') \ = \ t_n, \qquad \alpha_n'(\omega') \ = \ e_n
    \enqs
    and
    \beqs
    \pi'(\omega',dt,de) \ = \ \sum_{n\in\N} 1_{\{T_n'(\omega')<\infty\}} \delta_{(T_n'(\omega'),\alpha_n'(\omega'))}(dt,de),
    \enqs
    where $\delta_x$ denotes the Dirac measure at point $x$. Set $T_\infty'$ $:=$ $\lim_n T_n'$. Finally, define $\F'=(\Fc_s)_{t\geq0}$ as $\Fc_t=\cap_{s>t}\Gc_s$, where $\G'=(\Gc_s)_{t\geq0}$ is the canonical filtration, given by $\Gc_s=\sigma(\pi'(\cdot,F)\colon F\in\Bc([0,t])\otimes\Bc(E))$.
\item[(ii)] $(\Omega,\Fc,\F=(\Fc_t)_{t\geq0})$ is such that $\Omega$ $:=$ $\bar\Omega\times\Omega'$, $\Fc$ $:=$ $\bar\Fc\otimes\Fc_\infty'$, and  $\Fc_t$ $:=$ $\cap_{s>t}\bar\Fc_s\otimes\Fc_s'$. Moreover, we set $W(\omega)$ $:=$ $\bar W(\bar\omega)$, $B(\omega)$ $:=$ $\bar B(\bar\omega)$, and $\pi(\omega,\cdot)$ $:=$ $\pi'(\omega',\cdot)$. Finally, we set also $T_n(\omega)$ $:=$ $T_n'(\omega')$, $\alpha_n(\omega)$ $:=$ $\alpha_n'(\omega')$, and $T_\infty(\omega)$ $:=$ $T_\infty'(\omega')$.
\end{enumerate}

Let $\Pc_\infty$ denote the $\sigma$-field of $\F$-predictable subsets of $\R_+\times\Omega$. We recall that a random measure $\pi$ on $\R_+\times E$ is a transition kernel from $(\Omega,\Fc)$ into $(\R_+\times E,\Bc(\R_+)\otimes\Bc(E))$, satisfying $\pi(\omega,\{0\}\times E)=0$ for all $\omega\in\Omega$; moreover, an integer-valued random measure $\pi$ on $\R_+\times E$ is an optional and $\Pc_\infty\otimes\Bc(E)$-$\sigma$-finite, $\N\cup\{+\infty\}$-valued random measure such that $\pi(\omega,\{t\}\times E)\leq1$ for all $(t,\omega)\in[0,T]\times\Omega$, see Definition 1.13, Chapter II, in \cite{jacshiryaev03}.

\vspace{1mm}

We are given some measurable functions $b\colon\R^d\times \R^q \rightarrow \R^d$, $\sigma\colon\R^d\times \R^q \rightarrow \R^{d\times d}$, and $\beta\colon\R^d\times\R^q\times E\rightarrow\R^d$, where $E$ is a Borelian subset of $\R^k\backslash\{0\}$, equipped with its Borel $\sigma$-field $\Bc(E)$. Moreover, let $\lambda$ be a transition kernel from $(\R^q,\Bc(\R^q))$ into $(E,\Bc(E))$, namely $\lambda(a,\cdot)$ is a nonnegative measure on $(E,\Bc(E))$ for every $a\in\R^q$ and $\lambda(\cdot,E')$ is a Borel measurable function for every $E'\in\Bc(E)$. Furthermore, let $A$ be a compact subset of $\R^q$ such that there exists a surjection $h\colon\R^d\rightarrow A$ satisfying $h\in C^2(\R^d;A)$

\begin{Remark}
\rm{
The existence of such a function $h$ is guaranteed for the case $A=B_r(a)$, the ball of radius $r>0$ centered in $a\in\R^q$. As a matter of fact, consider the ball $B_1(0)$ of radius $1$ centered at zero. Define $\tilde h\colon\R_+\rightarrow[0,1]$ as follows
\[
\tilde h(\rho) \ = \
\begin{cases}
6\rho^5 - 15\rho^4 + 10\rho^3, \qquad &0\leq\rho\leq1, \\
1, &\rho>1.
\end{cases}
\]
Notice that $\tilde h(0)=0$ and $\tilde h(1)=1$, moreover $\tilde h'(0)=\tilde h''(0)=0$ and $\tilde h'(1)=\tilde h''(1)=0$. Then, we define $h(a)=\frac{a}{|a|}\tilde h(|a|)$, for $a\neq0$, and $h(0)=0$. In particular, we have
\[
h(a) \ = \ \big(6|a|^4 - 15|a|^3 + 10|a|^2\big)a1_{\{|a|\leq 1\}} + \frac{a}{|a|}1_{\{|a|>1\}},
\]
for all $a\in\R^q$.
\ep
}
\end{Remark}

For any $t\in[0,T]$ and $(x,a)\in\R^d\times\R^q$, we consider the forward stochastic differential equation in $\R^d\times\R^q$:
\beq
X_s &=& x + \int_t^s b(X_r,I_r) dr + \int_t^s \sigma(X_r,I_r) d W_r + \int_t^s\int_E \beta(X_{r^-},I_r,e) \tilde\pi(dr,de), \label{FSDEX}\\
I_s &=& h(a + B_s - B_t), \label{FSDEI}
\enq
for all $t\leq s\leq T$, where $\tilde\pi(ds,de) = \pi(ds,de) - \lambda(I_s,de)ds$ is the compensated martingale measure of $\pi$, which is an integer-valued random measure with compensator $\lambda(I_s,de)ds$.

As noticed in the introduction, the above SDE \eqref{FSDEX}-\eqref{FSDEI} is nonstandard, in the sense that the jump part of the driving factors in \eqref{FSDEX} is not given, but depends on the solution via its intensity. When the intensity $\lambda$ is absolutely continuous with respect to a given deterministic measure on $(E,\Bc(E))$, as in \cite{bechschw05}, \cite{crepey11}, and \cite{crepeymat08}, we can obtain \eqref{FSDEX}-\eqref{FSDEI} starting from a standard SDE via a change of intensity ``\`{a} la Girsanov''. On the other hand, in the present paper, we shall tackle the above SDE solving firstly equation \eqref{FSDEI}, then constructing the random measure $\pi(dt,de)$, and finally addressing \eqref{FSDEX}. The nontrivial part is the construction of $\pi$, which is essentially based on Theorem 3.6 in \cite{jac75}, and also on similar results in \cite{filovschm11}, Theorem 5.1, and \cite{cuchfilmayteich11}, Theorem A.4. Let us firstly introduce the following assumptions on the forward coefficients.

\vspace{3mm}

\ni {\bf (HFC)}
\begin{itemize}
\item[(i)] There exists a constant $C$ such that
\beqs
|b(x,a)-b(x',a')| + |\sigma(x,a)-\sigma(x',a')| &\leq& C \big(|x-x'| + |a-a'|\big),
\enqs
for all $x,x'\in\R^d$ and $a,a'\in\R^q$.
\item[(ii)] There exists a constant $C$ such that
\beqs
|\beta(x,a,e)| & \leq & C(1+|x|)(1\wedge|e|), \\
|\beta(x,a,e)-\beta(x',a',e)| & \leq & C \big(|x-x'| + |a-a'|\big)(1\wedge|e|),
\enqs
for all $x,x'\in\R^d$, $a,a'\in\R^q$, and $e\in E$.
\item[(iii)] The following integrability condition holds:
\beqs
\sup_{|a|\leq m} \int_E \big(1\wedge|e|^2\big) \lambda(a,de) & < & \infty, \qquad \forall\,m\in\N.
\enqs
\end{itemize}

Inspired by \cite{jacprotter82}, we give the definition of weak solution to equation \eqref{FSDEX}-\eqref{FSDEI}.

\begin{Definition}
\label{D:WeakSol}
A \textbf{weak solution} to equation \eqref{FSDEX}-\eqref{FSDEI} with initial condition $(t,x,a)\in[0,T]\times\R^d\times\R^q$ is a probability measure $\P$ on $(\Omega,\Fc)$ satisfying:
\begin{enumerate}
\item[\textup{(i)}] $\P(d\omega)$ $=$ $\bar\P(d\bar\omega)\otimes\P'(\bar\omega,d\omega')$, for some transition kernel $\P'$ from $(\bar\Omega,\bar\Fc)$ into $(\Omega',\Fc_\infty')$.
\item[\textup{(ii)}] Under $\P$, $\pi$ is an integer-valued random measure on $\R_+\times E$ with $\F$-compensator $1_{\{s<T_\infty\}}\lambda(I_s,de)ds$ and compensated martingale measure given by $\tilde\pi(ds,de)$ $=$ $\pi(ds,de)$ $-$ $1_{\{s<T_\infty\}}\lambda(I_s,de)ds$.
\item[\textup{(iii)}] We have
\beqs
X_s &=& x + \int_t^s b(X_r,I_r) dr + \int_t^s \sigma(X_r,I_r) dW_r + \int_t^s \int_E \beta(X_{r^-},I_r,e) \tilde\pi(dr,de), \\
I_s &=& h(a + B_s - B_t),
\enqs
for all $t \leq s \leq T$, $\P$ almost surely. Moreover, $(X_s,I_s)$ $=$ $(x,h(a))$ for $s<t$, and $(X_s,I_s)$ $=$ $(X_T,I_T)$ for $s>T$.
\end{enumerate}
\end{Definition}

Consider a probability measure $\P$ on $(\Omega,\Fc)$ satisfying condition (i) of Definition~\ref{D:WeakSol}. For every $(t,a)\in[0,T]\times\R^q$ let us denote $I^{t,a}$ $=$ $\{I_s^{t,a},\,s\geq0\}$ the unique process on $(\Omega,\Fc,\F,\P)$ satisfying $I_s^{t,a}=h(a+B_s-B_t)$ on $[t,T]$, with $I_s^{t,a}$ $=$ $h(a)$ for $s<t$ and $I_s^{t,a}$ $=$ $I_T^{t,a}$ for $s>T$. We notice that the notation $I^{t,a}$ can be misleading, since $a$ is not the initial point of $I^{t,a}$ at time $t$, indeed $I_t^{t,a}=h(a)$. Now we proceed to the construction of a probability measure on $(\Omega,\Fc)$ for which conditions (i) and (ii) of Definition~\ref{D:WeakSol} are satisfied. This result is based on Theorem 3.6 in \cite{jac75}, and we borrow also some ideas from \cite{filovschm11}, Theorem 5.1, and \cite{cuchfilmayteich11}, Theorem A.4.

\begin{Lemma}
\label{L:pi}
Under assumption \textup{\textbf{(HFC)}}, for every $(t,a)\in[0,T]\times\R^q$ there exists a unique probability measure on $(\Omega,\Fc)$, denoted by $\P^{t,a}$, satisfying conditions (i) and (ii) of Definition~\ref{D:WeakSol}, and also condition (ii)' given by:
\begin{enumerate}
\item[\textup{(ii)'}] $1_{\{s<T_\infty\}}\lambda(I_s^{t,a},de)ds$ is the $(\bar\Fc\otimes\Fc_s')_{s\geq0}$-compensator of $\pi$.
\end{enumerate}
\end{Lemma}
\textbf{Proof.}
The proof is essentially based on Theorem 3.6 in \cite{jac75}, after a reformulation of our problem in the setting of \cite{jac75}, which we now detail. Let $\hat\F=(\hat\Fc_s)_{s\geq0}$ where $\hat\Fc_s$ $:=$ $\bar\Fc\otimes\Fc_s'$. Notice that in $\hat\Fc_s$ we take $\bar\Fc$ instead of $\bar\Fc_s$. Indeed, in \cite{jac75} the $\sigma$-field $\bar\Fc$ represents the past information and is fixed throughout (we come back to this point later). Take $(t,a)\in[0,T]\times\R^q$ and consider the process $I^{t,a}$ $=$ $(I_s^{t,a})_{s\geq0}$. Set
\beqs
\nu(\omega,F) &=& \int_F 1_{\{s<T_\infty(\omega)\}} \lambda(I_s^{t,a}(\omega),de)ds
\enqs
for any $\omega\in\Omega$ and any $F\in\Bc(\R_+)\otimes\Bc(E)$. Now we show that $\nu$ satisfies the properties required in order to apply Theorem 3.6 in \cite{jac75}. In particular, since $\lambda$ is a transition kernel, we see that $\nu$ is a transition kernel from $(\Omega,\Fc)$ into $(\R_+\times E,\Bc(\R_+)\otimes\Bc(E))$; moreover, $\nu(\omega,\{0\}\times E)=0$ for all $\omega\in\Omega$, therefore $\nu$ is a random measure on $\R_+\times E$. Furthermore, for every $E'\in\Bc(E)$, the process $\nu((0,\cdot]\times E')$ $=$ $(\nu((0,s]\times E'))_{s\geq0}$ is $\hat\F$-predictable, hence $\nu$ is an $\hat\F$-predictable random measure. In addition, $\nu(\{s\}\times E)\leq1$, indeed $\nu$ is absolutely continuous with respect to the Lebesgue measure $ds$ and therefore $\nu(\{s\}\times E)=0$. Finally, we see by definition that $\nu([T_\infty,\infty)\times E)=0$. In conclusion, it follows from Theorem 3.6 in \cite{jac75} that there exists a unique probability measure on $(\Omega,\Fc)$, denoted by $\P^{t,a}$, satisfying condition (i) of Definition~\ref{D:WeakSol}, and for which $\nu$ is the $\hat\F$-compensator of $\pi$, i.e., the process
\beq
\label{E:nu-pi}
\big(\nu((0,s\wedge T_n]\times E') - \pi((0,s\wedge T_n]\times E')\big)_{s\geq0}
\enq
is a $(\P^{t,a},\hat\F)$-martingale, for any $E'\in\Bc(E)$ and any $n\in\N$. Therefore condition (ii)' is also satisfied.

To conclude, we need to prove that $\nu$ is also the $\F$-compensator of $\pi$. Since $\nu$ is an $\F$-predictable random measure, it follows from (2.6) in \cite{jac75} that it remains to prove that the process \eqref{E:nu-pi} is a $(\P^{t,a},\F)$-martingale. We solve this problem reasoning as in \cite{filovschm11}, Theorem 5.1, point (iv). Basically, for every $T\in\R_+$ we repeat the above construction with $\bar\Fc_T$ in place of $\bar\Fc$, changing what in \cite{jac75} is called the past information. More precisely, let $T\in\R_+$ and define $\hat\F^T=(\hat\Fc_s^T)_{s\geq0}$, where $\hat\Fc_s^T$ $:=$ $\bar\Fc_T\otimes\Fc_s'$. Let
\beqs
\nu^T(\omega,F) &=& \int_F 1_{\{s\leq T\}}1_{\{s<T_\infty(\omega)\}} \lambda(I_s^{t,a}(\omega),de)ds.
\enqs
Proceeding as before, we conclude that there exists a unique probability measure on $(\Omega,\bar\Fc_T\otimes\Fc_\infty')$, denoted by $\P^{t,a,T}$, whose restriction to $(\bar\Omega,\bar\Fc_T)$ coincides with the restriction of $\bar\P$ to this measurable space, and for which $\nu^T$ is the $\hat\F^T$-compensator of $\pi$, i.e.,
\beqs
\big(\nu^T((0,s\wedge T_n]\times E') - \pi((0,s\wedge T_n]\times E')\big)_{s\geq0}
\enqs
is a $(\P^{t,a,T},\hat\F^T)$-martingale, for any $E'\in\Bc(E)$ and any $n\in\N$. This implies that $\nu^T((0,T\wedge T_n]\times E') - \pi((0,T\wedge T_n]\times E')$ is $\hat\Fc_T^T$-measurable, and therefore $\Fc_T$-measurable. Notice that
\beqs
\nu^T((0,s\wedge T_n]\times E') &=& \nu((0,s\wedge T\wedge T_n]\times E'),
\enqs
hence $\nu((0,T\wedge T_n]\times E') - \pi((0,T\wedge T_n]\times E')$ is $\Fc_T$-measurable. As $T\in\R_+$ was arbitrary, we see that the process \eqref{E:nu-pi} is $\F$-adapted. Since \eqref{E:nu-pi} is a $(\P^{t,a},\hat\F)$-martingale, with $\Fc_s\subset\hat\Fc_s$, then it is also a $(\P^{t,a},\F)$-martingale. In other words, $\nu$ is the $\F$-compensator of $\pi$.
\ep

\begin{Remark}
{\rm
Notice that, under assumption \textbf{(HFC)} and if in addition $\lambda$ satisfies the integrability condition (which implies the integrability condition \textbf{(HFC)}(iii)):
\beq
\label{E:HFC(iii)'}
\sup_{|a|\leq m} \int_E \lambda(a,de) & < & \infty, \qquad \forall\,m\in\N,
\enq
then $T_\infty=\infty$, $\P^{t,a}$ a.s., and the compensator $\nu$ is given by
\beqs
\nu(\omega,F) &=& \int_F \lambda(I_s^{t,a}(\omega),de)ds
\enqs
for any $F\in\Bc(\R_+)\otimes\Bc(E)$ and for $\P^{t,a}$ almost every $\omega\in\Omega$. Indeed, we have (we denote by $\E^{t,a}$ the expectation with respect to $\P^{t,a}$)
\begin{align*}
\E^{t,a}\bigg[\sum_{n\in\N} 1_{\{T_n<\infty\}}\bigg] \ = \ \E^{t,a}\big[\pi(\R_+\times E)\big] \ &= \ \E^{t,a}\bigg[ \int_0^\infty\int_E \pi(ds,de)\bigg] \\
&= \ \E^{t,a}\bigg[\int_0^\infty\int_E \nu(ds,de)\bigg].
\end{align*}
Therefore, for $m\in\N$ large enough,
\beqs
\E^{t,a}\bigg[\sum_{n\in\N} 1_{\{T_n<\infty\}}\bigg] \ = \ \E^{t,a}\bigg[\int_0^\infty\int_E 1_{\{s<T_\infty\}} \lambda(I_s^{t,a},de)ds\bigg] \ \leq \ T\sup_{|a'|\leq m} \int_E  \lambda(a',de) \ < \ \infty,
\enqs
where we used condition \eqref{E:HFC(iii)'} and the fact that $\P^{t,a}$ almost every path of the process $I^{t,a}$ belongs to the compact set $\{h(a)\}\cup A$. Hence, $\P^{t,a}$ a.s.,
\beqs
\sum_{n\in\N} 1_{\{T_n<\infty\}} < \infty
\enqs
which means that $T_\infty=\infty$, $\P^{t,a}$ almost surely.
\ep
}
\end{Remark}

\begin{Lemma}
\label{L:ExistenceXI}
Under assumption \textup{\textbf{(HFC)}}, for every $(t,x,a)\in[0,T]\times\R^d\times\R^q$ there exists a unique $($up to indistinguishability$)$ process $X^{t,x,a}$ $=$ $\{X_s^{t,x,a},\,s\geq0\}$ on $(\Omega,\Fc,\F,\P^{t,a})$, solution to \eqref{FSDEX} on $[t,T]$, with $X_s^{t,x,a}$ $=$ $x$ for $s<t$ and $X_s^{t,x,a}$ $=$ $X_T^{t,x,a}$ for $s>T$. Moreover, for any $(t,x,a)\in[0,T]\times\R^d\times\R^q$ there exists a positive constant $C_a$ such that
\begin{equation}
\label{EstimateXI}
\E^{t,a}\Big[\sup_{t\leq s\leq T} \big( |X_s^{t,x,a}|^2 + |I_s^{t,a}|^2 \big) \Big] \ \leq \ C_a \big(1 + |x|^2+|h(a)|^2\big),
\end{equation}
where $C_a$ depends only on $T$, $|b(0,0)|$, $|\sigma(0,0)|$, the Lipschitz constants of $b$ and $\sigma$, and on the variable $a$ through the term $\sup_{t\leq s\leq T}\int_E (1\wedge|e|^2)\lambda(I_s^{t,a},de)<\infty$.
\end{Lemma}
\textbf{Proof.}
Since hypotheses (14.15) and (14.22) in \cite{jac79} are satisfied under \textbf{(HFC)}, the thesis follows from Theorem 14.23 in \cite{jac79}. Concerning estimate \eqref{EstimateXI}, taking the square in \eqref{FSDEX} (using the standard inequality $(x_1+\cdots+x_4)^2\leq 4(x_1^2+\cdots+x_4^2)$, for any $x_1,\ldots,x_4\in\R$) and then the supremum, we find
\begin{align}
\label{E:ProofWeak1}
\sup_{t\leq u\leq s}|X_u^{t,x,a}|^2 \ &\leq \ 4|x|^2 + 4\sup_{t\leq u\leq s}\bigg|\int_t^u b(X_r^{t,x,a},I_r^{t,a}) dr\bigg|^2 + 4\sup_{t\leq u\leq s}\bigg|\int_t^u \sigma(X_r^{t,x,a},I_r^{t,a}) dW_r\bigg|^2 \notag \\
&\quad \ + 4\sup_{t\leq u\leq s}\bigg|\int_t^u\int_E \beta(X_{r^-}^{t,x,a},I_r^{t,a},e) \tilde\pi(dr,de)\bigg|^2.
\end{align}
Notice that, from Cauchy-Schwarz inequality we have
\begin{equation}
\label{E:ProofWeak2}
\E^{t,a}\bigg[\sup_{t\leq u\leq s}\bigg|\int_t^u b(X_r^{t,x,a},I_r^{t,a}) dr\bigg|^2\bigg] \ \leq \ T\,\E^{t,a}\bigg[\int_t^s \big|b(X_r^{t,x,a},I_r^{t,a})\big|^2 dr\bigg].
\end{equation}
Moreover, from Burkholder-Davis-Gundy inequality there exists a positive constant $\bar C$ such that
\begin{equation}
\label{E:ProofWeak3}
\E^{t,a}\bigg[\sup_{t\leq u\leq s}\bigg|\int_t^u \sigma(X_r^{t,x,a},I_r^{t,a}) dW_r\bigg|^2\bigg] \ \leq \ \bar C\,\E^{t,a}\bigg[\int_t^s \text{tr}\big(\sigma\sigma\trans(X_r^{t,x,a},I_r^{t,a})\big) dr\bigg].
\end{equation}
Similarly, since the local martingale $M_u=\int_t^u\int_E \beta(X_{r^-}^{t,x,a},I_r^{t,a},e) \tilde\pi(dr,de)$, $t\leq u\leq s$, is such that $[M]_u = \int_t^u\int_E |\beta(X_{r^-}^{t,x,a},I_r^{t,a},e)|^2 \pi(dr,de)$, from Burkholder-Davis-Gundy inequality we obtain
\begin{align}
\label{E:ProofWeak4}
&\E^{t,a}\bigg[\sup_{t\leq u\leq s}\bigg|\int_t^u\int_E \beta(X_{r^-}^{t,x,a},I_r^{t,a},e) \tilde\pi(dr,de)\bigg|^2\bigg] \notag \\
&\leq \ \bar C\,\E^{t,a}\bigg[\int_t^s\int_E \big|\beta(X_{r^-}^{t,x,a},I_r^{t,a},e)\big|^2 \pi(dr,de)\bigg] \notag \\
&= \ \bar C\,\E^{t,a}\bigg[\int_t^s\int_E \big|\beta(X_{r^-}^{t,x,a},I_r^{t,a},e)\big|^2 \lambda(I_r^{t,a},de)dr\bigg].
\end{align}
In conclusion, taking the expectation in \eqref{E:ProofWeak1} and using \eqref{E:ProofWeak2}-\eqref{E:ProofWeak3}-\eqref{E:ProofWeak4}, we find (denoting $C_a$ a generic positive constant depending only on $T$, $|b(0,0)|$, $|\sigma(0,0)|$, the Lipschitz constants of $b$ and $\sigma$, and on the variable $a$ through the term $\sup_{t\leq s\leq T}\int_E (1\wedge|e|^2)\lambda(I_s^{t,a},de)<\infty$)
\[
\E^{t,a}\Big[\sup_{t\leq u\leq s}|X_u^{t,x,a}|^2\Big] \leq 4|x|^2 + C_a\bigg(1+\E^{t,a}\Big[\sup_{t\leq s\leq T}|I_s^{t,a}|^2\Big] + \int_t^s\E^{t,a}\Big[\sup_{t\leq u\leq r}|X_u^{t,x,a}|^2\Big]dr\bigg).
\]
Since the paths of $(I_s^{t,a})_{s\geq0}$ belong to the compact set $\{h(a)\}\cup A$, we have (here the constant $C_a$ can be chosen independent of $a$)
\[
\E^{t,a}\Big[\sup_{t\leq s\leq T}|I_s^{t,a}|^2\Big] \ \leq \ C_a\big(1+|h(a)|^2\big).
\]
Then, applying Gronwall's lemma to the map $r\mapsto\E^{t,a}[\sup_{t\leq u\leq r}|X_u^{t,x,a}|^2]$, we end up with estimate \eqref{EstimateXI}.
\ep

\section{BSDE with jumps and partially constrained diffusive part}
\label{S:BSDE}

\setcounter{equation}{0} \setcounter{Assumption}{0}
\setcounter{Theorem}{0} \setcounter{Proposition}{0}
\setcounter{Corollary}{0} \setcounter{Lemma}{0}
\setcounter{Definition}{0} \setcounter{Remark}{0}

Our aim is to derive a probabilistic representation formula, also called nonlinear Feynman-Kac formula, for the following nonlinear IPDE of HJB type:
\beq
-\frac{\partial u}{\partial t}(t,x) - \sup_{a \in A} \big( \Lc^a u(t,x) + f(x,a) \big) &=& 0, \quad\;\;\;\, (t,x)\in[0,T) \times \R^d,\label{HJB}\\
u(T,x) &=& g(x), \quad x\in\R^d, \label{condterminale}
\enq
where
\beqs
\Lc^a u(t,x) &=&b(x,a).D_x u(t,x) + \frac{1}{2}\text{tr}\big(\sigma\sigma\trans (x,a)D_x^2 u(t,x)\big) \\
& & + \int_E\big(u(t,x+\beta(x,a,e))-u(t,x)-\beta(x,a,e).D_x u(t,x)\big)\lambda(a,de),
\enqs
for all $(t,x,a) \in [0,T] \times \R^d \times \R^q$. Let us firstly introduce some additional notation. Fix a finite time horizon $T$ $<$ $\infty$ and set $\Pc_T$ the $\sigma$-field of $\F$-predictable subsets of $[0,T]\times\Omega$. For any $(t,a)\in[0,T]\times\R^q$, we denote:
 \begin{itemize}
 \item ${\bf L^p_{t,a}(}\Fc_s{\bf)}$, $p$ $\geq$ $1$, $s\geq0$, the set of $\Fc_s$-measurable random variables $X$ such that $\E^{t,a}[|X|^p]$ $<$ $\infty$.
 \item ${\bf S^2_{t,a}}$ the set  of real-valued c\`adl\`ag adapted processes $Y$ $=$
$(Y_s)_{t\leq s\leq T}$ such that
\[
\|Y\|_{_{{\bf S^2_{t,a}}}}^2 := \E^{t,a}\Big[ \sup_{t\leq s\leq T} |Y_s|^2 \Big] < \infty.
\]
\item ${\bf L^p_{t,a}(t,T)}$, $p$ $\geq$ $1$, the set of real-valued
adapted  processes  $(\phi_s)_{t\leq s\leq T}$ such that
\[
\|\phi\|_{_{\bf L^p_{t,a}(t,T)}}^p := \E\bigg[\int_t^T|\phi_s|^p ds\bigg] < \infty.
\]
\item ${\bf L^p_{t,a}(W)}$, $p$ $\geq$ $1$, the set of  $\R^d$-valued $\Pc_T$-measurable processes
$Z=(Z_s)_{t\leq s\leq T}$ such that
\[
\|Z\|_{_{\bf L^p_{t,a}(W)}}^p := \E\bigg[\bigg(\int_t^T |Z_s|^2 ds\bigg)^{\frac{p}{2}}\bigg] <\infty.
\]
\item ${\bf L^p_{t,a}(B)}$, $p$ $\geq$ $1$, the set of  $\R^q$-valued $\Pc_T$-measurable processes
$V=(V_s)_{t\leq s\leq T}$ such that
\[
\|V\|_{_{\bf L^p_{t,a}(B)}}^p := \E\bigg[\bigg(\int_t^T |V_s|^2 ds\bigg)^{\frac{p}{2}}\bigg] <\infty.
\]
\item ${\bf L_{t,a}^p(\tilde\pi)}$, $p$ $\geq$ $1$,  the set of
$\Pc_T\otimes\Bc(E)$-measurable maps $U\colon[t,T]\times\Omega\times E\rightarrow \R$ such that
\[
\|U\|_{_{{\bf L^p_{t,a}(\tilde\pi)}}}^p := \E\bigg[ \bigg(\int_t^T\int_E  |U_s(e)|^2 \lambda(I_s^{t,a},de)ds\bigg)^{\frac{p}{2}}\bigg] < \infty.
\]
\item ${\bf K^2_{t,a}}$  the  set of  nondecreasing predictable processes $K$ $=$ $(K_s)_{t\leq s\leq T}$ $\in$  ${\bf S^2_{t,a}}$ with $K_t$ $=$ $0$, so that
\beqs
\|K\|_{_{{\bf S^2_{t,a}}}}^2 &=&  \E\big[|K_T|^2\big].
\enqs
\end{itemize}

\begin{Remark}
{\rm
\emph{Equivalence relation in ${\bf L^p_{t,a}(\tilde\pi)}$.} When $U^1,U^2\in{\bf L^p_{t,a}(\tilde\pi)}$, with $U^1 = U^2$ we mean $\|U^1-U^2\|_{_{{\bf L^p_{t,a}(\tilde\pi)}}}=0$, i.e., $U^1 = U^2$ $ds\otimes d\P\otimes \lambda(I_s^{t,a},de)$ a.e. on $[t,T]\times\Omega\times E$, where $ds\otimes d\P\otimes \lambda(I_s^{t,a},de)$ is the measure on $([t,T]\times\Omega\times E,\Bc(t,T)\otimes\Fc\otimes\Bc(E))$ given by:
\beqs
ds\otimes d\P\otimes \lambda(I_s^{t,a},de)(F) &=& \E^{t,a}\bigg[\int_t^T\int_E  1_F(s,\omega,e) \lambda(I_s^{t,a}(\omega),de)ds\bigg],
\enqs
for all $F\in\Bc(t,T)\otimes\Fc\otimes\Bc(E)$. See also the beginning of Section 3 in \cite{conffuhrman12}.
\ep
}
\end{Remark}

The probabilistic representation formula is given in terms of the following BSDE with jumps and partially constrained diffusive part, for any $(t,x,a)\in[0,T]\times\R^d\times\R^q$, $\P^{t,a}$ a.s.,
\beq
Y_s & = & g(X_T^{t,x,a}) + \int_s^T f(X_r^{t,x,a},I_r^{t,a}) dr + K_T - K_s - \int_s^T Z_r dW_r \label{BSDEgen}  \\
& & - \int_s^T V_r dB_r - \int_s^T\int_E U_r(e) \tilde\pi(dr,de), \qquad\qquad   t \leq s \leq T \nonumber
\enq
and
\beq
V_s & = & 0 \qquad\qquad ds\otimes d\P^{t,a}\,a.e. \label{Ucons}
\enq
We look for the minimal solution $(Y,Z,V,U,K)\in{\bf S^2_{t,a}}\times{\bf L^2_{t,a}(W)}\times{\bf L^2_{t,a}(B)}\times{\bf L^2_{t,a}(\tilde\pi)}\times{\bf K^2_{t,a}}$ to \eqref{BSDEgen}-\eqref{Ucons}, in the sense that for any other solution $(\bar Y,\bar Z,\bar V,\bar U,\bar K)\in{\bf S^2_{t,a}}\times{\bf L^2_{t,a}(W)}\times{\bf L^2_{t,a}(B)}\times{\bf L^2_{t,a}(\tilde\pi)}\times{\bf K^2_{t,a}}$ to \eqref{BSDEgen}-\eqref{Ucons} we must have $Y\leq\bar Y$. We impose the following assumptions on the terminal condition $g :\R^d \rightarrow \R$ and on the generator $f\colon \R^d \times \R^q \rightarrow \R$.

\vspace{3mm}

\ni {\bf (HBC)} \quad There exists some continuity modulus $\rho$ (namely $\rho\colon[0,\infty)\rightarrow[0,\infty)$ is continuous, nondecreasing, subadditive, and $\rho(0)=0$) such that
\[
|f(x,a) - f(x',a')| + |g(x) - g(x')| \ \leq \ \rho(|x-x'| + |a-a'|),
\]
for all $x,x' \in \R^d$ and $a,a'\in \R^q$.

\begin{Proposition}
\label{Uniqueness}
Let assumptions \textup{\textbf{(HFC)}} and \textup{\textbf{(HBC)}} hold. For any $(t,x,a)\in[0,T]\times\R^d\times\R^q$, there exists at most one minimal solution on $(\Omega,\Fc,\F,\P^{t,a})$ to the BSDE \eqref{BSDEgen}-\eqref{Ucons}.
\end{Proposition}
\textbf{Proof.}
Let $(Y,Z,V,U,K)$ and $(\tilde Y,\tilde Z,\tilde V,\tilde U,\tilde K)$ be two minimal solutions to \eqref{BSDEgen}-\eqref{Ucons}. The uniqueness of the $Y$ component is clear by definition. Regarding the other components, taking the difference between the two backward equations we obtain
\[
0 \ = \ K_s - \tilde K_s - \int_t^s \big(Z_r - \tilde Z_r\big)dW_r - \int_t^s \big(V_r - \tilde V_r\big)dB_r - \int_t^s\int_E \big(U_r(e)-\tilde U_r(e)\big)\tilde\pi(dr,de),
\]
for all $t\leq s\leq T$, $\P^{t,a}$-almost surely.
Rewriting the above identity as follows
\[
\int_t^s \big(Z_r - \tilde Z_r\big)dW_r +  \int_t^s \big(V_r - \tilde V_r\big)dB_r \ = \ K_s - \tilde K_s - \int_t^s\int_E \big(U_r(e)-\tilde U_r(e)\big)\tilde\pi(dr,de),
\]
we see that the right-hand side is a finite variation process, while the left-hand side has not finite variation, unless $Z = \tilde Z$ and $V = \tilde V$. Therefore, we obtain the identity
\[
\int_t^s\int_E \big(U_r(e)-\tilde U_r(e)\big)\pi(dr,de) \ = \ \int_t^s\int_E \big(U_r(e)-\tilde U_r(e)\big)\lambda(I_r^{t,a},de)dr + K_s - \tilde K_s,
\]
where the right-hand side is a predictable process, therefore it has no totally inaccessible jumps (see, e.g., Proposition 2.24, Chapter I, in \cite{jacshiryaev03}); on the other hand, the left-hand side is a pure-jump process with totally inaccessible jumps, unless $U=\tilde U$. As a consequence, we must have $U=\tilde U$, from which it follows that $K=\tilde K$.
\ep

\vspace{3mm}

To guarantee the existence of the minimal solution to \eqref{BSDEgen}-\eqref{Ucons} we shall need the following result.

\begin{Lemma} \label{vpoly}
Let assumptions \textup{\textbf{(HFC)}} and \textup{\textbf{(HBC)}} hold. Then, for any initial condition $(t,x,a)$ $\in$ $[0,T]\times \R^d\times\R^q$, there exists a solution $\{(\bar Y_s^{t,x,a},\bar Z_s^{t,x,a},\bar V_s^{t,x,a},\bar U_s^{t,x,a},\bar K_s^{t,x,a}),\,t \leq s \leq T\}$ on $(\Omega,\Fc,\F,\P^{t,a})$ to the BSDE \eqref{BSDEgen}-\eqref{Ucons}, with $\bar Y_s^{t,x,a}$ $=$ $\bar v(s,X_s^{t,x,a})$ for some deterministic function $\bar v$ on $[0,T]\times \R^d$ satisfying a linear growth condition
\[
\sup_{(t,x)\in[0,T]\times\R^d}\frac{|\bar v(t,x)|}{1+|x|}\ < \ \infty.
\]
\end{Lemma}
\textbf{Proof.}
Let us consider the mollifier $\eta(x)=\bar c\exp(1/(|x|^2-1))1_{\{|x|<1\}}$, where $\bar c>0$ is such that $\int_{\R^d}\eta(x)dx=1$. Let us introduce the smooth function
\[
\bar v(t,x) \ = \ \bar C e^{\rho(T-t)}\bigg(1+\int_{\R^d}\eta(x-y)|y|dy\bigg), \qquad \forall\,(t,x)\in[0,T]\times\R^d\times\R^q,
\]
for some positive constants $\bar C$ and $\rho$ to be determined later. We claim that for $\bar C$ and $\rho$ large enough, the function $\bar v$ is a classical supersolution to \eqref{HJB}-\eqref{condterminale}. More precisely, $\bar C$ is such that $g(x) \leq \bar C(1+\int_{\{|y|<1\}}\eta(y)|x-y|dy)$, for all $x\in \R^d$, which follows from $\int_{\{|y|<1\}}\eta(y)|x-y|dy\geq||x|-1|$ and from the uniform continuity of $g$ (which implies the linear growth of $g$ itself). Furthermore, using the compactness of $A$, a straightforward calculation shows that
\beqs
-\frac{\partial \bar v}{\partial t}(t,x) - \sup_{a \in A} \big( \Lc^a \bar v(t,x) + f(x,a) \big) & \geq & (\rho-C)\bar v(t,x),
\enqs
for some positive constant $C$ depending only on $\bar C$ and the linear growth conditions of $b$, $\sigma$, $\beta$, and $f$. Then, we choose $\rho \geq C$. Let us now define the quintuple $(\bar Y,\bar Z,\bar U,\bar L,\bar K)$ as follows:
\begin{align*}
\bar Y_s &:= \bar v(s,X_s^{t,x,a}), \quad \text{for }t\leq s < T, \qquad \bar Y_T := g(X_T^{t,x,a}), \\
\bar Z_s &:= \sigma(X_{s^-}^{t,x,a},I_s^{t,a}) D_x\bar v(s,X_{s^-}^{t,x,a}), \quad t\leq s \leq T, \\
\bar V_s &:= 0, \quad t\leq s \leq T, \\
\bar U_s(e) &:= \bar v(s,X_{s^-}^{t,x,a} + \beta(X_{s^-}^{t,x,a},I_s^{t,a},e)) - \bar v(s,X_{s^-}^{t,x,a}), \quad t\leq s \leq T,\,e\in E, \\
\bar K_s &:= \int_t^s \Big(- \frac{\partial \bar v}{\partial t}(r,X_r^{t,x,a}) - \Lc^{I_r^{t,a}} \bar v(r,X_r^{t,x,a}) - f\big(X_r^{t,x,a},I_r^{t,a}\big)\Big)dr, \quad t\leq s < T, \\
\bar K_T &:= \bar K_{T^-} + \bar v(T,X_T^{t,x,a}) - g(X_T^{t,x,a}).
\end{align*}
We see that $(\bar Y,\bar Z,\bar V,\bar U,\bar K)$ lies in ${\bf S^2_{t,a}}\times{\bf L^2_{t,a}(W)}\times{\bf L^2_{t,a}(B)}\times{\bf L^2_{t,a}(\tilde\pi)}\times{\bf K^2_{t,a}}$. Moreover, by It\^o's formula applied to $\bar v(s,X_s^{t,x,a})$, we conclude that $(\bar Y,\bar Z,\bar V,\bar U,\bar K)$ solves \eqref{BSDEgen}, and the constraint \eqref{Ucons} is clearly satisfied.
\ep

\subsection{Existence of the minimal solution by penalization}

In this section we prove the existence of the minimal solution to \eqref{BSDEgen}-\eqref{Ucons}. We use a penalization approach and introduce the indexed sequence of BSDEs with jumps, for any $(t,x,a)\in[0,T]\times\R^d\times\R^q$, $\P^{t,a}$ a.s.,
\beq
Y_s^n &=& g(X_T^{t,x,a}) + \int_s^T f(X_r^{t,x,a},I_r^{t,a}) dr + K_T^n - K_s^n - \int_s^T  Z_r^n dW_r \nonumber \\
& & -  \int_s^T V_r^n dB_r - \int_s^T\int_E U_r^n(e)\tilde\pi(dr,de), \qquad\qquad t\leq s\leq T,  \label{penBSDE}
\enq
for $n$ $\in$ $\N$, where $K^n$ is the nondecreasing continuous process defined by
\[
K_s^n \; = \;  n \int_t^s |V_r^n| dr, \qquad t\leq s\leq T.
\]

\begin{Proposition}
\label{P:ExistenceY^n}
Under assumptions \textup{\textbf{(HFC)}} and \textup{\textbf{(HBC)}}, for every $(t,x,a)\in[0,T]\times\R^d\times\R^q$ and every $n\in\N$ there exists a unique solution $(Y^{n,t,x,a},Z^{n,t,x,a},V^{n,t,x,a},U^{n,t,x,a})$ $\in$ ${\bf S^2_{t,a}}\times{\bf L^2_{t,a}(W)}\times{\bf L^2_{t,a}(B)}\times{\bf L^2_{t,a}(\tilde\pi)}$ on $(\Omega,\Fc,\F,\P^{t,a})$ satisfying the BSDE with jumps \reff{penBSDE}.
\end{Proposition}
{\bf Proof.}
As usual, the proof is based on a fixed point argument. More precisely, let us consider the function $\Phi\colon\mathbf{L^2_{t,a}(t,T)}\times\mathbf{L^2_{t,a}(W)}\times\mathbf{L^2_{t,a}(B)}\times\mathbf{L^2_{t,a}(\tilde\pi)}\rightarrow\mathbf{L^2_{t,a}(t,T)}\times\mathbf{L^2_{t,a}(W)}\times\mathbf{L^2_{t,a}(B)}\times\mathbf{L^2_{t,a}(\tilde\pi)}$, mapping $(Y',Z',V',U')$ to $(Y,Z,V,U)$ defined by
\beq
Y_s &=& g(X_T^{t,x,a}) + \int_s^T f_n(X_r^{t,x,a},I_r^{t,a},V_r') dr - \int_s^T  Z_r dW_r \nonumber \\
& & -  \int_s^T V_r(a) dB_s - \int_s^T\int_E U_r(e)\tilde\pi(dr,de), \label{E:MartReprProofBSDE}
\enq
where
\[
f_n(x,a,v) \ = \ f(x,a) + n|v|.
\]
More precisely, the quadruple $(Y,Z,V,U)$ is constructed as follows: we consider the martingale $M_s=\E^{t,a}[g(X_T^{t,x,a})+\int_t^T f_n(X_r^{t,x,a},I_r^{t,a},V_r')dr|\Fc_s]$, which is square integrable under the assumptions on $g$ and $f$. From the martingale representation Theorem \ref{MartReprThm}, we deduce the existence and uniqueness of $(Z,V,U)\in\mathbf{L^2_{t,a}(W)}\times\mathbf{L^2_{t,a}(B)}\times\mathbf{L^2_{t,a}(\tilde\pi)}$ such that
\begin{equation}
\label{E:MartReprProofM}
M_s \ = \ M_t + \int_t^s  Z_r dW_r +  \int_t^s V_r dB_r + \int_t^s\int_E U_r(e)\tilde\pi(dr,de).
\end{equation}
We then define the process $Y$ by
\[
Y_s \ = \ \E^{t,a}\bigg[g(X_T^{t,x,a})+\int_s^T f_n(X_r^{t,x,a},I_r^{t,a},V_r')dr\bigg|\Fc_s\bigg] \ = \ M_s - \int_t^s f_n(X_r^{t,x,a},I_r^{t,a},V_r')dr.
\]
By using the representation \eqref{E:MartReprProofM} of $M$ in the previous relation, and noting that $Y_T=g(X_T^{t,x,a})$, we see that $Y$ satisfies \eqref{E:MartReprProofBSDE}. Using the conditions on $g$ and $f$, we deduce that $Y$ lies in $\mathbf{L^2_{t,a}(t,T)}$, and also in $\mathbf{S^2_{t,a}}$. Hence, $\Phi$ is a well-defined map. We then see that $(Y^{n,t,x,a},Z^{n,t,x,a},V^{n,t,x,a},U^{n,t,x,a})$ is a solution to the penalized BSDE \eqref{penBSDE} if and only if it is a fixed point of $\Phi$. To this end, for any $\alpha>0$ let us introduce the equivalent norm on $\mathbf{L^2_{t,a}(t,T)}\times\mathbf{L^2_{t,a}(W)}\times\mathbf{L^2_{t,a}(B)}\times\mathbf{L^2_{t,a}(\tilde\pi)}$:
\[
\|(Y,Z,V,U)\|_\alpha \ := \ \E^{t,a}\bigg[\int_t^T e^{\alpha(s-t)}\bigg(|Y_s|^2 + |Z_s|^2 + |V_s|^2 + \int_E |U_s(e)|^2 \lambda(I_s^{t,a},de)\bigg) ds\bigg].
\]
It can be shown, proceeding along the same lines as in the classical case (for which we refer, e.g., to Theorem 6.2.1 in \cite{pham09}), that there exists $\bar\alpha>0$ such that $\Phi$ is a contraction on $\mathbf{L^2_{t,a}(t,T)}\times\mathbf{L^2_{t,a}(W)}\times\mathbf{L^2_{t,a}(B)}\times\mathbf{L^2_{t,a}(\tilde\pi)}$ endowed with the equivalent norm $\|\cdot\|_{\bar\alpha}$. Then, the thesis follows from the Banach contraction mapping theorem.
\ep

\vspace{3mm}

We can now prove our main result of this section. Firstly, we need the following two lemmata.

\begin{Lemma}
\label{comparison}
Under assumptions \textup{\textbf{(HFC)}} and \textup{\textbf{(HBC)}}, for every $(t,x,a)\in[0,T]\times\R^d\times\R^q$ the sequence $(Y^{n,t,x,a})_n$ is nondecreasing and upper bounded by $\bar Y^{t,x,a}$, i.e., for all $n\in\N$,
\[
Y_s^{n,t,x,a} \ \leq \ Y_s^{n+1,t,x,a} \ \leq \ \bar Y_s^{t,x,a}
\]
for all $0 \leq s \leq T$, $\P^{t,a}$ almost surely.
\end{Lemma}
{\bf Proof.}
Fix $(t,x,a)\in[0,T]\times\R^d\times\R^q$ and $n\in\N$, and observe that
\[
f_n(x,a,v) \ \leq \ f_{n+1}(x,a,v),
\]
for all $(x,a,v)\in\R^d\times\R^q\times\R^q$. Then, the inequality $Y_s^{n,t,x,a}\leq Y_s^{n+1,t,x,a}$, for all $0\leq s\leq T$, $\P^{t,a}$ a.s., follows from the comparison Theorem A.1 in \cite{khapha12}. We should notice that Theorem A.1 in \cite{khapha12} is designed for BSDE with jumps driven by a Wiener process and a Poisson random measure, while in our case we have a general random measure $\pi$. Nevertheless, Theorem A.1 in \cite{khapha12} can be proved proceeding along the same lines as in \cite{khapha12} to encompass this more general case.\\
Similarly, since $\int_0^s |\bar V_r^{t,x,a}| dr=0$, it follows that $(\bar Y^{t,x,a},\bar Z^{t,x,a},\bar V^{t,x,a},\bar U^{t,x,a}, \bar K^{t,x,a})$ solves the BSDE \eqref{BSDEgen} with generator $f_n$, for any $n\in\N$, other than with generator $f$. Therefore, we can again apply the (generalized version, with the random measure $\pi$ in place of the Poisson random measure, of the) comparison Theorem A.1 in \cite{khapha12}, from which we deduce the thesis.
\ep

\begin{Lemma}
\label{bound}
Under assumptions \textup{\textbf{(HFC)}} and \textup{\textbf{(HBC)}}, there exists a positive constant $C$ such that, for all $(t,x,a)\in[0,T]\times\R^d\times\R^q$ and $n\in\N$,
\begin{align}
\label{EstimateBSDE}
&\|Y^{n,t,x,a}\|_{_{\bf S^2_{t,a}}}^2 + \|Z^{n,t,x,a}\|_{_{\bf L^2_{t,a}(W)}}^2 + \|V^{n,t,x,a}\|_{_{\bf L^2_{t,a}(B)}}^2 + \|U^{n,t,x,a}\|_{_{\bf L^2_{t,a}(\tilde\pi)}}^2 + \|K^{n,t,x,a}\|_{_{\bf S^2_{t,a}}}^2 \notag \\
&\leq \ C\bigg(\E^{t,a}\big[|g(X_T^{t,x,a})|^2\big] + \E^{t,a}\bigg[\int_t^T|f(X_s^{t,x,a},I_s^{t,a})|^2 ds\bigg] + \|\bar v(\cdot,X_\cdot^{t,x,a})\|_{_{\bf S^2_{t,a}}}^2\bigg).
\end{align}
\end{Lemma}
{\bf Proof.}
The proof is very similar to the proof of Lemma 3.3 in \cite{khapha12}, so it is not reported. We simply recall that the thesis follows applying It\^o's formula to $|Y^{n,t,x,a}_s|^2$ between $t$ and $T$, and exploiting Gronwall's lemma and Burkholder-Davis-Gundy inequality in an usual way.
\ep

\begin{Theorem}
\label{ThmExist}
Under assumptions \textup{\textbf{(HFC)}} and \textup{\textbf{(HBC)}}, for every $(t,x,a)\in[0,T]\times\R^d\times\R^q$ there exists a unique minimal solution $(Y^{t,x,a},Z^{t,x,a},V^{t,x,a},U^{t,x,a},K^{t,x,a})$ $\in$
${\bf S^2_{t,a}}\times{\bf L^2_{t,a}(W)}\times{\bf L^2_{t,a}(B)}\times{\bf L^2_{t,a}(\tilde\pi)}\times{\bf K^2_{t,a}}$ on $(\Omega,\Fc,\F,\P^{t,a})$ to the BSDE with jumps and partially constrained diffusive part \eqref{BSDEgen}-\eqref{Ucons}, where:
\begin{itemize}
\item[\textup{(i)}] $Y^{t,x,a}$ is the increasing limit of $(Y^{n,t,x,a})_n$.
\item[\textup{(ii)}] $(Z^{t,x,a},V^{t,x,a},U^{t,x,a})$ is the weak limit of $(Z^{n,t,x,a},V^{n,t,x,a},U^{n,t,x,a})_n$  in ${\bf L^2_{t,a}(W)}\times{\bf L^2_{t,a}(B)}\times{\bf L^2_{t,a}(\tilde\pi)}$.
\item[\textup{(iii)}] $K_s^{t,x,a}$ is the weak limit of $(K_s^{n,t,x,a})_n$ in ${\bf L^2_{t,a}}(\Fc_s)$, for any $t \leq s \leq T$.
\end{itemize}
\end{Theorem}
{\bf Proof.}
Let $(t,x,a)\in[0,T]\times\R^d\times\R^q$ be fixed. From Lemma \ref{comparison} it follows that $(Y^{n,t,x,a})_n$ converges increasingly to some adapted process $Y^{t,x,a}$. We see that $Y^{t,x,a}$ satisfies $\E[\sup_{t\leq s\leq T}|Y_s^{t,x,a}|^2]<\infty$ as a consequence of the uniform estimate for $(Y^{n,t,x,a})_n$ in Lemma \ref{bound} and Fatou's lemma. Moreover, by Lebesgue's dominated convergence theorem, the convergence also holds in ${\bf L^2_{t,a}(t,T)}$. Next, by the uniform estimates in Lemma \ref{bound}, the sequence $(Z^{n,t,x,a},V^{n,t,x,a},U^{n,t,x,a})_n$ is bounded in the Hilbert space ${\bf L^2_{t,a}(W)}\times{\bf L^2_{t,a}(B)}\times{\bf L^2_{t,a}(\tilde\pi)}$. Then, we can extract a subsequence which weakly converges to some $(Z^{t,x,a},V^{t,x,a},U^{t,x,a})$ in ${\bf L^2_{t,a}(W)}\times{\bf L^2_{t,a}(B)}\times{\bf L^2_{t,a}(\tilde\pi)}$. Thanks to the martingale representation Theorem \ref{MartReprThm}, for every stopping time $t \leq \tau \leq T$, the following weak convergences hold in ${\bf L^2_{t,a}(}\Fc_\tau{\bf)}$, as $n\rightarrow\infty$,
\begin{align*}
&\int_t^\tau Z^{n,t,x,a}_s dW_s \ \rightharpoonup \ \int_t^\tau Z^{t,x,a}_s dW_s, \qquad\qquad\qquad\qquad\quad \int_t^\tau V_s^{n,t,x,a} dB_s \ \rightharpoonup \ \int_t^\tau V_s^{t,x,a} dB_s, \\
&\hspace{2.5cm}\int_t^\tau\int_E U_s^{n,t,x,a}(e)\tilde\pi(ds,de) \ \rightharpoonup \ \int_t^\tau\int_E U_s^{t,x,a}(e)\tilde\pi(ds,de).
\end{align*}
Since
\beqs
K_\tau^{n,t,x,a} & = & Y^{n,t,x,a}_t - Y^{n,t,x,a}_\tau - \int_t^\tau f(X_s^{t,x,a},I_s^{t,a}) ds + \int_t^\tau Z^{n,t,x,a}_s dW_s \\
& & + \int_t^\tau V_s^{n,t,x,a} dB_s + \int_t^\tau \int_E U^{n,t,x,a}(e)\tilde\pi(ds,de),
\enqs
we also have the following weak convergence in ${\bf L^2_{t,a}(}\Fc_\tau{\bf)}$, as $n\rightarrow\infty$,
\beqs
K_\tau^{n,t,x,a} \ \rightharpoonup \ K_\tau^{t,x,a} & := & Y^{t,x,a}_t - Y^{t,x,a}_\tau - \int_t^\tau f(X_s^{t,x,a},I_s^{t,a}) ds \\
& & + \int_t^\tau Z^{t,x,a}_s dW_s + \int_t^\tau V_s^{t,x,a} dB_s + \int_t^\tau \int_E U^{t,x,a}(e)\tilde\pi(ds,de).
\enqs
Since the process $(K_s^{n,t,x,a})_{t\leq s\leq T}$ is nondecreasing and predictable and $K_t^{n,t,x,a} = 0$, the limit process $K^{t,x,a}$ remains nondecreasing and predictable with $\E^{t,a}[|K^{t,x,a}_T|^2] < \infty$ and $K^{t,x,a}_t=0$. Moreover, by Lemma 2.2 in \cite{peng00}, $K^{t,x,a}$ and $Y^{t,x,a}$ are c\`adl\`ag, therefore $Y^{t,x,a}\in\mathbf{S^2_{t,a}}$ and $K^{t,x,a}\in\mathbf{K^2_{t,a}}$. In conclusion, we have
\beqs
Y_t^{t,x,a} & = & g(X_T^{t,x,a}) + \int_t^T f(X_s^{t,x,a},I_s^{t,a}) ds + K^{t,x,a}_T - K^{t,x,a}_t - \int_t^T Z^{t,x,a}_s dW_s \\
& & - \int_t^T V_s^{t,x,a} dB_s - \int_t^T \int_E U^{t,x,a}(e)\tilde\pi(ds,de).
\enqs
It remains to show that the jump constraint \eqref{Ucons} is satisfied. To this end, we consider the 
functional $F\colon\mathbf{L^2_{t,a}(B)}\rightarrow\R$ given by
\[
F(V) \ := \ \E^{t,a}\bigg[\int_t^T |V_s| ds\bigg], \qquad \forall\,V\in\mathbf{L^2_{t,a}(B)}.
\]
Notice that $F(V^{n,t,x,a})=\E^{t,a}[K_T^{n,t,x,a}]/n$, for any $n\in\N$. From estimate \eqref{EstimateBSDE}, we see that $F(V^{n,t,x,a})\rightarrow0$ as $n\rightarrow\infty$. Since $F$ is convex and strongly continuous in the strong topology of $\mathbf{L^2_{t,a}(B)}$, then $F$ is lower semicontinuous in the weak topology of $\mathbf{L^2_{t,a}(B)}$, see, e.g., Corollary 3.9 in \cite{brezis10}. Therefore, we find
\[
F(V^{t,x,a}) \ \leq \ \liminf_{n\rightarrow\infty}F(V^{n,t,x,a}) \ = \ 0,
\]
which implies the validity of the jump constraint \eqref{Ucons}. Hence, $(Y^{t,x,a},Z^{t,x,a},V^{t,x,a},U^{t,x,a},K^{t,x,a})$ is a solution to the BSDE with jumps and partially constrained diffusive part \eqref{BSDEgen}-\eqref{Ucons}. From Lemma \ref{comparison}, we also see that $Y^{t,x,a} = \lim Y^{n,t,x,a}$ is the minimal solution to \eqref{BSDEgen}-\eqref{Ucons}. Finally, the uniqueness of the solution $(Y^{t,x,a},Z^{t,x,a},V^{t,x,a},U^{t,x,a},K^{t,x,a})$ follows from Proposition \ref{Uniqueness}.
\ep

\section{Nonlinear Feynman-Kac formula}
\label{S:HJB}

\setcounter{equation}{0} \setcounter{Assumption}{0}
\setcounter{Theorem}{0} \setcounter{Proposition}{0}
\setcounter{Corollary}{0} \setcounter{Lemma}{0}
\setcounter{Definition}{0} \setcounter{Remark}{0}

We know from Theorem \ref{ThmExist} that, under \textbf{(HFC)} and \textbf{(HBC)}, there exists a unique minimal solution $(Y^{t,x,a},Z^{t,x,a},V^{t,x,a},U^{t,x,a},K^{t,x,a})$ on $(\Omega,\Fc,\F,\P^{t,a})$ to \eqref{BSDEgen}-\eqref{Ucons}. As we shall see below, this minimal solution admits the representation $Y_s^{t,x,a} = v(s,X_s^{t,x,a},I_s^{t,a})$, where $v\colon[0,T]\times\R^d\times\R^q\rightarrow\R$ is the deterministic function defined as
\beq
\label{v}
v(t,x,a) & := & Y_t^{t,x,a}, \qquad (t,x,a)\in[0,T]\times\R^d\times\R^q.
\enq
Our aim is to prove that the function $v$ given by \eqref{v} does not depend on the variable $a$ in the interior of $A$, and it is related to the fully nonlinear partial differential equation of HJB type \eqref{HJB}-\eqref{condterminale}. Notice that we do not know a priori whether the function $v$ is continuous. Therefore, we shall adopt the definition of discontinuous viscosity solution to \eqref{HJB}-\eqref{condterminale}. Firstly, we impose the following conditions on $h$ and $A$.
\vspace{3mm}

\ni {\bf (H$A$)}  \hspace{7mm}  There exists a compact set $A_h\subset\R^q$ such that $h(A_h)=A$. Moreover, the interior set $\mathring{A}_h$ of $A_h$ is connected, and $A_h$ $=$ Cl$(\mathring{A}_h)$, the closure of its interior. Furthermore, $h(\mathring{A}_h)=\mathring{A}$.

\vspace{3mm}

We also impose some conditions on $\lambda$, which will imply the validity of a comparison theorem for viscosity sub and supersolutions to the fully nonlinear IPDE of HJB type \eqref{HJB}-\eqref{condterminale} and also for penalized IPDE \eqref{eqvn}-\eqref{condtermvn}. To this end, let us define, for every $\delta>0$ and $(t,x,a)\in[0,T]\times\R^d\times\R^q$,
\[
I_a^{1,\delta}(t,x,\varphi) \ = \ \int_{E\cap\{|e|\leq\delta\}}\big(\varphi(t,x+\beta(x,a,e))-\varphi(t,x)-\beta(x,a,e).D_x \varphi(t,x)\big)\lambda(a,de),
\]
for any $\varphi\in C^{1,2}([0,T]\times\R^d)$, and
\[
I_a^{2,\delta}(t,x,q,u) \ = \ \int_{E\cap\{|e|>\delta\}}\big(u(t,x+\beta(x,a,e))-u(t,x)-\beta(x,a,e).q\big)\lambda(a,de),
\]
for any $q\in\R^d$ and any locally bounded function $u$. Let us impose the following continuity conditions on $I_a^{1,\delta}$ and $I_a^{2,\delta}$. Notice that, whenever $I_a^{1,\delta}$ and $I_a^{2,\delta}$ do not depend on $a$, then \textbf{(H$\lambda$)}(i)-(ii) are consequences of Lebesgue's dominated convergence theorem, while \textbf{(H$\lambda$)}(iii) follows from Fatou's lemma.

\vspace{3mm}

\ni {\bf (H$\lambda$)}
\begin{itemize}
\item[(i)] Let $\eps>0$ and define $\varphi_\eps(e) = 1\wedge|e|^2\wedge\eps$, $e\in E$. Then
\[
\sup_{a\in A} I_a^{1,\delta}(t,x,\varphi_\eps) \ \overset{\eps\rightarrow0^+}{\longrightarrow} \ 0,
\]
for any $(t,x)\in[0,T]\times\R^d$ and $\delta>0$.
\item[(ii)] Let $\varphi\in C^{1,2}([0,T]\times\R^d)$. If $(t_k,x_k,a_k) \rightarrow (t^*,x^*,a^*)$ as $k$ goes to infinity, then
\[
\lim_{k\rightarrow\infty} I_{a_k}^{1,\delta}(t_k,x_k,\varphi) \ = \ I_{a^*}^{1,\delta}(t^*,x^*,\varphi),
\]
for any $\delta>0$.
\item[(iii)] Let $u\colon[0,T]\times\R^d\rightarrow\R$ be usc (resp. lsc) and locally bounded. If $(t_k,x_k,q_k,a_k) \rightarrow (t^*,x^*,q^*,a^*)$ and $u(t_k,x_k) \rightarrow u(t^*,x^*)$, as $k$ goes to infinity, then 
\begin{align*}
\limsup_{k\rightarrow\infty} I_{a_k}^{2,\delta}(t_k,x_k,q_k,u) \ &\leq \ I_{a^*}^{2,\delta}(t^*,x^*,q^*,u) \\
\Big(\text{resp. }\liminf_{k\rightarrow\infty} I_{a_k}^{2,\delta}(t_k,x_k,q_k,u) \ &\geq \ I_{a^*}^{2,\delta}(t^*,x^*,q^*,u)\Big)
\end{align*}
for any $\delta>0$.
\end{itemize}

\vspace{3mm}

For a locally bounded function $u$ on $[0,T)\times \R^k$, we define its lower semicontinuous (lsc for short) envelope $u_*$, and upper semicontinuous (usc for short) envelope $u^*$, by
\beqs
u_*(t,\xi) = \liminf_{\substack{(s,\eta)\rightarrow (t,\xi) \\ s < T}} u(s,\xi) \quad \text{ and } \quad u^*(t,\xi) = \limsup_{\substack{(s,\eta)\rightarrow (t,\xi) \\ s < T}} u(s,\xi)
\enqs
for all $(t,\xi) \in [0,T]\times \R^k$.

\begin{Definition}$($\!Viscosity solution to \reff{HJB}-\reff{condterminale}$)$
\label{D:viscosity}
\quad
\begin{itemize}
\item[\textup{(i)}] A lsc $($resp. usc$)$ function $u$ on $[0,T]\times\R^d$ is called a \textbf{viscosity supersolution} $($resp. \textbf{viscosity subsolution}$)$ to \reff{HJB}-\reff{condterminale} if
\beqs
u(T,x) \ \geq \ (resp. \; \leq) \  g(x)
\enqs
for any $x\in \R^d$, and
\beqs
-\frac{\partial \varphi}{\partial t}(t,x) - \sup_{a \in A} \big( \Lc^a \varphi(t,x) + f(x,a) \big) \ \geq \ (resp. \; \leq) \  0
\enqs
for any $(t, x) \in [0,T)\times \R^d$ and any $\varphi \in C^{1,2}([0,T]\times \R^d)$ such that
\beqs
(u - \varphi)(t,x) \  = \ \min_{[0,T]\times \R^d} (u - \varphi) \quad (resp. \ \max_{[0,T]\times \R^d} (u - \varphi)).
\enqs
\item[\textup{(ii)}] A locally bounded function $u$ on $[0,T)\times\R^d$ is called a \textbf{viscosity solution} to \reff{HJB}-\reff{condterminale} if $u_*$ is a viscosity supersolution and $u^*$ is a viscosity subsolution to \reff{HJB}-\reff{condterminale}.
\end{itemize}
\end{Definition}

We can now state the main result of this paper.

\begin{Theorem}
\label{ThmMain}
Assume that conditions \textup{\textbf{(HFC)}}, \textup{\textbf{(HBC)}}, \textup{\textbf{(H$A$)}}, and \textup{\textbf{(H$\lambda$)}} hold. Then, the function $v$ in \eqref{v} does not depend on the variable $a$ on $[0,T)\times\R^d\times\mathring{A}$$:$
\beqs
v(t,x,a) & = & v(t,x,a'), \qquad \forall\,a,a'\in\mathring{A},
\enqs
for all $(t,x)\in[0,T)\times\R^d$. Let us then define by misuse of notation the function $v$ on $[0,T)\times\R^d$ by
\beqs
v(t,x) & = & v(t,x,a), \qquad (t,x)\in[0,T)\times\R^d,
\enqs
for any $a\in\mathring{A}$. Then $v$ is a viscosity solution to \eqref{HJB}-\eqref{condterminale}.
\end{Theorem}

The rest of the paper is devoted to the proof of Theorem \ref{ThmMain}.

\subsection{Viscosity property of the penalized BSDE}

For every $n\in\N$, let us introduce the deterministic function $v_n$ defined on $[0,T]\times \R^d \times \R^q$ by
\beq
\label{defvn}
v_n(t,x,a) \ := \ Y^{n,t,x,a}_t, \qquad (t,x,a) \in [0,T]\times \R^d \times \R^q,
\enq
where $(Y^{n,t,x,a},Z^{n,t,x,a},V^{n,t,x,a},U^{n,t,x,a})$ is the unique solution to the BSDE with jumps \reff{penBSDE}, see Proposition \ref{P:ExistenceY^n}. As we shall see in Proposition \ref{P:vn_cont}, the identification $Y_s^{n,t,x,a}=v_n(s,X_s^{t,x,a},I_s^{t,a})$ holds. Therefore, sending $n$ to infinity, it follows from the convergence results of the penalized BSDE, Theorem \ref{ThmExist}, that the minimal solution to the BSDE with jumps and partially constrained diffusive part \eqref{BSDEgen}-\eqref{Ucons} can be written as $Y_s^{t,x,a}$ $=$ $v(s,X_s^{t,x,a},I_s^{t,a})$, $t\leq s \leq T$, where $v$ is the deterministic function defined in \reff{v}.

Now, notice that, from the uniform estimate \reff{EstimateBSDE}, the linear growth conditions of $g$, $f$, and $\bar v$, and estimate \eqref{EstimateXI}, it follows that $v_n$, and thus also $v$ by passing to the limit, satisfies the following linear growth condition: there exists some positive constant $C_v$ such that, for all $n\in\N$,
\beq
\label{E:vn_v_bdd}
|v_n(t,x,a)| + |v(t,x,a)| \ \leq\  C_v\big(1+|x|+|h(a)|\big), \qquad \forall\,(t,x,a) \in [0,T]\times \R^d \times \R^q.
\enq
As expected, for every $n\in\N$, the function $v_n$ in \eqref{defvn} is related to a parabolic semi-linear penalized IPDE. More precisely, let us introduce the function $v_n^h\colon[0,T]\times\R^d\times\R^q\rightarrow\R$ given by
\begin{equation}
\label{v_n^h}
v_n^h(t,x,a) \ := \ v_n(t,x,h(a)), \qquad (t,x,a)\in[0,T]\times\R^d\times\R^q.
\end{equation}
Then, the function $v_n^h$ is related to the semi-linear penalized IPDE:
\beq
-\frac{\partial v_n^h}{\partial t}(t,x,a) - \Lc^{h(a)} v_n^h(t,x,a) - f(x,h(a)) & &\label{eqvn}\\
- \frac{1}{2}\text{tr}\big(D_a^2 v_n^h(t,x,a)\big) - n \big|D_a v_n^h(t,x,a)\big| &=& 0, \;\;  \text{on }[0,T)\times \R^d \times \R^q,\nonumber\\
v_n^h(T,\cdot,\cdot) & = & g, \;\; \text{on }\R^d\times\R^q. \label{condtermvn}
\enq
Let us provide the definition of discontinuous viscosity solution to equation \eqref{eqvn}-\eqref{condtermvn}.

\begin{Definition}$($\!Viscosity solution to \eqref{eqvn}-\eqref{condtermvn}$)$
\label{D:viscosity_n}
\quad
\begin{itemize}
\item[\textup{(i)}] A lsc $($resp. usc$)$ function $u$ on $[0,T]\times\R^d\times\R^q$ is called a \textbf{viscosity supersolution} $($resp. \textbf{viscosity subsolution}$)$ to \eqref{eqvn}-\eqref{condtermvn} if
\beqs
u(T,x,a) \ \geq \ (resp. \; \leq) \  g(x)
\enqs
for any $(x,a)\in \R^d\times\R^q$, and
\begin{align*}
-\frac{\partial \varphi}{\partial t}(t,x,a) - \Lc^{h(a)} \varphi(t,x,a) - f(x,h(a)) & \\
- \frac{1}{2}\textup{tr}\big(D_a^2 \varphi(t,x,a)\big) - n \big|D_a \varphi(t,x,a)\big|& \ \geq \ 0 \quad \big(resp. \quad \leq \ 0\big)
\end{align*}
for any $(t,x,a) \in [0,T)\times \R^d\times\R^q$ and any $\varphi \in C^{1,2}([0,T]\times (\R^d\times\R^q))$ such that
\begin{equation}
\label{vn-phi}
(u - \varphi)(t,x,a) \ = \ \min_{[0,T]\times \R^d\times\R^q} (u - \varphi) \qquad (resp. \ \max_{[0,T]\times \R^d\times\R^q} (u - \varphi)).
\end{equation}
\item[\textup{(ii)}] A locally bounded function $u$ on $[0,T)\times\R^d\times\R^q$ is called a \textbf{viscosity solution} to \eqref{eqvn}-\eqref{condtermvn} if $u_*$ is a viscosity supersolution and $u^*$ is a viscosity subsolution to \eqref{eqvn}-\eqref{condtermvn}.
\end{itemize}
\end{Definition}

Then, we have the following result, which states that the penalized BSDE with jumps \eqref{penBSDE} provides a viscosity solution to the penalized IPDE \reff{eqvn}-\reff{condtermvn}.

\begin{Proposition}
\label{P:vn_cont}
Let assumptions \textup{\textbf{(HFC)}}, \textup{\textbf{(HBC)}}, \textup{{\bf (H$A$)}}, and \textup{\textbf{(H$\lambda$)}} hold. Then, the function $v_n^h$ in \eqref{v_n^h} is a viscosity solution to \reff{eqvn}-\reff{condtermvn}. Moreover, $v_n^h$ is continuous on $[0,T]\times\R^d\times\R^q$.
\end{Proposition}
\textbf{Proof} We divide the proof into three steps.\\
\emph{Step 1. Identification $Y_s^{n,t,x,a}=v_n(s,X_s^{t,x,a},I_s^{t,a})=v_n^h(s,X_s^{t,x,a},a+B_s-B_t)$.} Inspired by the proof of Theorem 4.1 in \cite{elkaroui_peng_quenez}, we shall prove the identification $Y_s^{n,t,x,a}=v_n(s,X_s^{t,x,a},I_s^{t,a})$ using the Markovian property of $(X,I)$ studied in Appendix \ref{AppSubS:Markov} and the construction of $(Y^{n,t,x,a},Z^{n,t,x,a},U^{n,t,x,a},L^{n,t,x,a})$ based on Proposition \ref{P:ExistenceY^n}. More precisely, for any $(t,x,a)\in[0,T]\times\R^d\times\R^q$, from Proposition \ref{P:ExistenceY^n} we know that there exists a sequence $(Y^{n,k,t,x,a},Z^{n,k,t,x,a},V^{n,k,t,x,a},U^{n,k,t,x,a})\in\mathbf{L^2_{t,a}(t,T)}\times\mathbf{L^2_{t,a}(W)}\times\mathbf{L^2_{t,a}(B)}\times\mathbf{L^2_{t,a}(\tilde\pi)}$, converging to $(Y^{n,t,x,a},Z^{n,t,x,a},V^{n,t,x,a},U^{n,t,x,a})$ in $\mathbf{L^2_{t,a}(t,T)}\times\mathbf{L^2_{t,a}(W)}\times\mathbf{L^2_{t,a}(B)}\times\mathbf{L^2_{t,a}(\tilde\pi)}$, such that $(Y^{n,0,t,x,a},Z^{n,0,t,x,a},V^{n,0,t,x,a},U^{n,0,t,x,a})\equiv(0,0,0,0)$ and
\begin{align*}
Y_s^{n,k+1,t,x,a} \ &= \ g(X_T^{t,x,a}) + \int_s^T f(X_r^{t,x,a},I_r^{t,a}) dr - \int_s^T\int_E U_r^{n,k+1,t,x,a}(e)\tilde\pi(dr,de) \\
&\quad \ - \int_s^T Z_r^{n,k+1,t,x,a} dW_r - \int_s^T V_r^{n,k+1,t,x,a} dB_r + n\int_s^T \big|V_r^{n,k,t,x,a}\big| dr,
\end{align*}
for all $t\leq s\leq T$, $\P^{t,a}$ almost surely. Let us define $v_{n,k}(t,x,a):=Y_t^{n,k,t,x,a}$. We begin noting that, for $k=1$ we have
\[
Y_s^{n,1,t,x,a} \ = \ \E^{t,a}\bigg[g(X_T^{t,x,a}) + \int_s^T f(X_r^{t,x,a},I_r^{t,a}) dr\bigg|\Fc_s\bigg].
\]
Then, we see from Proposition \ref{P:MarkovFamily} that $Y_s^{n,1,t,x,a}=v_{n,1}(s,X_s^{t,x,a},I_s^{t,a})$, $d\P^{t,a}\otimes ds$-almost everywhere. Proceeding as in Lemma 4.1 of \cite{elkaroui_peng_quenez} (in particular, relying on Theorem 6.27 in \cite{cinlar_jacod_protter_sharpe80}), we also deduce that there exists a Borel measurable function $\tilde v_{n,1}$ such that $V_s^{n,1,t,x,a}=\tilde v_{n,1}(s,X_{s^-}^{t,x,a},I_s^{t,a})$, $d\P^{t,a}\otimes ds$ almost everywhere. Since $V^{n,1,t,x,a}\in\mathbf{L^2_{t,a}(B)}$, we notice that
\begin{equation}
\label{E:IntCond_tilde_v}
\E^{t,a}\bigg[\int_t^T |\tilde{v}_{n,1}(s,X_{s^-}^{t,x,a},I_s^{t,a})|^2 ds\bigg] \ < \ \infty.
\end{equation}
Let us now prove the inductive step: let $k\geq1$ be an integer and suppose that $Y_s^{n,k,t,x,a}=v_{n,k}(s,X_s^{t,x,a},I_s^{t,a})$ and $V_s^{n,k,t,x,a}=\tilde v_{n,k}(s,X_{s^-}^{t,x,a},I_s^{t,a})$, $d\P^{t,a}\otimes ds$-almost everywhere, with $\E^{t,a}[\int_t^T |\tilde{v}_{n,k}(s,X_{s^-}^{t,x,a},I_s^{t,a})|^2 ds]<\infty$. Then, we have
\begin{align*}
Y_s^{n,k+1,t,x,a} \ &= \ \E^{t,a}\bigg[g(X_T^{t,x,a}) + \int_s^T f(X_r^{t,x,a},I_r^{t,a}) dr + n\int_s^T \big|\tilde v_{n,k}(r,X_{r^-}^{t,x,a},I_r^{t,a})\big| dr\bigg|\Fc_s\bigg].
\end{align*}
Using again Proposition \ref{P:MarkovFamily} (notice that, by a monotone class argument, we can extend Proposition \ref{P:MarkovFamily} to Borel measurable functions verifying an integrability condition of the type \eqref{E:IntCond_tilde_v}) we see that $Y_s^{n,k+1,t,x,a}=v_{n,k+1}(s,X_s^{t,x,a},I_s^{t,a})$, $d\P^{t,a}\otimes ds$ almost everywhere. Now, we notice that it can be shown that $\E[\sup_{t\leq s\leq T}|Y_s^{n,k,t,x,a}-Y_s^{n,t,x,a}|]\rightarrow0$, as $k$ tends to infinity (e.g., proceeding as in Remark (b) after Proposition 2.1 in \cite{elkaroui_peng_quenez}). Therefore, $v_{n,k}(t,x,a)\rightarrow v_n(t,x,a)$ as $k$ tends to infinity, for all $(t,x,a)\in[0,T]\times\R^d\times\R^q$, from which it follows the validity of the identification $Y_s^{n,t,x,a}=v_n(s,X_s^{t,x,a},I_s^{t,a})=v_n^h(s,X_s^{t,x,a},a+B_s-B_t)$, $d\P^{t,a}\otimes ds$ almost everywhere.\\
\emph{Step 2. Viscosity property of $v_n^h$.} We shall divide the proof into two substeps.\\
\emph{Step 2a. $v_n^h$ is a viscosity solution to \eqref{eqvn}.} We now prove the viscosity supersolution property of $v_n^h$ to \eqref{eqvn}. A similar argument would show that $v_n^h$ it is a viscosity subsolution to \reff{eqvn}. Let $(\bar t,\bar x,\bar a)\in[0,T)\times\R^d\times\R^q$ and $\varphi\in C^{1,2}([0,T]\times(\R^d\times\R^q))$ such that
\begin{align}
\label{E:max_vn-psi}
0 \ = \ ((v_n^h)_*-\varphi)(\bar t,\bar x,\bar a) \ = \ \min_{[0,T]\times\R^d\times\R^q}((v_n^h)_*-\varphi).
\end{align}
Let us proceed by contradiction, assuming that
\beqs
-\frac{\partial \varphi}{\partial t}(\bar t,\bar x,\bar a) - \Lc^{h(\bar a)} \varphi(\bar t,\bar x,\bar a) - f(\bar x,h(\bar a)) & & \\
- \frac{1}{2}\textup{tr}\big(D_a^2 \varphi(\bar t,\bar x,\bar a)\big) - n \big|D_a \varphi(\bar t,\bar x,\bar a)\big| & =: & -2\eps \ < \ 0.
\enqs
Using the continuity of $b$, $\sigma$, $\beta$, $f$, and $h$, we find $\delta>0$ such that
\beq
\label{E:eps/2}
-\frac{\partial \varphi}{\partial t}(t,x,a) - \Lc^{h(a)} \varphi(t,x,a) - f(x,h(a)) & & \notag \\
- \frac{1}{2}\textup{tr}\big(D_a^2 \varphi(t,x,a)\big) - n \big|D_a \varphi(t,x,a)\big| & =: & -2\eps \ < \ 0.
\enq
for any $(t,x,a)\in[0,T]\times\R^d\times\R^q$ with $|t-\bar t|,|x-\bar x|,|a-\bar a|<\delta$. Define
\[
\tau := \inf\big\{s\geq\bar t \colon |X_s^{\bar t,\bar x,\bar a}-\bar x| > \delta,\, |B_s-B_{\bar t}| > \delta\big\}\wedge(\bar t+\delta)\wedge T.
\]
Since $X^{\bar t,\bar x,\bar a}$ is c\`{a}dl\`{a}g, it is in particular right-continuous at time $\bar t$. Therefore, $\tau>\bar t$, $\P^{\bar t,\bar a}$ almost surely. Then, an application of It\^o's formula to $\varphi(s,X_s^{\bar t,\bar x,\bar a},\bar a+B_s-B_{\bar t})$ between $\bar t$ and $\tau$, using also \eqref{E:eps/2}, yields
\begin{align}
\label{E:Inequality_psi}
&\varphi(\tau,X_\tau^{\bar t,\bar x,\bar a},\bar a+B_s-B_{\bar t}) \ \geq \ \varphi(\bar t,\bar x,\bar a) - n \int_{\bar t}^\tau\big|D_a\varphi(r,X_r^{\bar t,\bar x,\bar a},\bar a+B_r-B_{\bar t})\big| dr \notag \\
&- \int_{\bar t}^\tau f(X_r^{\bar t,\bar x,\bar a},I_r^{\bar t,\bar a})dr + \eps(\tau-\bar t) + \int_{\bar t}^\tau D_x \varphi(r,X_r^{\bar t,\bar x,\bar a},\bar a+B_r-B_{\bar t})\sigma(X_r^{\bar t,\bar x,\bar a},I_r^{\bar t,\bar a})dW_r \notag \\
&+ \int_{\bar t}^\tau D_a\varphi(r,X_r^{\bar t,\bar x,\bar a},\bar a+B_r-B_{\bar t})dB_r \\
&+ \int_{\bar t}^\tau \int_E \big(\varphi(r,X_{r^-}^{\bar t,\bar x,\bar a}+\beta(X_{r^-}^{\bar t,\bar x,\bar a},I_r^{\bar t,\bar a},e),\bar a+B_r-B_{\bar t}) - \varphi(r,X_{r^-}^{\bar t,\bar x,\bar a},\bar a+B_r-B_{\bar t})\big) \tilde\pi(dr,de). \notag
\end{align}
Writing the BSDE \eqref{penBSDE} from $\bar t$ to $\tau$, using the identification $Y^{n,\bar t,\bar x,\bar a}_s$ $=$ $v_n^h(s,X^{\bar t,\bar x,\bar a}_s,\bar a+B_s-B_{\bar t})$ and the inequality $(v_n^h)_*(\bar t,\bar x,\bar a)\leq v_n^h(\bar t,\bar x,\bar a)$,  we find
\begin{align}
\label{E:BSDEvn}
&(v_n^h)_*(\bar t,\bar x,\bar a) \ \leq \ v_n^h(\tau,X_\tau^{\bar t,\bar x,\bar a},\bar a+B_\tau-B_{\bar t}) + \int_{\bar t}^\tau f(X_r^{\bar t,\bar x,\bar a},I_r^{\bar t,\bar a}) dr + n\int_{\bar t}^\tau \big|V_r^{n,\bar t,\bar x,\bar a}\big| dr \notag \\
& - \int_{\bar t}^\tau  Z_r^{n,\bar t,\bar x,\bar a} dW_r - \int_{\bar t}^\tau V_r^{n,\bar t,\bar x,\bar a} dB_r - \int_{\bar t}^\tau\int_E U_r^{n,\bar t,\bar x,\bar a}(e) \tilde\pi(dr,de).
\end{align}
Plugging \eqref{E:BSDEvn} into \eqref{E:Inequality_psi}, we obtain
\begin{align}
\label{E:psi-vn}
&\varphi(\tau,X_\tau^{\bar t,\bar x,\bar a},\bar a+B_\tau-B_{\bar t}) - v_n^h(\tau,X_\tau^{\bar t,\bar x,\bar a},\bar a+B_\tau-B_{\bar t}) \\
&\geq \ \varphi(\bar t,\bar x,\bar a) - (v_n^h)_*(\bar t,\bar x,\bar a) + \eps(\tau-\bar t) \notag \\
&+ \int_{\bar t}^\tau D_x \varphi(r,X_r^{\bar t,\bar x,\bar a},\bar a+B_r-B_{\bar t})\sigma(X_r^{\bar t,\bar x,\bar a},I_r^{\bar t,\bar a})dW_r - \int_{\bar t}^\tau  Z_r^{n,\bar t,\bar x,\bar a} dW_r \notag \\
& - n \int_{\bar t}^\tau\big|D_a\varphi(r,X_r^{\bar t,\bar x,\bar a},\bar a+B_r-B_{\bar t})\big| dr +  n\int_{\bar t}^\tau \big|V_r^{n,\bar t,\bar x,\bar a}\big| dr - \int_{\bar t}^\tau V_r^{n,\bar t,\bar x,\bar a} dB_r \notag \\
&+ \int_{\bar t}^\tau D_a\varphi(r,X_r^{\bar t,\bar x,\bar a},\bar a+B_r-B_{\bar t})dB_r- \int_{\bar t}^\tau\int_E U^{n,\bar t,\bar x,\bar a}_r(e) \tilde\pi(dr,de) \notag \\
&+ \int_{\bar t}^\tau \int_E \big(\varphi(r,X_{r^-}^{\bar t,\bar x,\bar a}+\beta(X_{r^-}^{\bar t,\bar x,\bar a},I_r^{\bar t,\bar a},e),\bar a+B_r-B_{\bar t}) - \varphi(r,X_{r^-}^{\bar t,\bar x,\bar a},\bar a+B_r-B_{\bar t})\big) \tilde\pi(dr,de). \notag
\end{align}
Let us introduce the process $\alpha\colon[\bar t,T]\times\Omega\rightarrow\R^q$ given by
\begin{align*}
\alpha_r \ &= \ n\frac{|D_a\varphi(r,X_r^{\bar t,\bar x,\bar a},\bar a+B_r-B_{\bar t})| - |V_r^{n,\bar t,\bar x,\bar a}|}{|D_a\varphi(r,X_r^{\bar t,\bar x,\bar a},\bar a+B_r-B_{\bar t}) - V_r^{n,\bar t,\bar x,\bar a}|}\cdot \\
&\quad \ \cdot \frac{D_a\varphi(r,X_r^{\bar t,\bar x,\bar a},\bar a+B_r-B_{\bar t}) - V_r^{n,\bar t,\bar x,\bar a}}{|D_a\varphi(r,X_r^{\bar t,\bar x,\bar a},\bar a+B_r-B_{\bar t}) - V_r^{n,\bar t,\bar x,\bar a}|} 1_{\{|D_a\varphi(r,X_r^{\bar t,\bar x,\bar a},\bar a+B_r-B_{\bar t}) - V_r^{n,\bar t,\bar x,\bar a}|\neq0\}}
\end{align*}
for all $\bar t\leq r\leq T$. Notice that $\alpha$ is bounded, moreover
\[
n\big(|D_a\varphi(r,X_r^{\bar t,\bar x,\bar a},\bar a+B_r-B_{\bar t})| - |V_r^{n,\bar t,\bar x,\bar a}|\big) \ = \ \big(D_a\varphi(r,X_r^{\bar t,\bar x,\bar a},\bar a+B_r-B_{\bar t}) - V_r^{n,\bar t,\bar x,\bar a}\big)\alpha_r.
\]
Consider now the probability measure $\P^{\bar t,\bar a,\alpha}$ equivalent to $\P^{\bar t,\bar a}$ on $(\Omega,\Fc_T)$, with Radon-Nikodym density given by
\[
\frac{d\P^{\bar t,\bar a,\alpha}}{d\P^{\bar t,\bar a}}\bigg|_{\Fc_s} \ = \ \Ec\bigg(\int_t^\cdot\alpha_rdB_r - \frac{1}{2}\int_t^\cdot|\alpha_r|^2dr\bigg)_{\!\!s}
\]
for all $\bar t\leq s\leq T$, where $\Ec(\cdot)$ is the Dol\'eans-Dade exponential. Notice that the stochastic integrals with respect to $W$ and $\tilde\pi$ in \eqref{E:psi-vn} remain martingales with respect to $\P^{\bar t,\bar a,\alpha}$, while the effect of the measure $\P^{\bar t,\bar a,\alpha}$ is to render the process $B_r-B_t-\int_t^r\alpha_udu$ a Brownian motion. As a consequence, taking the expectation with respect to $\P^{\bar t,\bar a,\alpha}$ in \eqref{E:psi-vn} we end up with (recalling that $\varphi(\bar t,\bar x,\bar a)=(v_n^h)_*(\bar t,\bar x,\bar a)$)
\begin{align*}
&\E^{\P^{\bar t,\bar a,\alpha}}\big[\varphi(\tau,X_\tau^{\bar t,\bar x,\bar a},\bar a+B_\tau-B_{\bar t}) - (v_n^h)_*(\tau,X_\tau^{\bar t,\bar x,\bar a},\bar a+B_\tau-B_{\bar t})\big] \\
&\geq \ \E^{\P^{\bar t,\bar a,\alpha}}\big[\varphi(\tau,X_\tau^{\bar t,\bar x,\bar a},\bar a+B_\tau-B_{\bar t}) - v_n^h(\tau,X_\tau^{\bar t,\bar x,\bar a},\bar a+B_\tau-B_{\bar t})\big] \ \geq \ \eps\E^{\P^{\bar t,\bar a,\alpha}}[\tau-\bar t].
\end{align*}
Since $\tau>\bar t$, $\P^{\bar t,\bar a}$-a.s., it follows that $\tau>\bar t$, $\P^{\bar t,\bar a,\alpha}$-a.s., therefore $\E^{\P^{\bar t,\bar a,\alpha}}[\tau-\bar t]>0$. This implies that there exists $B\in\Fc_\tau$ such that $(\varphi(\tau,X_\tau^{\bar t,\bar x,\bar a},\bar a+B_\tau-B_{\bar t}) - (v_n^h)_*(\tau,X_\tau^{\bar t,\bar x,\bar a},\bar a+B_\tau-B_{\bar t}))1_B>0$ and $\P^{\bar t,\bar a,\alpha}(B)>0$. This is a contradiction with \eqref{E:max_vn-psi}.\\
\emph{Step 2b. $v_n^h$ is a viscosity solution to \eqref{condtermvn}.} As in step 2a, we shall only prove the viscosity supersolution property of $v_n^h$ to \eqref{condtermvn}, since the viscosity subsolution of $v_n^h$ to \eqref{condtermvn} can be proved similarly. Let $(\bar x,\bar a)\in\R^d\times\R^q$. Our aim is to show that
\begin{equation}
\label{E:SuperSolTermCond}
(v_n^h)_*(T,\bar x,\bar a) \ \geq \ g(\bar x).
\end{equation}
Notice that there exists $(t_k,x_k,a_k)_k\subset[0,T)\times\R^d\times\R^q$ such that
\[
\big(t_k,x_k,a_k,v_n^h(t_k,x_k,a_k)\big) \ \overset{k\rightarrow\infty}{\longrightarrow} \ \big(\bar t,\bar x,\bar a,(v_n^h)_*(\bar t,\bar x,\bar a)\big).
\]
Recall that $v_n^h(t_k,x_k,a_k)=Y_{t_k}^{n,t_k,x_k,a_k}$ and
\begin{align}
\label{E:SuperSolTermCond2}
Y_{t_k}^{n,t_k,x_k,a_k} \ &= \ \E^{t_k,a_k}\big[g(X_T^{t_k,x_k,a_k})\big] + \int_{t_k}^T \E^{t_k,a_k}\big[f(X_s^{t_k,x_k,a_k},I_s^{t_k,a_k})\big] ds \notag \\
&\quad \ + n\int_{t_k}^T \E^{t_k,a_k}\big[\big|V_s^{n,t_k,x_k,a_k}\big|\big] ds.
\end{align}
Now we observe that, from classical convergence results of diffusion processes with jumps, see, e.g., Theorem 4.8, Chapter IX, in \cite{jacshiryaev03}, we have that the law of $(X^{t',x',a'},I^{t',a'})$ weakly converges to the law of $(X^{t,x,a},I^{t,a})$. As a consequence, we obtain
\[
\E^{t_k,a_k}\big[g(X_T^{t_k,x_k,a_k})\big] \
\overset{k\rightarrow\infty}{\longrightarrow} \ g(\bar x).
\]
Moreover, from estimate \eqref{EstimateXI} and \eqref{EstimateBSDE}, it follows by Lebesgue's dominated convergence theorem that the two integrals in time in \eqref{E:SuperSolTermCond2} go to zero as $k\rightarrow\infty$. In conclusion, letting $k\rightarrow\infty$ in \eqref{E:SuperSolTermCond2} we deduce that $(v_n^h)_*(T,\bar x,\bar a)=g(\bar x)$, therefore \eqref{E:SuperSolTermCond} holds. Notice that, from this proof, we also have that, for any $(x,a)\in\R^d\times\R^q$, $v_n^h(t',x',a')\rightarrow v_n^h(T,x,a)=g(x)$, as $(t',x',a')\rightarrow(T,x,a)$, with $t'<T$. In other words, $v_n^h$ is continuous at $T$.\\
\emph{Step 3. Continuity of $v_n^h$ on $[0,T]\times\R^d\times\R^q$.} The continuity of $v_n^h$ at $T$ was proved in step 2b. On the other hand, the continuity of $v_n^h$ on $[0,T)\times\R^d\times\R^q$ follows from the comparison theorem for viscosity solutions to equation \eqref{eqvn}-\eqref{condtermvn}. We notice, however, that a comparison theorem for equation \eqref{eqvn}-\eqref{condtermvn} does not seem to be at disposal in the literature. Indeed, Theorem 3.5 in \cite{barles97} applies to semilinear PDEs in which a L\'evy measure appears, instead in our case $\lambda$ depends on $a$. We can not even apply our comparison Theorem \ref{CompThm}, designed for equation \eqref{HJB}-\eqref{condterminale}, since in Theorem \ref{CompThm} the variable $a$ is a parameter while in equation \eqref{eqvn} is a state variable. Moreover, in \eqref{eqvn} there is also a nonlinear term in the gradient $D_av_n^h$, i.e., we need a comparison theorem for an equation with a generator $f$ depending also on $z$. Nevertheless, we observe that, under assumption \textup{\textbf{(H$\lambda$)}} we can easily extend Theorem 3.5 in \cite{barles97} to our case and, since the proof is very similar to that of Theorem 3.5 in \cite{barles97}, we do not prove it here to alleviate the presentation.
\ep

\subsection{The non dependence of the function $v$ on the variable $a$}
\label{SubS:NonDep_a}

In the present subsection, our aim is to prove that the function $v$ does not depend on the variable $a$. This is indeed a consequence of the constraint \eqref{Ucons} on the component $V$ of equation \eqref{BSDEgen}. If $v$ were smooth enough, then, for any $(t,x,a)\in[0,T]\times\R^d\times\R^q$, we could express the process $V^{t,x,a}$ as follows (we use the notations $h(a)=(h_i(a))_{i=1,\ldots,q}$, $D_a h(a)=(D_{a_j}h_i(a))_{i,j=1,\ldots,q}$, and finally $D_hv$ to denote the gradient of $v$ with respect to its last argument)
\[
V_s^{t,x,a} \ = \ D_hv(s,X_s^{t,x,a},I_s^{t,a})D_ah(a+B_s-B_t), \qquad t\leq s\leq T.
\]
Therefore, from the constraint \eqref{Ucons} we would find
\[
\E^{t,a}\bigg[\int_t^{t+\delta}|D_hv(s,X_s^{t,x,a},I_s^{t,a})D_ah(a+B_s-B_t)|ds\bigg] \ = \ 0,
\]
for any $\delta>0$. By sending $\delta$ to zero in the above equality divided by $\delta$, we would obtain
\[
|D_hv(t,x,h(a))D_ah(a)| \ = \ 0.
\]
Let us consider the function $v^h\colon[0,T]\times\R^d\times\R^q\rightarrow\R$ given by
\begin{equation}
\label{v^h}
v^h(t,x,a) \ := \ v(t,x,h(a)), \qquad (t,x,a)\in[0,T]\times\R^d\times\R^q.
\end{equation}
Then $|D_a v^h|\equiv0$, so that the function $v^h$ is constant with respect to $a$. Since $h(\R^q)=A$, we have that $v$ does not depend on the variable $a$ on $A$.

Unfortunately, we do not know if $v$ is regular enough in order to justify the above passages. Therefore, we shall rely on viscosity solutions techniques to derive the non dependence of $v$ on the variable $a$. To this end, let us introduce the following first-order PDE:
\begin{equation}
\label{PDEa}
-\,|D_av^h(t,x,a)| \ = \ 0, \qquad (t,x,a)\in[0,T)\times\R^d\times\R^q.
\end{equation}
\begin{Lemma}
\label{L:PDEa}
Let assumptions \textup{\textbf{(HFC)}}, \textup{\textbf{(HBC)}}, \textup{{\bf (H$A$)}}, and \textup{\textbf{(H$\lambda$)}} hold. The function $v^h$ in \eqref{v^h} is a viscosity supersolution to \eqref{PDEa}$:$ for any $(t,x,a)\in[0,T)\times\R^d\times\R^q$ and any function $\varphi\in C^{1,2}([0,T]\times(\R^d\times\R^q))$ such that
\[
(v^h-\varphi)(t,x,a) \ = \ \min_{[0,T]\times\R^d\times\R^q}(v^h-\varphi)
\]
we have
\[
-\,|D_av^h(t,x,a)| \ \geq \ 0.
\]
\end{Lemma}
\textbf{Proof.}
We know that $v^h$ is the pointwise limit of the nondecreasing sequence of functions $(v_n^h)_n$. By continuity of $v_n^h$, the function $v^h$ is lower semicontinuous and we have (see, e.g., page 91 in \cite{bar94}):
\[
v^h(t,x,a) \ = \ v_*^h(t,x,a) \ = \ \liminf_{n\rightarrow\infty}\!{}_*\,v_n^h(t,x,a),
\]
for all $(t,x,a)\in[0,T)\times\R^d\times\R^q$, where
\[
\liminf_{n\rightarrow\infty}\!{}_*\,v_n^h(t,x,a) \ = \ \liminf_{\substack{n\rightarrow\infty\\
(t',x',a')\rightarrow(t,x,a)\\ t'<T}}v_n^h(t',x',a'), \qquad (t,x,a)\in[0,T)\times\R^d\times\R^q.
\]
Let $(t,x,a)\in[0,T)\times\R^d\times\R^q$ and $\varphi\in C^{1,2}([0,T]\times(\R^d\times\R^q))$ such that
\[
(v^h-\varphi)(t,x,a) \ = \ \min_{[0,T]\times\R^d\times\R^q}(v^h-\varphi).
\]
We may assume, without loss of generality, that this minimum is strict. Up to a suitable negative perturbation of $\varphi$ for large values of $x$ and $a$, we can assume, without loss of generality, that there exists a bounded sequence $(t_n,x_n,a_n)\in[0,T]\times\R^d\times\R^q$ such that
\[
(v_n^h-\varphi)(t_n,x_n,a_n) \ = \ \min_{[0,T]\times\R^d\times\R^q}(v_n^h-\varphi).
\]
Then, it follows that, up to a subsequence,
\begin{equation}
\label{E:ProofPDEa}
\big(t_n,x_n,a_n,v_n^h(t_n,x_n,a_n)\big) \ \longrightarrow \ \big(t,x,a,v^h(t,x,a)\big), \qquad \text{as }n\rightarrow\infty.
\end{equation}
Now, from the viscosity supersolution property of $v_n^h$ at $(t_n,x_n,a_n)$ with the test function $\varphi$, we have
\begin{align*}
-\frac{\partial \varphi}{\partial t}(t_n,x_n,a_n) - \Lc^{h(a_n)} \varphi(t_n,x_n,a_n) - f(x_n,h(a_n))& \\
- \frac{1}{2}\textup{tr}\big(D_a^2 \varphi(t_n,x_n,a_n)\big) - n \big|D_a \varphi(t_n,x_n,a_n)\big|& \ \geq \ 0,
\end{align*}
which implies
\begin{align*}
\big|D_a \varphi(t_n,x_n,a_n)\big| \ &\leq \ \frac{1}{n}\bigg(-\frac{\partial \varphi}{\partial t}(t_n,x_n,a_n) - \Lc^{h(a_n)} \varphi(t_n,x_n,a_n) \\
&\quad \ - f(x_n,h(a_n)) - \frac{1}{2}\textup{tr}\big(D_a^2 \varphi(t_n,x_n,a_n)\big)\bigg).
\end{align*}
Sending $n$ to infinity, we get from \eqref{E:ProofPDEa} and the continuity of $b$, $\sigma$, $\beta$, $f$, and $h$:
\[
\big|D_a \varphi(t,x,a)\big| \ = \ 0,
\]
which is the thesis.
\ep

\vspace{3mm}

We can now state the main result of this subsection.
\begin{Proposition}
\label{P:NonDep_v_a}
Let assumptions \textup{\textbf{(HFC)}}, \textup{\textbf{(HBC)}}, \textup{{\bf (H$A$)}}, and \textup{\textbf{(H$\lambda$)}} hold. Then, the function $v$ in \eqref{v} does not depend on its last argument on $[0,T)\times\R^d\times\mathring{A}$:
\[
v(t,x,a) \ = \ v(t,x,a'), \qquad a,a'\in\mathring{A},
\]
for any $(t,x)\in[0,T)\times\R^d$.
\end{Proposition}
\textbf{Proof.}
From Lemma \ref{L:PDEa}, we have that $v^h$ is a viscosity supersolution to the first-order PDE:
\[
-\,\big|D_av^h(t,x,a)\big| \ = \ 0, \qquad (t,x,a)\in[0,T)\times\R^d\times\mathring{A}_h,
\]
where $A_h$ was introduced in assumption \textbf{(H$A$)}. Then, from Proposition 5.2 in \cite{khapha12} we conclude that $v^h$ does not depend on the variable $a$ in $\mathring{A}_h$:
\[
v^h(t,x,a) \ = \ v^h(t,x,a'), \qquad (t,x)\in[0,T)\times\R^d,\,a,a'\in\mathring{A}_h.
\]
Since, from assumption \textbf{(H$A$)} we have $h(\mathring{A}_h)=\mathring{A}$, we deduce the thesis.
\ep

\subsection{Viscosity properties of the function $v$}

From Proposition \ref{P:NonDep_v_a}, by misuse of notation, we can define the function $v$ on $[0,T)\times\R^d$ by
\[
v(t,x) \ = \ v(t,x,a), \qquad (t,x)\in[0,T)\times\R^d,
\]
for some $a\in\mathring{A}$. Since $h(\mathring{A}_h)=\mathring{A}$, we also have
\[
v(t,x) \ = \ v^h(t,x,a), \qquad (t,x,a)\in[0,T)\times\R^d,
\]
for some $a\in\mathring{A}_h$. Moreover, from estimate \eqref{E:vn_v_bdd} we deduce the linear growth condition for $v$ (recall that $h(a)\in A$ and $A$ is a compact set, so that $h$ is a bounded function):
\begin{equation}
\label{E:LinGrowth_v}
\sup_{(t,x)\in[0,T)\times\R^d}\frac{|v(t,x)|}{1+|x|} \ < \ \infty.
\end{equation}
The present subsection is devoted to the remaining part of the proof of Theorem \ref{ThmMain}, namely that $v$ is a viscosity solution to \eqref{HJB}-\eqref{condterminale}.

\vspace{3mm}

\noindent\textbf{Proof of the viscosity supersolution property to \eqref{HJB}.} We know that $v$ is the pointwise limit of the nondecreasing sequence of functions $(v_n^h)_n$, so that $v$ is lower semicontinuous and we have
\begin{equation}
\label{E:ProofMainThm1}
v(t,x) \ = \ v_*(t,x) \ = \ \liminf_{n\rightarrow\infty}\!{}_*\,v_n^h(t,x,a),
\end{equation}
for all $(t,x,a)\in[0,T)\times\R^d\times\mathring{A}_h$. Let $(t,x)\in[0,T)\times\R^d$ and $\varphi\in C^{1,2}([0,T]\times\R^d)$ such that
\[
(v-\varphi)(t,x) \ = \ \min_{[0,T]\times\R^d}(v-\varphi).
\]
From the linear growth condition \eqref{E:LinGrowth_v} on $v$, we can assume, without loss of generality, that $\varphi$ satisfies $\sup_{(t,x)\in[0,T]\times\R^d}|\varphi(t,x)|/(1+|x|)<\infty$. Fix some $a\in\mathring{A}_h$ and define, for any $\eps>0$, the test function
\[
\varphi^\eps(t',x',a') \ = \ \varphi(t',x') - \eps\big(|t'-t|^2 + |x'-x|^2 + |a'-a|^2\big),
\]
for all $(t',x',a')\in[0,T]\times\R^d\times\R^q$. Notice that $\varphi^\eps\leq\varphi$ with equality if and only if $(t',x',a')=(t,x,a)$, therefore $v-\varphi^\eps$ has a strict global minimum at $(t,x,a)$. From the linear growth condition on the continuous functions $v_n^h$ and $\varphi$, there exists a bounded sequence $(t_n,x_n,a_n)_n$ (we omit the dependence in $\eps$) in $[0,T)\times\R^d\times\R^q$ such that
\[
(v_n^h-\varphi^\eps)(t_n,x_n,a_n) \ = \ \min_{[0,T]\times\R^d\times\R^q}(v_n^h-\varphi^\eps).
\]
By standard arguments, we obtain that, up to a subsequence,
\[
\big(t_n,x_n,a_n,v_n^h(t_n,x_n,a_n)\big) \ \longrightarrow \ \big(t,x,a,v(t,x)\big), \qquad \text{as }n\rightarrow\infty.
\]
Now, from the viscosity supersolution property of $v_n^h$ at $(t_n,x_n,a_n)$ with the test function $\varphi_\eps$, we have
\begin{align*}
-\frac{\partial \varphi^\eps}{\partial t}(t_n,x_n,a_n) - \Lc^{h(a_n)} \varphi^\eps(t_n,x_n,a_n) - f(x_n,h(a_n))& \\
- \frac{1}{2}\textup{tr}\big(D_a^2 \varphi^\eps(t_n,x_n,a_n)\big) - n \big|D_a \varphi^\eps(t_n,x_n,a_n)\big|& \ \geq \ 0.
\end{align*}
Therefore
\[
-\frac{\partial \varphi^\eps}{\partial t}(t_n,x_n,a_n) - \Lc^{h(a_n)} \varphi^\eps(t_n,x_n,a_n) - f(x_n,h(a_n)) - \frac{1}{2}\textup{tr}\big(D_a^2 \varphi^\eps(t_n,x_n,a_n)\big) \ \geq \ 0.
\]
Sending $n$ to infinity in the above inequality, we obtain, from the definition of $\varphi^\eps$,
\[
-\frac{\partial \varphi^\eps}{\partial t}(t,x,a) - \Lc^{h(a)} \varphi^\eps(t,x,a) - f(x,h(a)) + \eps \ \geq \ 0.
\]
Sending $\eps$ to zero, recalling that $\varphi^\eps(t,x,a)=\varphi(t,x)$, we find
\[
-\frac{\partial \varphi}{\partial t}(t,x) - \Lc^{h(a)} \varphi(t,x) - f(x,h(a)) \ \geq \ 0.
\]
Since $a\in\mathring{A}_h$ and $h(\mathring{A}_h)=\mathring{A}$, the above equation can be rewritten in an equivalent way as follows
\[
-\frac{\partial \varphi}{\partial t}(t,x) - \Lc^a \varphi(t,x) - f(x,a) \ \geq \ 0,
\]
where $a$ is arbitrarily chosen in $\mathring{A}$. As a consequence, using assumption \textbf{(H$A$)} and the continuity of the coefficients $b$, $\sigma$, $\beta$, and $f$ in the variable $a$, we end up with
\[
-\frac{\partial \varphi}{\partial t}(t,x) - \sup_{a\in A}\Big[\Lc^a \varphi(t,x) - f(x,a)\Big] \ \geq \ 0,
\]
which is the viscosity supersolution property.
\ep

\vspace{3mm}

\noindent\textbf{Proof of the viscosity subsolution property to \eqref{HJB}.} Since $v$ is the pointwise limit of the nondecreasing sequence $(v_n^h)_n$, we have
(see, e.g., page 91 in \cite{bar94}):
\begin{equation}
\label{E:ProofMainThm2}
v^*(t,x) \ = \ \limsup_{n\rightarrow\infty}\!{}_*\,v_n^h(t,x,a),
\end{equation}
for all $(t,x,a)\in[0,T)\times\R^d\times\mathring{A}_h$, where
\[
\limsup_{n\rightarrow\infty}\!{}_*\,v_n^h(t,x,a) \ = \ \limsup_{\substack{n\rightarrow\infty\\
(t',x',a')\rightarrow(t,x,a)\\ t'<T,\,a'\in\mathring{A}_h}}v_n^h(t',x',a'), \qquad (t,x,a)\in[0,T)\times\R^d\times\R^q.
\]
Let $(t,x)\in[0,T)\times\R^d$ and $\varphi\in C^{1,2}([0,T]\times\R^d)$ such that
\[
(v^*-\varphi)(t,x) \ = \ \max_{[0,T]\times\R^d}(v^*-\varphi).
\]
We may assume, without loss of generality, that this maximum is strict and that $\varphi$ satisfies a linear growth condition $\sup_{(t,x)\in[0,T]\times\R^d}|\varphi(t,x)|/(1+|x|)<\infty$. Fix $a\in\mathring{A}_h$ and consider a sequence $(t_n,x_n,a_n)_n$ in $[0,T)\times\R^d\times\mathring{A}_h$ such that
\[
\big(t_n,x_n,a_n,v_n(t_n,x_n,a_n)\big) \ \longrightarrow \ \big(t,x,a,v^*(t,x)\big), \qquad \text{as }n\rightarrow\infty.
\]
Let us define for $n\geq1$ the function $\varphi_n\in C^{1,2}([0,T]\times(\R^d\times\R^q))$ by
\[
\varphi_n(t',x',a') \ = \ \varphi(t',x') + n\big(|t'-t_n|^2 + |x'-x_n|^2\big),
\]
for all $(t',x',a')\in[0,T]\times\R^d\times\R^q$. From the linear growth condition on $v_n^h$ and $\varphi$, we can find a sequence $(\bar t_n,\bar x_n,\bar a_n)_n$ in $[0,T)\times\R^d\times A_h$ such that
\[
(v_n^h-\varphi_n)(\bar t_n,\bar x_n,\bar a_n) \ = \ \max_{[0,T]\times\R^d\times A_h}(v_n^h-\varphi_n).
\]
By standard arguments, we obtain that, up to a subsequence,
\[
n\big(|\bar t_n-t_n|^2 + |\bar x_n-x_n|^2\big) \ \overset{n\rightarrow\infty}{\longrightarrow} \ 0.
\]
As a consequence, up to a subsequence, we have
\[
(\bar t_n,\bar x_n,\bar a_n) \ \overset{n\rightarrow\infty}{\longrightarrow} \ (t,x,\bar a),
\]
for some $\bar a\in A_h$. Now, from the viscosity subsolution property of $v_n^h$ at $(\bar t_n,\bar x_n,\bar a_n)$ with the test function $\varphi_n$, we have:
\begin{align*}
-\frac{\partial \varphi_n}{\partial t}(\bar t_n,\bar x_n,\bar a_n) - \Lc^{h(\bar a_n)} \varphi_n(\bar t_n,\bar x_n,\bar a_n) - f(\bar x_n,h(\bar a_n))& \\
- \frac{1}{2}\textup{tr}\big(D_a^2 \varphi_n(\bar t_n,\bar x_n,\bar a_n)\big) - n \big|D_a \varphi_n(\bar t_n,\bar x_n,\bar a_n)\big|& \ \leq \ 0.
\end{align*}
Therefore, using the definition of $\varphi_n$,
\[
-\frac{\partial \varphi_n}{\partial t}(\bar t_n,\bar x_n,\bar a_n) - \Lc^{h(\bar a_n)} \varphi_n(\bar t_n,\bar x_n,\bar a_n) - f(\bar x_n,h(\bar a_n)) \ \leq \ 0.
\]
Sending $n$ to infinity in the above inequality, we obtain
\[
-\frac{\partial \varphi}{\partial t}(t,x) - \Lc^{h(\bar a)} \varphi(t,x) - f(x,h(\bar a)) \ \leq \ 0.
\]
Setting $\tilde a=h(\bar a)$, the above equation can be rewritten in an equivalent way as follows
\[
-\frac{\partial \varphi}{\partial t}(t,x) - \Lc^{\tilde a} \varphi(t,x) - f(x,\tilde a) \ \leq \ 0.
\]
As a consequence, we have
\[
-\frac{\partial \varphi}{\partial t}(t,x) - \sup_{a\in A}\Big[\Lc^a \varphi(t,x) - f(x,a)\Big] \ \leq \ 0,
\]
which is the viscosity subsolution property.
\ep

\vspace{3mm}

\noindent\textbf{Proof of the viscosity supersolution property to \eqref{condterminale}.} Let $x\in\R^d$. From \eqref{E:ProofMainThm1}, we can find a sequence $(t_n,x_n,a_n)_n$ valued in $[0,T)\times\R^d\times\R^q$ such that
\[
\big(t_n,x_n,a_n,v_n^h(t_n,x_n,a_n)\big) \ \longrightarrow \ \big(T,x,a,v_*(T,x)\big), \qquad \text{as }n\rightarrow\infty,
\]
for some $a\in\mathring{A}_h$. Since the sequence $(v_n^h)_n$ is nondecreasing and $v_n^h(T,\cdot,\cdot)=g$, we have
\[
v_*(T,x) \ \geq \ \lim_{n\rightarrow\infty} v_1^h(t_n,x_n,a_n) \ = \ g(x).
\]
\ep

\vspace{3mm}

\noindent\textbf{Proof of the viscosity subsolution property to \eqref{condterminale}.} Let $x\in\R^d$. From \eqref{E:ProofMainThm2}, for every $\eps>0$ and $a\in\mathring{A}_h$ there exist $N\in\mathbb N$ and $\delta>0$ such that
\begin{equation}
\label{E:ProofMainThm3}
\big|v_n^h(t',x',a')-v^*(T,x)\big| \ \leq \ \eps,
\end{equation}
for all $n\geq N$ and $|t'-T|,|x'-x|,|a'-a|\leq\delta$, with $t'<T$ and $a'\in\mathring{A}_h$. Now, we recall that $v_n^h(T,x,a)=g(x)$, therefore, from the continuity of $v_n^h$, for every $n\in\N$, there exists $\delta_n>0$ such that
\begin{equation}
\label{E:ProofMainThm4}
\big|v_n^h(t',x',a') - g(x)\big| \ \leq \ \eps,
\end{equation}
for all $|t'-T|,|x'-x|,|a'-a|\leq\delta_n$, with $a'\in\mathring{A}_h$. Combining \eqref{E:ProofMainThm3} with \eqref{E:ProofMainThm4}, we end up with
\[
v^*(T,x) \ \leq \ g(x) + 2\eps.
\]
From the arbitrariness of $\eps$, we get the thesis.
\ep

\appendix

\setcounter{equation}{0} \setcounter{Assumption}{0}
\setcounter{Theorem}{0} \setcounter{Proposition}{0}
\setcounter{Corollary}{0} \setcounter{Lemma}{0}
\setcounter{Definition}{0} \setcounter{Remark}{0}

\renewcommand\thesection{Appendices}

\section{}

\renewcommand\thesection{\Alph{subsection}}

\renewcommand\thesubsection{\Alph{subsection}.}

\subsection{Martingale representation theorem}

We present here a martingale representation theorem, which is one of the fundamental result to derive our nonlinear Feynman-Kac representation formula. It is indeed a direct consequence of Theorem 4.29, Chapter III, in \cite{jacshiryaev03}, which is however designed for local (instead of square integrable) martingales. 

\begin{Theorem}
\label{MartReprThm}
Let $(t,a)\in[0,T]\times\R^q$ and $M$ $=$ $(M_s)_{t \leq s \leq T}$ be a c\`{a}dl\`{a}g square integrable $\F$-martingale, with $M_t$ constant. Then, there exist $Z\in{\bf L_{t,a}^2(W)}$, $V\in{\bf L_{t,a}^2(B)}$, and $U\in{\bf L_{t,a}^2(\tilde\pi)}$ such that
\[
M_s \ = \ M_t + \int_t^s Z_r dW_r + \int_t^s V_r dB_r + \int_t^s \int_E U_r(e) \tilde\pi(dr,de),
\]
for all $t \leq s \leq T$, $\P^{t,a}$ almost surely.
\end{Theorem}
{\bf Proof.}
Since  $M$ is a local martingale, we know from Theorem 4.29, Chapter III, in \cite{jacshiryaev03}, that
\[
M_s \ = \ M_t + \int_t^s Z_r dW_r + \int_t^s V_r dB_r + \int_t^s \int_E U_r(e) \tilde\pi(dr,de),
\]
for some predictable processes $(Z_s)_{t\leq s\leq T}$, $(V_s)_{t\leq s\leq T}$, and $(U_s)_{t\leq s\leq T}$, satisfying
\begin{align*}
&\E^{t,a}\bigg[\int_t^{T\wedge\tau_n^Z}|Z_s|^2 ds\bigg] \ < \ \infty, \qquad\qquad
\E^{t,a}\bigg[\int_t^{T\wedge\tau_n^V}|V_s|^2 ds\bigg] \ < \ \infty, \\
&\hspace{2cm}\E^{t,a}\bigg[\int_t^{T\wedge\tau_n^U}\int_E |U_s(e)|^2 \lambda(I_s^{t,a},de)ds\bigg] \ < \ \infty,
\end{align*}
for all $n\in\N$, where $(\tau_n^Z)_{n\in\N}$, $(\tau_n^V)_{n\in\N}$, and $(\tau_n^U)_{n\in\N}$ are nondecreasing sequences of $\F$-stopping times valued in $[t,T]$, converging pointwisely $\P^{t,a}$ a.s. to $T$. It remains to show that $Z\in{\bf L_{t,a}^2(W)}$, $V\in{\bf L_{t,a}^2(B)}$, and $U\in{\bf L_{t,a}^2(\tilde\pi)}$. To this end, set $\tau_n:=\tau_n^Z\wedge\tau_n^V\wedge\tau_n^U$, for every $n\in\N$. Notice that $\tau_n$ is an $\F$-stopping time valued in $[t,T]$, converging pointwisely $\P^{t,a}$ a.s. to $T$. Then, applying It\^o's formula to $M_s^2$ between $t$ and $\tau_n$, we find
\begin{align}
\label{E:MartReprThmProof1}
M_{\tau_n}^2 \ &= \ M_t^2 + 2\int_t^{\tau_n}M_sZ_sdW_s + 2\int_t^{\tau_n}M_sV_sdB_s + 2\int_t^{\tau_n}\int_E M_sU_s(e)\tilde\pi(ds,de) \notag \\
&\quad \ + \int_t^{\tau_n}|Z_s|^2ds + \int_t^{\tau_n}|V_s|^2ds + \int_t^{\tau_n}\int_E|U_s(e)|^2\pi(ds,de).
\end{align}
Observe that the local martingale $(\int_t^{s\wedge\tau_n}M_rZ_rdW_r)_{t\leq s\leq T}$ satisfies, using Burkholder-Davis-Gundy inequality and the fact that $\E^{t,a}[\sup_{t\leq s\leq T}|M_s|^2]<\infty$ (which is a consequence of Doob's inequality),
\[
\E^{t,a}\bigg[\sup_{t\leq s\leq T}\bigg|\int_t^{s\wedge\tau_n}M_rZ_rdW_r\bigg|\bigg] \ < \ \infty.
\]
In particular, $(\int_t^{s\wedge\tau_n}M_sZ_sdW_s)_{t\leq s\leq T}$ is a martingale. Similarly, $(\int_t^{s\wedge\tau_n}M_rV_rdB_r)_{t\leq s\leq T}$ and $(\int_t^{s\wedge\tau_n}\int_E M_rU_r(e)\tilde\pi(dr,de))_{t\leq s\leq T}$ are martingales. Therefore, taking the expectation in \eqref{E:MartReprThmProof1} yields
\begin{align}
\label{E:MartReprThmProof2}
\E^{t,a}\big[M_{\tau_n}^2\big] \ &= \ M_t^2 + \E^{t,a}\bigg[\int_t^{\tau_n}|Z_s|^2ds\bigg] + \E^{t,a}\bigg[\int_t^{\tau_n}|V_s|^2ds\bigg] \notag \\
&\quad \ + \E^{t,a}\bigg[\int_t^{\tau_n}\int_E|U_s(e)|^2\pi(ds,de)\bigg].
\end{align}
Recall that
\[
\E^{t,a}\bigg[\int_t^{\tau_n}\int_E|U_s(e)|^2\pi(ds,de)\bigg] \ = \ \E^{t,a}\bigg[\int_t^{\tau_n}\int_E|U_s(e)|^2\lambda(I_s^{t,a},de)ds\bigg].
\]
Moreover, we have $\E^{t,a}[M_{\tau_n}^2]\leq\E^{t,a}[\sup_{t\leq s\leq T}M_s^2]<\infty$. Therefore, from \eqref{E:MartReprThmProof2} it follows that there exists a positive constant $C$, independent of $n$, such that
\[
\E^{t,a}\bigg[\int_t^{\tau_n}|Z_s|^2ds\bigg] + \E^{t,a}\bigg[\int_t^{\tau_n}|V_s|^2ds\bigg] + \E^{t,a}\bigg[\int_t^{\tau_n}\int_E|U_s(e)|^2\lambda(I_s^{t,a},de)ds\bigg] \ \leq \ C.
\]
Letting $n\rightarrow\infty$, by Fatou's lemma we conclude that $Z\in{\bf L_{t,a}^2(W)}$, $V\in{\bf L_{t,a}^2(B)}$, and $U\in{\bf L_{t,a}^2(\tilde\pi)}$.
\ep

\setcounter{equation}{0} \setcounter{Assumption}{0}
\setcounter{Theorem}{0} \setcounter{Proposition}{0}
\setcounter{Corollary}{0} \setcounter{Lemma}{0}
\setcounter{Definition}{0} \setcounter{Remark}{0}

\subsection{Characterization of $\pi$ and Markov property of $(X,I)$}
\label{AppSubS:Markov}

In the following lemma, inspired by the results concerning Poisson random measures (see, e.g., Proposition 1.12, Chapter XII, in \cite{revuzyor99}), we present a characterization of $\pi$ in terms of Fourier and Laplace functionals. This shows that $\pi$ is a conditionally Poisson random measure (also known as doubly stochastic Poisson random measure or Cox random measure) relative to $\sigma(I_z;z\geq0)$.

\begin{Proposition}[Fourier and Laplace functionals of $\pi$]
\label{P:Fourier}
Assume that \textup{\textbf{(HFC)}} holds and fix $(t,a)\in[0,T]\times\R^q$. Let $\ell\colon\R_+\times E\rightarrow\R$ be a $\Bc(\R_+)\otimes\Bc(E)$-measurable function such that $\int_0^\infty\int_E|\ell_u(e)|\lambda(I_u^{t,a},de)du<\infty$, $\P^{t,a}$ a.s., then, for every $s\leq\infty$,
\[
\E^{t,a}\Big[e^{i\int_0^s\int_E\ell_u(e)\pi(du,de)}\Big|\sigma(I_z^{t,a};z\geq0)\Big] \ = \ e^{\int_0^s\int_E(e^{i\ell_u(e)}-1)\lambda(I_u^{t,a},de)du}, \qquad \P^{t,a}\,a.s.
\]
If $\ell$ is nonnegative, then the following equality holds:
\[
\E^{t,a}\Big[e^{-\int_0^s\int_E\ell_u(e)\pi(du,de)}\Big|\sigma(I_z^{t,a};z\geq0)\Big] \ = \ e^{-\int_0^s\int_E(1-e^{-\ell_u(e)})\lambda(I_u^{t,a},de)du}, \qquad \P^{t,a}\,a.s.
\]
In particular, if $(F_k)_{1\leq k\leq n}$, with $n\in\N\backslash\{0\}$, is a finite sequence of pairwise disjoint Borel measurable sets from $\R_+\times E$, with $\int_{F_k}\lambda(I_u^{t,a},de)du<\infty$, $\P^{t,a}$ a.s., then
\[
\E^{t,a}\Big[e^{i\sum_{k=1}^n\theta_k\pi(F_k)}\Big|\sigma(I_z^{t,a};z\geq0)\Big] \ = \ \Prod_{k=1}^n e^{\int_{F_k}(e^{i\theta_k}-1)\lambda(I_u^{t,a},de)du}, \qquad \P^{t,a}\,a.s.
\]
for all $\theta_1,\ldots,\theta_n\in\R$. In other words, $\pi(F_1),\ldots,\pi(F_n)$ are conditionally independent relative to $\sigma(I_z^{t,a};z\geq0)$.
\end{Proposition}
\textbf{Proof.}
Let $J_s = \int_0^s\int_E\ell_u(e)\pi(du,de)$, for any $s\geq0$, and define
\[
\phi(s) \ = \ \E^{t,a}\big[e^{iJ_s}\big|\sigma(I_z^{t,a};z\geq0)\big], \qquad \forall\,s\geq0.
\]
Applying It\^o's formula to the process $e^{iJ_s}$, we find
\[
e^{iJ_s} \ = \ 1 + \int_0^s\int_E e^{iJ_{u^-}}\big(e^{i\ell_u(e)} - 1\big) \pi(du,de).
\]
Taking the conditional expectation with respect to $\sigma(I_u^{t,a};u\geq0)$, we get
\begin{align*}
\E^{t,a}\big[e^{iJ_s}\big|\sigma(I_z^{t,a};z\geq0)\big] \ &= \ 1 + \E^{t,a}\bigg[\int_0^s\int_E e^{iJ_{u^-}}\big(e^{i\ell_u(e)} - 1\big) \lambda(I_u^{t,a},de)du\bigg|\sigma(I_z^{t,a};z\geq0)\bigg] \\
&= \ 1 + \int_0^s\int_E \E^{t,a}\big[e^{iJ_{u^-}}\big|\sigma(I_z^{t,a};z\geq0)\big]\big(e^{i\ell_u(e)} - 1\big) \lambda(I_u^{t,a},de)du.
\end{align*}
In terms of $\phi$ this reads
\[
\phi(s) \ = \ 1 + \int_0^s \phi(u^-)\psi(u)du, \qquad \P^{t,a}\,a.s.,
\]
where
\[
\psi(u) \ = \ \int_E \big(e^{i\ell_u(e)} - 1\big) \lambda(I_u^{t,a},de), \qquad \P^{t,a}\,a.s.
\]
Notice that $\psi$ belongs to $\mathbf{L^1}(\R_+)$, as a consequence of the integrability condition on $f$. We see then that $\phi$ is continuous, so that
\[
\phi(s) \ = \ e^{\int_0^s\psi(u)du}, \qquad \P^{t,a}\,a.s.,
\]
which yields the first formula of the lemma. The second formula is proved similarly.
\ep

\vspace{3mm}

We shall now study the Markov properties of the pair $(X,I)$ in the following two propositions.

\begin{Proposition}
\label{P:Markov}
Under assumption \textup{\textbf{(HFC)}}, for every $(t,x,a)\in[0,T]\times\R^d\times\R^q$ the stochastic process $(X_s^{t,x,a},I_s^{t,a})_{s\geq0}$ on $(\Omega,\Fc,\F,\P^{t,a})$ is Markov with respect to $\F$$:$ for every $r,s\in\R_+$, $r \leq s$, and for every Borel measurable and bounded function $h\colon\R^d\times\R^q\rightarrow\R$ we have
\beqs
\E^{t,a}\big[h(X_s^{t,x,a},I_s^{t,a})\big|\Fc_r\big] & = & \E^{t,a}\big[h(X_s^{t,x,a},I_s^{t,a})\big|\sigma(X_r^{t,x,a},I_r^{t,a})\big], \qquad \P^{t,a}\,a.s.
\enqs
\end{Proposition}
\textbf{Proof.}
Fix $(t,x,a)\in[0,T]\times\R^d\times\R^q$. Notice that it is enough to show the Markov property for $t \leq r \leq s \leq T$. Therefore, let $r\in[t,T]$ and consider, on $(\Omega,\Fc,\F,\P^{t,a})$, the following equation for $\tilde X$:
\beq
\label{E:tildeX}
\tilde X_s &=& X_r^{t,x,a} + \int_r^s b(\tilde X_u,I_u^{t,a}) du + \int_r^s \sigma(\tilde X_u,I_u^{t,a}) dW_u \\
& & + \; \int_r^s \int_E \beta(\tilde X_{u^-},I_u^{t,a},e) \tilde\pi(du,de), \notag
\enq
for all $s\in[r,T]$, $\P^{t,a}$ a.s., where $\tilde\pi(du,de)$ $=$ $\pi(du,de)$ $-$ $1_{\{u<T_\infty\}}\lambda(I_u^{t,a},de)du$. Under assumption \textbf{(HFC)}, it is known (see, e.g., Theorem 14.23 in \cite{jac79}) that there exists a unique solution to equation \eqref{E:tildeX}, which is clearly given by the process $(X_s^{t,x,a})_{s\in[r,T]}$. We recall that this solution is constructed using an iterative procedure, which relies on a recursively defined sequence of processes $(\tilde X^{(n)})_n$, see, e.g., Lemma 14.20 in \cite{jac79}. More precisely, we set $\tilde X^{(0)}\equiv0$ and then we define $\tilde X^{(n+1)}$ from $\tilde X^{(n)}$ as follows:
\beqs
\tilde X_s^{(n+1)} &=& X_r^{t,x,a} + \int_r^s b(\tilde X_u^{(n)},I_u^{t,a}) du + \int_r^s \sigma(\tilde X_u^{(n)},I_u^{t,a}) dW_u \\
& & + \; \int_r^s \int_E \beta(\tilde X_{u^-}^{(n)},I_{u^-}^{t,a},e) \tilde\pi(du,de), \notag
\enqs
for all $s\in[r,T]$, $\P^{t,a}$ a.s., for every $n\in\N$. It can be shown that $\tilde X^{(n)}$ converges uniformly towards the solution $X^{t,x,a}$ of \eqref{E:tildeX} on $[r,T]$, $\P^{t,a}$ a.s., namely $\sup_{s\in[r,T]}|\tilde X_s^{(n)}-X_s^{t,x,a}|\rightarrow0$ as $n$ tends to infinity, $\P^{t,a}$ almost surely. This shows that $X_s^{t,x,a}$ (and also $(X_s^{t,x,a},I_s^{t,a})$) is $\tilde\F$-adapted, where $\tilde\F = (\tilde\Fc_s)_{s\in[r,T]}$ is the augmentation of the filtration $\tilde\G = (\tilde\Gc_s)_{s\in[r,T]}$ given by:
\[
\tilde\Gc_s \ = \ \sigma(X_r^{t,x,a},I_r^{t,a}) \vee \Fc_{[r,s]}^W \vee \Fc_{[r,s]}^B \vee \Fc_{[r,s]}^\pi,
\]
where $\Fc_{[r,s]}^W = \sigma(W_u-W_r;r \leq u \leq s)$, $\Fc_{[r,s]}^B = \sigma(B_u-B_r;r \leq u \leq s)$, and $\Fc_{[r,s]}^\pi = \sigma(\pi(F);F\in\Bc([r,s])\otimes\Bc(E))$. Since $\Fc_{[r,s]}^W$ and $\Fc_{[r,s]}^B$ are independent with respect to $\Fc_r$, it is enough to prove that $\Fc_{[r,s]}^\pi$ and $\Fc_r$ are conditionally independent relative to $\sigma(X_r^{t,x,a},I_r^{t,a})$. To prove this, take $C\in\Fc_r$ and a $\Bc(\R_+)\otimes\Bc(E)$-measurable function $\ell\colon\R_+\times E\rightarrow\R$ such that $\int_0^\infty\int_E|\ell_u(e)|\lambda(I_u^{t,a},de)du<\infty$, $\P^{t,a}$ almost surely. Then, the thesis follows if we prove that
\begin{align}
\label{E:IndipCell}
&\E^{t,a}\Big[e^{i\theta_1 1_C + i\theta_2 \int_r^s\int_E\ell_u(e)\pi(du,de)}\Big|\sigma(X_r^{t,x,a},I_r^{t,a})\Big] \\
&= \ \E^{t,a}\big[e^{i\theta_1 1_C}\big|\sigma(X_r^{t,x,a},I_r^{t,a})\big] \E^{t,a}\Big[e^{i\theta_2 \int_r^s\int_E \ell_u(e) \pi(du,de)}\Big|\sigma(X_r^{t,x,a},I_r^{t,a})\Big], \qquad \P^{t,a}\,a.s., \notag
\end{align}
for all $\theta_1,\theta_2\in\R$. Firstly, let us prove that $1_C$ and $\int_r^s\int_E\ell_u(e)\pi(du,de)$ are conditionally independent relative to $\sigma(I_z^{t,a};z\geq r)$, i.e.,
\begin{align}
\label{E:IndipCellz}
&\E^{t,a}\Big[e^{i\theta_1 1_C + i\theta_2 \int_r^s\int_E\ell_u(e)\pi(du,de)}\Big|\sigma(I_z^{t,a};z\geq r)\Big] \\
&= \ \E^{t,a}\big[e^{i\theta_1 1_C}\big|\sigma(I_z^{t,a};z\geq r)\big] e^{\int_r^s\int_E(e^{i\ell_u(e)\theta_2}-1)\lambda(I_u^{t,a},de)du}, \qquad \P^{t,a}\,a.s. \notag
\end{align}
Proceeding as in Proposition \ref{P:Fourier}, let $J_s = \int_r^s\int_E\ell_u(e)\pi(du,de)$ and
\[
\phi(s) \ = \ \E^{t,a}\big[e^{i\theta_1 1_C + i\theta_2 J_s}\big|\sigma(I_z^{t,a};z\geq r)\big], \qquad \forall\,s\geq r.
\]
Applying It\^o's formula to the process $e^{iJ_s}$, we find
\begin{align*}
&\E^{t,a}\big[e^{i\theta_1 1_C + i\theta_2 J_s}\big|\sigma(I_z^{t,a};z\geq r)\big] \ = \ \E^{t,a}\big[e^{i\theta_1 1_C}\big|\sigma(I_z^{t,a};z\geq r)\big] \\
& + \ \E^{t,a}\bigg[\int_r^s\int_E e^{i\theta_1 1_C + i\theta_2 J_{u^-}}\big(e^{i\ell_u(e)\theta_2} - 1\big) \lambda(I_u^{t,a},de)du\bigg|\sigma(I_z^{t,a};z\geq0)\bigg] \\
&= \ \E^{t,a}\big[e^{i\theta_1 1_C}\big|\sigma(I_z^{t,a};z\geq r)\big] \\
&+ \ \int_r^s\int_E \E^{t,a}\big[e^{i\theta_1 1_C + i\theta_2 J_{u^-}}\big|\sigma(I_z^{t,a};z\geq r)\big]\big(e^{i\ell_u(e)\theta_2} - 1\big) \lambda(I_u^{t,a},de)du.
\end{align*}
In terms of $\phi$ this reads
\[
\phi(s) \ = \ 1 + \int_r^s \phi(u^-)\psi(u)du, \qquad \P^{t,a}\,a.s.,
\]
where
\[
\psi(u) \ = \ \int_E \big(e^{i\ell_u(e)\theta_2} - 1\big) \lambda(I_u^{t,a},de), \qquad \P^{t,a}\,a.s.
\]
Notice that $\psi$ belongs to $\mathbf{L^1}(\R_+)$, as a consequence of the integrability condition on $f$. We see then that $\phi$ is continuous, so that
\[
\phi(s) \ = \ \E^{t,a}\big[e^{i\theta_1 1_C}\big|\sigma(I_z^{t,a};z\geq r)\big] e^{\int_r^s\psi(u)du}, \qquad \P^{t,a}\,a.s.,
\]
which yields \eqref{E:IndipCellz}. Let us come back to \eqref{E:IndipCell}. We have, using \eqref{E:IndipCellz},
\[
\E^{t,a}\Big[e^{i\theta_1 1_C + i\theta_2 \int_r^s\int_E\ell_u(e)\pi(du,de)}\Big|\sigma(X_r^{t,x,a},I_r^{t,a})\Big] \ = \ \E^{t,a}[Y_1 Y_2|\sigma(X_r^{t,x,a},I_r^{t,a})],
\]
where
\begin{align*}
Y_1 \ &= \ \E^{t,a}\big[e^{i\theta_1 1_C}\big|\sigma(I_z^{t,a};z\geq r) \vee \sigma(X_r^{t,x,a},I_r^{t,a}) \big], \\
Y_2 \ &= \ \E^{t,a}\Big[e^{i\theta_2 \int_r^s\int_E\ell_u(e)\pi(du,de)} \Big| \sigma(I_z^{t,a};z\geq r) \vee \sigma(X_r^{t,x,a},I_r^{t,a}) \Big].
\end{align*}
Since $(I_z^{t,a})_{z\geq0}$ is Markov with respect to $\F$, we have that $\Fc_r$ and $\sigma(I_z^{t,a};z\geq r)$ are independent relative to $\sigma(I_r^{t,a})$. Therefore, $Y_1$ can be written as
\[
Y_1 \ = \ \E^{t,a}\big[e^{i\theta_1 1_C}\big|\sigma(X_r^{t,x,a},I_r^{t,a}) \big].
\]
It follows that $Y_1$ is $\sigma(X_r^{t,x,a},I_r^{t,a})$-measurable, so that
\[
\E^{t,a}\Big[e^{i\theta_1 1_C + i\theta_2 \int_r^s\int_E\ell_u(e)\pi(du,de)}\Big|\sigma(X_r^{t,x,a},I_r^{t,a})\Big] \ = \ Y_1\E^{t,a}[Y_2|\sigma(X_r^{t,x,a},I_r^{t,a})], \qquad \P^{t,a}\,a.s.,
\]
which proves \eqref{E:IndipCell}.
\ep

\begin{Proposition}
\label{P:MarkovFamily}
Under assumption \textup{\textbf{(HFC)}}, the family $(\Omega,\Fc,(X^{t,x,a},I^{t,a}),\P^{t,a})_{t,x,a}$ is Markovian with respect to $\F$ and satisfies, for every $(t,x,a)\in[0,T]\times\R^d\times\R^q$, $r,s\in\R_+$ with $r \leq s$, and for every Borel measurable and bounded function $h\colon\R^d\times\R^q\rightarrow\R$,
\begin{equation}
\label{E:MarkovPropertyFamily}
\E^{t,a}\big[h(X_s^{t,x,a},I_s^{t,a})\big|\Fc_r\big] \ = \ \int_{\R^d\times\R^q} h(x',a') p\big(r,(X_r^{t,x,a},I_r^{t,a}),s,dx'da'\big), \qquad \P^{t,a}\,a.s.
\end{equation}
where $p$ is the Markovian transition function given by
\[
p\big(r,(x',a'),s,\Gamma\big) \ = \ \P^{r,a'}\big((X_s^{r,x',a'},I_s^{r,a'})\in\Gamma\big),
\]
for every $r,s\in\R_+$, $r \leq s$, $(x',a')\in\R^d\times\R^q$, and every Borelian set $\Gamma\subset\R^d\times\R^q$. 
\end{Proposition}
\begin{Remark}
\label{R:TwoSol}
{\rm
For the proof of Proposition \ref{P:MarkovFamily} we shall need to consider simultaneously two distinct solutions $\{(X_s^{t,x,a},I_s^{t,a}),\,s\geq0\}$ and $\{(X_s^{t',x',a'},I_s^{t',a'}),\,s\geq0\}$, for $(t,x,a),(t',x',a')\in[0,T]\times\R^d\times\R^q$. According to Lemma \ref{L:ExistenceXI}, $\{(X_s^{t,x,a},I_s^{t,a}),\,s\geq0\}$ is defined on $(\Omega,\Fc,\F,\P^{t,a})$ and $\{(X_s^{t',x',a'},I_s^{t',a'}),\,s\geq0\}$ on $(\Omega,\Fc,\F,\P^{t',a'})$, respectively. However, we can construct a single probability space supporting both solutions. More precisely, we can construct a single probability space supporting both the random measure with compensator $1_{\{s<T_\infty\}}\lambda(I_s^{t,a},de)ds$ and the random measure with compensator $1_{\{s<T_\infty\}}\lambda(I_s^{t',a'},de)ds$, proceeding as follows.

Let $\Omega''$ be a copy of $\Omega'$, with corresponding canonical marked point process denoted by $(T_n'',\alpha_n'')_{n\in\N}$, canonical random measure $\pi''$, $T_\infty'' := \lim_n T_n''$, and filtration $\F''=(\Fc_s')_{t\geq0}$. Define $(\hat\Omega,\hat\Fc,\hat\F=(\hat\Fc_t)_{t\geq0})$ with $\hat\Omega$ $:=$ $\Omega\times\Omega''$, $\hat\Fc$ $:=$ $\Fc\otimes\Fc_\infty''$, and  $\hat\Fc_t$ $:=$ $\cap_{s>t}\Fc_s\otimes\Fc_s''$. Moreover, set $\hat W(\hat\omega)$ $:=$ $W(\omega)$, $\hat B(\hat\omega)$ $:=$ $B(\omega)$, $\hat\pi'(\hat\omega,\cdot)$ $:=$ $\pi(\omega,\cdot)$, and $\hat\pi''(\hat\omega,\cdot)$ $:=$ $\pi''(\omega'',\cdot)$. Set also $\hat T_\infty'(\hat\omega)$ $:=$ $T_\infty(\omega)$ and $\hat T_\infty''(\hat\omega)$ $:=$ $T_\infty''(\omega'')$. Let $\P^{t,a,t',a'}$ be the probability measure on $(\hat\Omega,\hat\Fc)$ given by $\P^{t,a,t',a'}(d\hat\omega)=\bar\P(d\bar\omega)\otimes\P^{',t,a}(\bar\omega,d\omega')\otimes\P^{'',t',a'}(\bar\omega,d\omega'')$. Finally, set $(\hat X^{t,x,a},\hat I^{t,a})(\hat\omega) := (X^{t,x,a},I^{t,a})(\bar\omega,\omega')$ and $(\hat X^{t',x',a'},\hat I^{t',a'})(\hat\omega) := (X^{t',x',a'},I^{t',a'})(\bar\omega,\omega'')$. Then $(\hat X^{t,x,a},\hat I^{t,a})$ solves \eqref{FSDEX}-\eqref{FSDEI} on $[t,T]$ starting from $(x,a)$ at time $t$, and  $(\hat X^{t',x',a'},\hat I^{t',a'})$ solves \eqref{FSDEX}-\eqref{FSDEI} on $[t',T]$ starting from $(x',a')$ at time $t'$.
\ep
}
\end{Remark}
\textbf{Proof (of Proposition \ref{P:MarkovFamily}).}
We begin noting that from Proposition \ref{P:Markov} the left-hand side of \eqref{E:MarkovPropertyFamily} is equal to $\E^{t,a}[h(X_s^{t,x,a},I_s^{t,a})|\sigma(X_r^{t,x,a},I_r^{t,a})]$, $\P^{t,a}$ almost surely. Let us now divide the proof into two steps.\\
\emph{Step 1. $(X_r^{t,x,a},I_r^{t,a})$ is a discrete random variable.} Suppose that
\[
(X_r^{t,x,a},I_r^{t,a}) \ = \ \sum_{i\geq1} (x_i,a_i) 1_{\Gamma_i},
\]
for some $(x_i,a_i)\in\R^d\times\R^q$ and a Borel partition $(\Gamma_i)_{i\geq1}$ of $\R^d\times\R^q$ satisfying $\P(\Gamma_i)>0$, for any $i\geq1$. In this case, \eqref{E:MarkovPropertyFamily} becomes
\begin{equation}
\label{E:MarkovPropertyFamily2}
\E^{t,a}\big[h(X_s^{t,x,a},I_s^{t,a})\big|\sigma(X_r^{t,x,a},I_r^{t,a})\big] \ = \ \sum_{i\geq1} 1_{\Gamma_i} \E^{r,a_i}\big[h(X_s^{r,x_i,a_i},I_s^{r,a_i})\big], \qquad \P^{t,a}\,a.s.
\end{equation}
Now notice that the process $(\hat X_s^{t,x,a}1_{\Gamma_i})_{s\geq r}$ satisfies on $(\hat\Omega,\hat\Fc,\hat\F,\P^{t,a,r,a_i})$ (using the same notation as in Remark \ref{R:TwoSol})
\begin{align*}
\hat X_s^{t,x,a}1_{\Gamma_i} \ &= \ x_i1_{\Gamma_i} + \int_r^s b_i(\hat X_u^{t,x,a}1_{\Gamma_i},\hat I_u^{t,a}1_{\Gamma_i})dr + \int_r^s \sigma_i(\hat X_u^{t,x,a}1_{\Gamma_i},\hat I_u^{t,a}1_{\Gamma_i}) d\hat W_u \\
&\quad \ + \int_r^s\int_E \beta(\hat X_{u^-}^{t,x,a}1_{\Gamma_i},\hat I_{u^-}^{t,a}1_{\Gamma_i},e) \tilde{\hat \pi}_i(du,de),
\end{align*}
with $b_i=b1_{\Gamma_i}$, $\sigma_i=\sigma1_{\Gamma_i}$, and $\tilde{\hat \pi}_i$ is the compensated martingale measure associated to the random measure $\hat \pi_i$, which has $1_{\Gamma_i}\lambda(\hat I_{s^-}^{t,a}1_{\Gamma_i},de)ds$, $s\geq r$, as compensator. Similarly, the process $(\hat X_s^{r,x_i,a_i}1_{\Gamma_i})_{s\geq r}$ satisfies on $(\hat\Omega,\hat\Fc,\hat\F,\P^{t,a,r,a_i})$
\begin{align*}
\hat X_s^{r,x_i,a_i}1_{\Gamma_i} \ &= \ x_i1_{\Gamma_i} + \int_r^s b_i(\hat X_u^{r,x_i,a_i}1_{\Gamma_i},\hat I_u^{r,a_i}1_{\Gamma_i})dr + \int_r^s \sigma_i(\hat X_u^{r,x_i,a_i}1_{\Gamma_i},\hat I_u^{r,a_i}1_{\Gamma_i}) d\hat W_u \\
&\quad \ + \int_r^s\int_E \beta(\hat X_{u^-}^{r,x_i,a_i}1_{\Gamma_i},\hat I_{u^-}^{r,a_i}1_{\Gamma_i},e) \tilde{\hat \pi}_i'(du,de),
\end{align*}
where $\tilde{\hat \pi}_i'$ is the compensated martingale measure associated to the random measure $\hat \pi_i'$, which has $1_{\Gamma_i}\lambda(\hat I_{s^-}^{r,a_i}1_{\Gamma_i},de)ds$, $s\geq r$, as compensator. Since the two processes $(\hat I_s^{t,a}1_{\Gamma_i})_{s\geq r}$ and $(\hat I_s^{r,a_i}1_{\Gamma_i})_{s\geq r}$ have the same law, we see that $(\hat X_s^{t,x,a}1_{\Gamma_i})_{s\geq r}$ and $(\hat X_s^{r,x_i,a_i}1_{\Gamma_i})_{s\geq r}$ solve the same equation, and, from uniqueness, they have the same law, as well. This implies (denoting $\E^{t,a,r,a_i}$ the expectation with respect to $\P^{t,a,r,a_i}$)
\[
\E^{t,a,r,a_i}\big[h(\hat X_s^{t,x,a},\hat I_s^{t,a})1_{\Gamma_i}\big] \ = \ \E^{t,a,r,a_i}\big[h(\hat X_s^{r,x_i,a_i},\hat I_s^{r,a_i})1_{\Gamma_i}\big].
\]
Notice that
\[
\E^{t,a,r,a_i}\big[h(\hat X_s^{t,x,a},\hat I_s^{t,a})1_{\Gamma_i}\big] \ = \ \E^{t,a}\big[h(X_s^{t,x,a},I_s^{t,a})1_{\Gamma_i}\big]
\]
and
\begin{align*}
\E^{t,a,r,a_i}\big[h(\hat X_s^{r,x_i,a_i},\hat I_s^{r,a_i})1_{\Gamma_i}\big] \ &= \ \E^{t,a,r,a_i}\big[\E^{t,a,r,a_i}\big[h(\hat X_s^{r,x_i,a_i},\hat I_s^{r,a_i})1_{\Gamma_i}\big|\Fc_r\big]\big] \\
&= \ \E^{t,a,r,a_i}\big[\E^{t,a,r,a_i}\big[h(\hat X_s^{r,x_i,a_i},\hat I_s^{r,a_i})\big|\Fc_r\big]1_{\Gamma_i}\big] \\
&= \ \E^{t,a,r,a_i}\big[\E^{t,a,r,a_i}\big[h(\hat X_s^{r,x_i,a_i},\hat I_s^{r,a_i})\big]1_{\Gamma_i}\big] \\
&= \ \E^{t,a}\big[\E^{r,a_i}\big[h(X_s^{r,x_i,a_i},I_s^{r,a_i})\big]1_{\Gamma_i}\big].
\end{align*}
In other words, we have
\[
\E^{t,a}\big[h(X_s^{t,x,a},I_s^{t,a})1_{\Gamma_i}\big] \ = \ \E^{t,a}\big[\E^{r,a_i}\big[h(X_s^{r,x_i,a_i},I_s^{r,a_i})\big]1_{\Gamma_i}\big],
\]
from which \eqref{E:MarkovPropertyFamily2} follows.\\
\emph{Step 2. General case.} From estimate \eqref{EstimateXI}, we see that $(X_r^{t,x,a},I_r^{t,a})$ is square integrable, so that there exists a sequence $(X_r^{t,x,a,n},I_r^{t,a,n})_n$ of square integrable discrete random variables converging to $(X_r^{t,x,a},I_r^{t,a})$ pointwisely $\P^{t,a}$ a.s. and in $L^2(\Omega,\Fc,\P^{t,a};\R^d\times\R^q)$. The sequence $(X_r^{t,x,a,n},I_r^{t,a,n})_n$ can be chosen in such a way that $(X_r^{t,x,a,n+1},I_r^{t,a,n+1})$ is a better approximation of $(X_r^{t,x,a},I_r^{t,a})$ than $(X_r^{t,x,a,n},I_r^{t,a,n})$, in other words such that $\sigma(X_r^{t,x,a,n},I_r^{t,a,n})\subset\sigma(X_r^{t,x,a,n+1},I_r^{t,a,n+1})$. Let us denote $(X_s^{t,x,a,n},I_s^{t,a,n})$ the solution to \eqref{FSDEX}-\eqref{FSDEI} starting at time $r$ from $(X_r^{t,x,a,n},I_r^{t,a,n})$. Notice that, from classical convergence results of diffusion processes with jumps (see, e.g., Theorem 4.8, Chapter IX, in \cite{jacshiryaev03}), it follows that $(X_s^{t,x,a,n},I_s^{t,a,n})$ converges weakly to $(X_s^{t,x,a},I_s^{t,a})$. From Step~1, for any $n$ we have
\begin{equation}
\label{E:MarkovPropertyFamilyProof1}
\E^{t,a}\big[h(X_s^{t,x,a,n},I_s^{t,a,n})\big|\sigma(X_r^{t,x,a,n},I_r^{t,a,n})\big] \ = \ p\big(r,(X_r^{t,x,a,n},I_r^{t,a,n}),s,h\big), \qquad \P^{t,a}\,a.s.
\end{equation}
where
\[
p(r,(x',a'),s,h) \ = \ \E^{r,a'}\big[h(X_s^{r,x',a',n},I_s^{r,a',n})\big],
\]
for every $r,s\in\R_+$, $r \leq s$, $(x',a')\in\R^d\times\R^q$, and every Borel measurable and bounded function $h\colon\R^d\times\R^q\rightarrow\R$. Let us suppose that $h$ is bounded and continuous. Since the sequence $(\E^{t,a}[h(X_s^{t,x,a,n},I_s^{t,a,n})|\sigma(X_r^{t,x,a,n},I_r^{t,a,n})])_n$ is uniformly bounded in $L^2(\Omega,\Fc,\P^{t,a})$, there exists a subsequence $(\E^{t,a}[h(X_s^{t,x,a,n_k},I_s^{t,a,n_k})|\sigma(X_r^{t,x,a,n_k},I_r^{t,a,n_k})])_k$ which converges weakly to some $Z\in L^2(\Omega,\Fc,\P^{t,a})$. For any $N\in\N$ and $\Gamma_N\in\sigma(X_r^{t,x,a,N},I_r^{t,a,N})$, we have, by definition of conditional expectation,
\[
\E^{t,a}\big[\E^{t,a}\big[h(X_s^{t,x,a,n_k},I_s^{t,a,n_k})\big|\sigma(X_r^{t,x,a,n_k},I_r^{t,a,n_k})\big]1_{\Gamma_N}\big] \ = \ \E^{t,a}\big[h(X_s^{t,x,a,n_k},I_s^{t,a,n_k})1_{\Gamma_N}\big],
\]
for all $n_k\geq N$. Letting $k\rightarrow\infty$, we deduce
\[
\E^{t,a}\big[Z1_{\Gamma_N}\big] \ = \ \E^{t,a}\big[h(X_s^{t,x,a},I_s^{t,a})1_{\Gamma_N}\big].
\]
Since $\sigma(X_r^{t,x,a},I_r^{t,a})=\vee_n\sigma(X_r^{t,x,a,n},I_r^{t,a,n})$, it follows that
\[
Z \ = \ \E^{t,a}[h(X_s^{t,x,a},I_s^{t,a})|\sigma(X_r^{t,x,a},I_r^{t,a})], \qquad \P^{t,a}\,a.s.
\]
Notice that every convergent subsequence of $(\E^{t,a}[h(X_s^{t,x,a,n},I_s^{t,a,n})|\sigma(X_r^{t,x,a,n},I_r^{t,a,n})])_n$ has to converge to $\E^{t,a}[h(X_s^{t,x,a},I_s^{t,a})|\sigma(X_r^{t,x,a},I_r^{t,a})]$, so that the whole sequence converges. On the other hand, when $h$ is bounded and continuous, it follows again from classical convergence results of diffusion processes with jumps (see, e.g., Theorem 4.8, Chapter IX, in \cite{jacshiryaev03}), that $p=p(r,(x',a'),s,h)$ is continuous in $(x',a')$. Since $(X_r^{t,x,a,n},I_r^{t,a,n})_n$ converges pointwisely $\P^{t,a}$ a.s. to $(X_r^{t,x,a},I_r^{t,a})$, letting $n\rightarrow\infty$ in \eqref{E:MarkovPropertyFamilyProof1} we obtain
\begin{equation}
\label{E:MarkovPropertyFamilyProof2}
\E^{t,a}\big[h(X_s^{t,x,a},I_s^{t,a})\big|\sigma(X_r^{t,x,a},I_r^{t,a})\big] \ = \ p\big(r,(X_r^{t,x,a},I_r^{t,a}),s,h\big), \qquad \P^{t,a}\,a.s.
\end{equation}
for any $h$ bounded and continuous. Using a monotone class argument, we conclude that \eqref{E:MarkovPropertyFamilyProof2} remains true for any $h$ bounded and Borel measurable.
\ep

\setcounter{equation}{0} \setcounter{Assumption}{0}
\setcounter{Theorem}{0} \setcounter{Proposition}{0}
\setcounter{Corollary}{0} \setcounter{Lemma}{0}
\setcounter{Definition}{0} \setcounter{Remark}{0}

\subsection{Comparison theorem for equation \eqref{HJB}-\eqref{condterminale}}

We shall prove a comparison theorem for viscosity sub and supersolutions to the fully nonlinear IPDE of HJB type \eqref{HJB}-\eqref{condterminale}. Inspired by Definition 2 in \cite{barlesimbert08}, we begin recalling the following result concerning an equivalent definition of viscosity super and subsolution to \reff{HJB}-\reff{condterminale}, whose standard proof is not reported.

\begin{Lemma}
\label{L:viscosityI1I2}
Let assumption \textup{\textbf{(HFC)}}, \textup{\textbf{(HBC)}}, and \textup{\textbf{(H$\lambda$)}} hold. A locally bounded and lsc $($resp. usc$)$ function $u$ on $[0,T]\times\R^d$ is a viscosity supersolution $($resp. viscosity subsolution$)$ to \eqref{HJB}-\eqref{condterminale} if and only if
\beqs
u(T,x) \ \geq \ (resp. \; \leq) \  g(x)
\enqs
for any $x\in \R^d$, and, for any $\delta>0$,
\begin{align*}
-\,\frac{\partial \varphi}{\partial t}(t,x) - \sup_{a \in A} \bigg[ b(x,a).D_x\varphi(t,x) + \frac{1}{2}\textup{tr}\big(\sigma\sigma\trans (x,a)D_x^2\varphi(t,x)\big) + I_a^{1,\delta}(t,x,\varphi) & \\
+ \; I_a^{2,\delta}(t,x,D_x\varphi(t,x),u) + f\big(x,a\big) \bigg] \ \geq & \ \;(resp. \; \leq) \ 0, 
\end{align*}
for any $(t, x) \in [0,T)\times \R^d$ and any $\varphi \in C^{1,2}([0,T]\times \R^d)$ such that
\beqs
(u - \varphi)(t,x) \  = \ \min_{[0,T]\times \R^d} (u - \varphi) \quad (resp. \ \max_{[0,T]\times \R^d} (u - \varphi)).
\enqs
\end{Lemma}
%
%

As in \cite{barlesimbert08}, see Definition 4, for the proof of the comparison theorem it is useful to adopt another equivalent definition of viscosity solution to equation \eqref{HJB}-\eqref{condterminale}, see Lemma \ref{L:viscosityJets} below, where we mix test functions and sub/superjets. We first recall the definition of sub and superjets.

\begin{Definition}
Let $u\colon[0,T] \times \R^d \rightarrow \R$ be a lsc $($resp. usc$)$ function.\\
\textup{(i)} We denote by $\Pc^{2,-} u(t,x)$ the \textbf{parabolic subjet} $($resp. $\Pc^{2,+} u(t,x)$ the \textbf{parabolic superjet}$)$ of $u$ at $(t,x) \in [0,T) \times \R^d$, as the set of triples $(p, q, M) \in \R \times \R^d \times \S^d$ $($we denote by $\S^d$ the set of $d\times d$ symmetric matrices$)$ satisfying
\beqs
u(s,y) \;\, \geq &(resp. \ \leq)& u(t,x) +p(s-t) +   q.(y-x) + \frac1{2} (y-x).M(y-x) \\
& &+ \; o\big(|s-t|+|y-x|^2\big), \qquad\qquad\; \text{as }(s,y)\rightarrow(t,x).
\enqs
\textup{(ii)} We denote by $\bar \Pc^{2,-} u(t,x)$ the \textbf{parabolic limiting subjet} $($resp. $\bar \Pc^{2,+} u(t,x)$ the \textbf{parabolic limiting superjet}$)$ of $u$ at $(t,x) \in [0,T) \times \R^d$, as the set of triples $(p, q, M) \in \R \times \R^d \times \S^d$ such that
\beqs
(p,q,M) &=& \lim_{n\rightarrow\infty} (p_n, q_n, M_n)
\enqs
with $(p_n, q_n, M_n) \in \Pc^{2,-} u(t_n,x_n)$ $($resp. $\Pc^{2,+} u(t_n,x_n)$$)$, where
\beqs
(t,x,u(t,x)) &=& \lim_{n\rightarrow\infty} (t_n, x_n, u(t_n,x_n)).
\enqs
\end{Definition}

\begin{Lemma}
\label{L:viscosityJets}
Let assumption \textup{\textbf{(HFC)}}, \textup{\textbf{(HBC)}}, and \textup{\textbf{(H$\lambda$)}} hold. A locally bounded and lsc $($resp. usc$)$ function $u$ on $[0,T]\times\R^d$ is a viscosity supersolution $($resp. viscosity subsolution$)$ to \eqref{HJB}-\eqref{condterminale} if and only if
\beqs
u(T,x) \ \geq \ (resp. \; \leq) \  g(x)
\enqs
for any $x\in \R^d$, and, for any $\delta>0$,
\beqs
-\;p - \sup_{a \in A} \bigg[ b(x,a).q + \frac{1}{2}\textup{tr}\big(\sigma\sigma\trans (x,a)M\big) + I_a^{1,\delta}(t,x,\varphi)  & & \\
+ \; I_a^{2,\delta}(t,x,q,u) + f\big(x,a\big) \bigg] & \geq & (resp. \; \leq) \;\;\;  0,
\enqs
for any $(t, x) \in [0,T)\times \R^d$, $(p,q,M)\in\bar\Pc^{2,-}u(t,x)$ $($resp. $(p,q,M)\in\bar\Pc^{2,+}u(t,x)$$)$, and any $\varphi \in C^{1,2}([0,T]\times \R^d)$, with $\frac{\partial \varphi}{\partial t}(t,x) = p$, $D_x\varphi(t,x) = q$, and $D_x^2 \varphi(t,x) \leq M$ $($resp. $D_x^2 \varphi(t,x) \geq M$$)$, such that
\beqs
(u - \varphi)(t,x) \  = \ \min_{[0,T]\times \R^d} (u - \varphi) \quad (resp. \ \max_{[0,T]\times \R^d} (u - \varphi)).
\enqs
\end{Lemma}
\textbf{Proof.}
Using Lemma \ref{L:viscosityI1I2}, we see that the \emph{if} part is true. We have to prove the \emph{only if} part. In particular, we prove the equivalence for the supersolution case only, since the subsolution case can be proved similarly.

Let $u$ be locally bounded and lsc on $[0,T]\times\R^d$ and suppose that $u$ is a viscosity supersolution to \eqref{HJB}-\eqref{condterminale}. Fix $\delta>0$, $(t,x)\in[0,T)\times\R^d$, $(p,q,M)\in\bar\Pc^{2,-} u(t,x)$ and $\varphi\in C^{1,2}([0,T]\times\R^d)$, with $\frac{\partial\varphi}{\partial t} = p$, $D_x\varphi(t,x) = q$, and $D_x^2\varphi(t,x) \leq M$, such that
\beqs
(u - \varphi)(t,x) \  = \ \min_{[0,T]\times \R^d} (u - \varphi).
\enqs
By classical results (see, e.g., Lemma 4.1, Chapter V, in \cite{flemso}), there exists a function $\psi\colon[0,T]\times\R^d\rightarrow\R$, $\psi\in C^{1,2}([0,T]\times\R^d)$, such that $\psi(t,x)=u(t,x)$, $\frac{\partial\psi}{\partial t}(t,x)=p$, $D_x\psi(t,x)=q$, $D_x^2\psi(t,x)=M$, and $\psi\leq u$ on $[0,T]\times\R^d$. For any $\eps>0$, we define $\psi_\eps$ as follows:
\[
\psi_\eps(s,y) \ = \ \chi_\eps(s,y)\psi(s,y) + (1 - \chi_\eps(s,y))\varphi(s,y), \qquad (s,y)\in[0,T]\times\R^d,
\]
where $\chi_\eps$ is a smooth function satisfying:
\begin{align*}
&0 \ \leq \ \chi_\eps(s,y) \ \leq \ 1, \qquad \text{if }(s,y)\in[0,T]\times\R^d, \\
&\chi_\eps(s,y) = 1, \qquad\qquad\quad\, \text{if }(s,y)\in([0,T]\cap\{|s-t|<\eps\})\times(\R^d\cap\{|y-x|<\eps\}), \\
&\chi_\eps(s,y) = 0, \qquad\qquad\quad\, \text{if }(s,y)\in([0,T]\cap\{|s-t|>2\eps\})\times(\R^d\cap\{|y-x|>2\eps\}).
\end{align*}
Notice that $\psi_\eps\in C^{1,2}([0,T]\times\R^d)$ and $\min_{[0,T]\times\R^d}(u-\psi_\eps) = (u-\psi_\eps)(t,x)$. Moreover, $\psi_\eps = \psi$ in a neighborhood of $(t,x)$. As a consequence, from Lemma \ref{L:viscosityI1I2} we have
\beq
\label{E:Proof_Equivalence}
-\;p - \sup_{a \in A} \bigg[ b(x,a).q + \frac{1}{2}\textup{tr}\big(\sigma\sigma\trans (x,a)M\big) + I_a^{1,\delta}(t,x,\psi_\eps)  & & \\
+ \; I_a^{2,\delta}(t,x,q,u) + f\big(x,a\big) \bigg] & \geq & 0. \notag
\enq
Let us assume, for a moment, the validity of the following result:
\begin{equation}
\label{E:psi_eps-varphi}
\sup_{a\in A}\big|I_a^{1,\delta}(t,x,\psi_\eps) - I_a^{1,\delta}(t,x,\varphi)\big| \ \overset{\eps\rightarrow0^+}{\longrightarrow} \ 0.
\end{equation}
Then, by sending $n\rightarrow\infty$ in \eqref{E:Proof_Equivalence}, we obtain the thesis
\beqs
-\;p - \sup_{a \in A} \bigg[ b(x,a).q + \frac{1}{2}\textup{tr}\big(\sigma\sigma\trans (x,a)M\big) + I_a^{1,\delta}(t,x,\varphi)  & & \\
+ \; I_a^{2,\delta}(t,x,q,u) + f\big(x,a\big) \bigg] & \geq & 0. \notag
\enqs
Therefore, it remains to prove \eqref{E:psi_eps-varphi}. Notice that
\begin{align}
\label{E:Ipsi-varphi}
&\sup_{a\in A}\big|I_a^{1,\delta}(t,x,\psi_\eps) - I_a^{1,\delta}(t,x,\varphi)\big| \\
&= \ \sup_{a\in A}\bigg|\int_{E\cap\{|e|\leq\delta\}} \chi_\eps(t,x+\beta(x,a,e))\big(\psi(t,x+\beta(x,a,e)) - \varphi(t,x+\beta(x,a,e))\big) \lambda(a,de)\bigg|. \notag
\end{align}
From the regularity of $\psi$ and $\varphi$, we have
\[
\big|\psi(t,x+\beta(x,a,e)) - \varphi(t,x+\beta(x,a,e))\big| \ \leq \ |\beta(x,a,e)|^2 \sup_{|y-x| \leq r_{\delta,x}}|D_x^2(\psi-\varphi)(t,y)|,
\]
where $r_{\delta,x} := \sup_{(a,e)\in A\times(E\cap\{|e|\leq\delta\})}|\beta(x,a,e)|$. In particular, \eqref{E:Ipsi-varphi} becomes (in the sequel we shall denote by $C$ a generic positive constant depending only on $\delta$ and $x$)
\begin{align*}
&\sup_{a\in A}\big|I_a^{1,\delta}(t,x,\psi_\eps) - I_a^{1,\delta}(t,x,\varphi)\big| \\
&\leq \ C\sup_{a\in A}\int_{E\cap\{|e|\leq\delta\}} \chi_\eps(t,x+\beta(x,a,e))|\beta(x,a,e)|^2 \lambda(a,de).
\end{align*}
Observe that $\chi_\eps(t,x+\beta(x,a,e))|\beta(x,a,e)|^2 \leq |\beta(x,a,e)|^2 1_{\{|\beta(x,a,e)| \leq 2\eps\}}$. Since $\beta(x,a,e) \leq C(1\wedge|e|^2)$, we find
\begin{equation}
\label{E:Ipsi-varphi3}
\sup_{a\in A}\big|I_a^{1,\delta}(t,x,\psi_\eps) - I_a^{1,\delta}(t,x,\varphi)\big| \ \leq \ C\sup_{a\in A}\int_{E\cap\{|e|\leq\delta\}} 1\wedge|e|^2\wedge(4\eps^2) \lambda(a,de).
\end{equation}
It follows from assumption \textbf{(H$\lambda$)}(i) that the right-hand side of \eqref{E:Ipsi-varphi3} goes to zero as $\eps\rightarrow0^+$, from which we deduce \eqref{E:psi_eps-varphi}.
\ep

\vspace{3mm}

We can now state the main result of this appendix.

\begin{Theorem}
\label{CompThm}
Assume that \textup{\textbf{(HFC)}}, \textup{\textbf{(HBC)}}, and \textup{\textbf{(H$\lambda$)}} hold. Let $u$ be a usc viscosity subsolution to \eqref{HJB}-\eqref{condterminale} and $w$ a lsc viscosity supersolution to \eqref{HJB}-\eqref{condterminale}, satisfying a linear growth condition
\beq
\label{linear}
\sup_{(t,x)\in[0,T]\times\R^d} \frac{|u(t,x)| + |w(t,x)|}{1+|x|} & < & \infty.
\enq
If $u(T,x) \leq w(T,x)$ for all $x\in\R^d$, then $u \leq w$ on $[0,T]\times\R^d$. 
\end{Theorem}
\textbf{Proof}
We shall argue by contradiction, assuming that
\beq
\label{vmoinsu}
\sup_{[0,T] \times \R^d} (u-w) & > & 0.
\enq
\emph{Step 1.} For some $\rho > 0$ to be chosen later, set
\beqs
\tilde u(t,x) \ = \ e^{\rho t}u(t,x), \qquad \tilde w(t,x) \ = \ e^{\rho t}w(t,x), \qquad (t,x)\in[0,T]\times\R^d.
\enqs
Let us consider the following equation:
\beq
\rho \tilde v -\frac{\partial \tilde v}{\partial t} - \sup_{a \in A} \big( \Lc^a \tilde v + \tilde f(\cdot,a) \big) &=& 0, \quad\;\;\; \text{on}\  [0,T) \times \R^d, \label{HJB2} \\
\tilde v(T,x) &=& \tilde g(x), \quad x\in \R^d, \label{condterm2}
\enq
where
\beqs
\tilde f(t,x,a) \ = \ e^{\rho t} f(x,a), \qquad \tilde g(x) \ = \ e^{\rho T} g(x),
\enqs
for all $(t,x,a) \in [0,T] \times \R^d \times A$. Then $\tilde u$ (resp. $\tilde w$) is a viscosity subsolution (resp. supersolution) to \eqref{HJB2}-\eqref{condterm2} (the definition of viscosity sub/supersolution to \eqref{HJB2}-\eqref{condterm2} is an obvious adaptation of Definition \ref{D:viscosity}). Indeed, concerning the subsolution property of $\tilde u$, let $(t, x) \in [0,T)\times \R^d$ and $\tilde\varphi \in C^{1,2}([0,T]\times \R^d)$ such that
\beqs
(\tilde u - \tilde\varphi)(t,x) \  = \ \max_{[0,T]\times \R^d} (\tilde u - \tilde\varphi).
\enqs
We can suppose $\tilde u(t,x) = \tilde\varphi(t,x)$, without loss of generality. Set $\varphi(s,y) = e^{-\rho s}\tilde\varphi(s,y)$, for all $(s,y)\in[0,T]\times\R^d$. Then $u(t,x) = \varphi(t,x)$. Moreover, since $\tilde u-\tilde\varphi\leq0$ on $[0,T]\times\R^d$, we see that $\max_{[0,T]\times\R^d}(u-\varphi)=0$. The claimed viscosity subsolution property of $\tilde u$ to \eqref{HJB2} then follows from the viscosity subsolution property of $u$ to \eqref{HJB}. Similarly, we can show the viscosity supersolution property of $\tilde w$.\\
\noindent\emph{Step 2.} Denote, for all $(t,s,x,y) \in [0,T]^2 \times \R^{2d}$, and for any $n\in\N\backslash\{0\}$ and $\gamma>0$,
\[
\Phi_{n,\gamma} (t,s,x,y) \ = \ \tilde u(t,x) - \tilde w(s,y) - n\frac{|t-s|^2}{2} - n\frac{|x-y|^2}{2} - \gamma\big(|x|^2 + |y|^2\big).
\]
By the linear growth assumption on $u$ and $w$, for each $n$ and $\gamma$, there exists $(t_{n,\gamma},s_{n,\gamma},x_{n,\gamma},y_{n,\gamma}) \in [0,T]^2\times\R^{2d}$ attaining the maximum of $\Phi_{n,\gamma}$ on $[0,T]^2\times\R^{2d}$. Notice that $\Phi_{n,\gamma}(t_{n,\gamma},s_{n,\gamma},x_{n,\gamma},y_{n,\gamma})\geq0$, for $\gamma$ small enough. Indeed, from \eqref{vmoinsu} we see that there exists $(\hat t,\hat x)\in[0,T)\times\R^d$ such that $\tilde u(\hat t,\hat x) - \tilde w(\hat t,\hat x) =: \eta > 0$. Then
\[
\Phi_{n,\gamma}(t_{n,\gamma},s_{n,\gamma},x_{n,\gamma},y_{n,\gamma}) \ \geq \ \Phi_{n,\gamma}(\hat t,\hat t,\hat x,\hat x) = \eta - 2\gamma|\hat x|^2,
\]
therefore it is enough to take $\gamma \leq \eta/(2|\hat x|^2)$. From $\Phi_{n,\gamma}(t_{n,\gamma},s_{n,\gamma},x_{n,\gamma},y_{n,\gamma})\geq0$ it follows that
\begin{equation}
\label{E:UniqProof1}
n\frac{|t_{n,\gamma}-s_{n,\gamma}|^2}{2} + n\frac{|x_{n,\gamma}-y_{n,\gamma}|^2}{2} + \gamma\big(|x_{n,\gamma}|^2 + |y_{n,\gamma}|^2\big) \ \leq \ \tilde u(t_{n,\gamma},x_{n,\gamma}) - \tilde w(s_{n,\gamma},y_{n,\gamma}).
\end{equation}
On the other hand, from the linear growth condition \eqref{linear} of $u$ and $w$, we deduce that there exists a constant $C>0$ such that (recalling the standard inequality $ab \leq a^2/(2\gamma) + \gamma b^2/2$, for any $a,b\in\R$ and $\gamma>0$)
\begin{align}
\label{E:UniqProof2}
\tilde u(t,x) - \tilde w(s,y) \ &\leq \ C\big(1 + |x| + |y|\big) \\
&\leq \ C + \frac{C^2}{\gamma} + \frac{\gamma}{2}\big(|x|^2 + |y|^2\big), \qquad \forall\,(t,s,x,y)\in[0,T]^2\times\R^{2d}. \notag
\end{align}
Combining \eqref{E:UniqProof1} with \eqref{E:UniqProof2}, we obtain
\begin{align*}
n\frac{|t_{n,\gamma}-s_{n,\gamma}|^2}{2} + n\frac{|x_{n,\gamma}-y_{n,\gamma}|^2}{2} + \gamma\big(|x_{n,\gamma}|^2 + |y_{n,\gamma}|^2\big) \ \leq \ \tilde u(t_{n,\gamma},x_{n,\gamma}) - \tilde w(s_{n,\gamma},y_{n,\gamma})& \\
\leq \ C + \frac{C^2}{\gamma} + \frac{\gamma}{2}\big(|x_{n,\gamma}|^2 + |y_{n,\gamma}|^2\big),&
\end{align*}
which implies
\begin{align}
\label{E:UniqProof3}
n\frac{|t_{n,\gamma}-s_{n,\gamma}|^2}{4} + n\frac{|x_{n,\gamma}-y_{n,\gamma}|^2}{4} + \frac{\gamma}{2}\big(|x_{n,\gamma}|^2 + |y_{n,\gamma}|^2\big) \ \leq \ C + \frac{C^2}{\gamma}.
\end{align}
From \eqref{E:UniqProof3} it follows that, for each $\gamma$, there exists $(t_\gamma,x_\gamma)\in[0,T]\times\R^d$ such that
\begin{align}
(t_{n,\gamma},s_{n,\gamma},x_{n,\gamma},y_{n,\gamma}) \ &\overset{n\rightarrow\infty}{\longrightarrow} \ (t_\gamma,t_\gamma,x_\gamma,x_\gamma), \label{E:UniqProof4} \\
n|x_{n,\gamma}-x_\gamma|^2 + n|y_{n,\gamma}-y_\gamma|^2 \ &\overset{n\rightarrow\infty}{\longrightarrow} \ 0, \label{E:UniqProof5} \\
\tilde u(t_{n,\gamma},x_{n,\gamma}) - \tilde w(s_{n,\gamma},y_{n,\gamma}) \ &\overset{n\rightarrow\infty}{\longrightarrow} \ \tilde u(t_\gamma,x_\gamma) - \tilde w(s_\gamma,y_\gamma).\label{E:UniqProof6}
\end{align}
As a matter of fact, we see from \eqref{E:UniqProof3} that, for every $\gamma$, there exists a constant $C_\gamma>0$ such that $|x_{n,\gamma}|,|y_{n,\gamma}| \leq C_\gamma$. Moreover, we obviously have $|t_{n,\gamma}|,|s_{n,\gamma}| \leq T$. Therefore, from Bolzano-Weierstrass theorem, there exist a subsequence $((t_{n_k,\gamma},s_{n_k,\gamma},x_{n_k,\gamma},y_{n_k,\gamma}))_k$ and $(t_\gamma,t_\gamma',x_\gamma,x_\gamma')\in[0,T]^2\times\R^{2d}$ such that $(t_{n_k,\gamma},s_{n_k,\gamma},x_{n_k,\gamma},y_{n_k,\gamma})$ converges to $(t_\gamma,t_\gamma',x_\gamma,x_\gamma')$ as $k$ goes to infinity. Combining this latter result with $\limsup_{n\rightarrow\infty}(|t_{n,\gamma}-s_{n,\gamma}|^2 + |x_{n,\gamma}-y_{n,\gamma}|^2)=0$, which follows from \eqref{E:UniqProof3}, we finally obtain \eqref{E:UniqProof4}. On the other hand, to prove \eqref{E:UniqProof5}-\eqref{E:UniqProof6}, notice that we have (recalling that $\tilde u-\tilde w$ is usc)
\begin{align*}
\tilde u(t_\gamma,x_\gamma) - \tilde w(s_\gamma,y_\gamma) - 2\gamma|x_\gamma|^2 \ &\leq \ \liminf_{n\rightarrow\infty} \Phi_{n,\gamma}(t_{n,\gamma},s_{n,\gamma},x_{n,\gamma},y_{n,\gamma}) \\
&\leq \ \limsup_{n\rightarrow\infty} \Phi_{n,\gamma}(t_{n,\gamma},s_{n,\gamma},x_{n,\gamma},y_{n,\gamma}) \\
&\leq \ \tilde u(t_\gamma,x_\gamma) - \tilde w(s_\gamma,y_\gamma) - 2\gamma|x_\gamma|^2.
\end{align*}
This implies that
\begin{align*}
\tilde u(t_\gamma,x_\gamma) - \tilde w(s_\gamma,y_\gamma) &= \lim_{n\rightarrow\infty} \bigg(\tilde u(t_{n,\gamma},x_{n,\gamma}) - \tilde w(s_{n,\gamma},y_{n,\gamma}) - n\frac{|t_{n,\gamma}-s_{n,\gamma}|^2}{2} - n\frac{|x_{n,\gamma}-y_{n,\gamma}|^2}{2} \bigg) \\
&\leq \liminf_{n\rightarrow\infty} \big(\tilde u(t_{n,\gamma},x_{n,\gamma}) - \tilde w(s_{n,\gamma},y_{n,\gamma})\big) \\
&\leq \limsup_{n\rightarrow\infty} \big(\tilde u(t_{n,\gamma},x_{n,\gamma}) - \tilde w(s_{n,\gamma},y_{n,\gamma})\big) = \tilde u(t_\gamma,x_\gamma) - \tilde w(s_\gamma,y_\gamma),
\end{align*}
which proves \eqref{E:UniqProof5} and \eqref{E:UniqProof6}.

Finally, we derive a useful inequality. More precisely, for any $\xi,\xi'\in\R^d$, from the maximum property $\Phi_{n,\gamma}(t_{n,\gamma},s_{n,\gamma},x_{n,\gamma}+d,y_{n,\gamma}+d') \leq \Phi_{n,\gamma}(t_{n,\gamma},s_{n,\gamma},x_{n,\gamma},y_{n,\gamma})$ we get
\begin{align}
\label{E:IneqI2delta}
&\tilde u(t_{n,\gamma},x_{n,\gamma} + d) - \tilde u(t_{n,\gamma},x_{n,\gamma}) - nd.(x_{n,\gamma}-y_{n,\gamma}) \notag \\
&\leq \ \tilde w(s_{n,\gamma},y_{n,\gamma} + d') - \tilde w(s_{n,\gamma},y_{n,\gamma}) - nd'.(x_{n,\gamma}-y_{n,\gamma}) \notag \\
&+ n \frac{|d-d'|^2}{2} + \gamma \big( |x_{n,\gamma} + d|^2 - |x_{n,\gamma}|^2 + |y_{n,\gamma} + d'|^2 - |y_{n,\gamma}|^2 \big).
\end{align}
\emph{Step 3.} Let us prove that, if $\gamma$ is small enough, then $t_\gamma<T$, so that $t_{n,\gamma},s_{n,\gamma}<T$, up to a subsequence. We proceed by contradiction, assuming $t_\gamma=T$. From \eqref{E:tilde_u-tilde_w>0} we obtain the contradiction (recalling that $\tilde u-\tilde w$ is usc)
\[
0 \ < \ \limsup_{n\rightarrow\infty}\big(\tilde u(t_{n,\gamma},x_{n,\gamma}) - \tilde w(s_{n,\gamma},y_{n,\gamma})\big) \ \leq \ \tilde u (T,x_\gamma) - \tilde w(T,x_\gamma) \ \leq \ 0.
\]
Consider, as in step~3, $(\hat t,\hat x)\in[0,T)\times\R^d$ such that $\tilde u(\hat t,\hat x) - \tilde w(\hat t,\hat x) =: \eta > 0$. Then, from the inequality $\Phi_{n,\gamma}(t_{n,\gamma},s_{n,\gamma},x_{n,\gamma},y_{n,\gamma}) \geq \Phi_{n,\gamma}(\hat t,\hat t,\hat x,\hat x)$, we obtain
\[
\tilde u(t_{n,\gamma},x_{n,\gamma}) - \tilde w(s_{n,\gamma},y_{n,\gamma}) \ \geq \ \tilde u(\hat t,\hat x) - \tilde w(\hat t,\hat x) - 2\gamma|\hat x|^2.
\]
Set $\gamma^* := (\tilde u(\hat t,\hat x) - \tilde w(\hat t,\hat x))/(4|\hat x|^2) \wedge 1$ if $|\hat x|^2 > 0$, and $\gamma^* := 1$ if $|\hat x|^2 = 0$. Then, for any $0 < \gamma \leq \gamma^*$, we have
\begin{equation}
\label{E:tilde_u-tilde_w>0}
\tilde u(t_{n,\gamma},x_{n,\gamma}) - \tilde w(s_{n,\gamma},y_{n,\gamma}) \ \geq \ \frac{\tilde u(\hat t,\hat x) - \tilde w(\hat t,\hat x)}{2} \ > \ 0,
\end{equation}
from which we obtain the contradiction (recalling that $\tilde u-\tilde w$ is usc)
\[
0 \ < \ \limsup_{n\rightarrow\infty}\big(\tilde u(t_{n,\gamma},x_{n,\gamma}) - \tilde w(s_{n,\gamma},y_{n,\gamma})\big) \ \leq \ \tilde u (T,x_\gamma) - \tilde w(T,x_\gamma) \ \leq \ 0.
\]
\emph{Step 4.} We shall apply the nonlocal Jensen-Ishii's lemma (see Lemma 1 in \cite{barlesimbert08}). To this end, let $\gamma\in(0,\gamma^*]$ and define
\[
\varphi_n(t,s,x,y) = n\frac{|t-s|^2}{2} + n\frac{|x-y|^2}{2} + \gamma\big(|x|^2 + |y|^2\big) - \Phi_{n,\gamma} (t_{n,\gamma},s_{n,\gamma},x_{n,\gamma},y_{n,\gamma}),
\]
for all $(t,s,x,y)\in\R^{2+2d}$ and for any $n\in\N\backslash\{0\}$. Then $(t_n,s_n,x_n,y_n) := (t_{n,\gamma},s_{n,\gamma},x_{n,\gamma},y_{n,\gamma})$ is a zero global maximum point  for $\tilde u(t,x) - \tilde w(s,y) - \varphi_n(t,s,x,y)$ on $[0,T]^2\times\R^{2d}$. Set
\begin{align*}
(p_n,q_n) \ &:= \ \bigg(\frac{\partial\varphi_n}{\partial t}(t_n,s_n,x_n,y_n),D_x\varphi_n(t_n,s_n,x_n,y_n)\bigg), \\
(-p_n',-q_n') \ &:= \ \bigg(\frac{\partial\varphi_n}{\partial s}(t_n,s_n,x_n,y_n),D_{y}\varphi_n(t_n,s_n,x_n,y_n)\bigg).
\end{align*}
Then, for any $\hat r > 0$, it follows from the nonlocal Jensen-Ishii's lemma that there exists $\hat\alpha(\hat r) > 0$ such that, for any $0 < \alpha \leq \hat\alpha(\hat r)$, we have: there exist sequences (to alleviate the notation, we omit the dependence of the sequences on $\alpha$) $(t_{n,k},s_{n,k},x_{n,k},y_{n,k})$ $\rightarrow$ $(t_n,s_n,x_n,y_n)$, $(t_{n,k},s_{n,k},x_{n,k},y_{n,k})\in[0,T)^2\times\R^{2d}$, $(p_{n,k},p_ {n,k}',q_ {n,k},q_ {n,k}') \rightarrow (p_n,p_n',q_n,q_n')$, matrices $N_{n,k},N_{n,k}'\in\S^d$, with $(N_{n,k},N_{n,k}')$ converging to some $(M_{n,\alpha},M_{n,\alpha}')$, and a sequence of functions $\varphi_{n,k}\in C^{1,2}([0,T]^2\times\R^{2d})$ such that:
\begin{enumerate}
\item[(i)] $(t_{n,k},s_{n,k},x_{n,k},y_{n,k})$ is a global maximum point of $\tilde u - \tilde w - \varphi_{n,k}$;
\item[(ii)] $\tilde u(t_{n,k},x_{n,k}) \rightarrow \tilde u(t_n,x_n)$ and $\tilde w(s_{n,k},y_{n,k}) \rightarrow \tilde w(s_n,y_n)$, as $k$ tends to infinity;
\item[(iii)] $(p_{n,k},q_{n,k},N_{n,k})\in\Pc^{2,+}\tilde u(t_{n,k},x_{n,k})$, $(p_{n,k}',q_{n,k}',N_{n,k}')\in\Pc^{2,-}\tilde w(s_{n,k},y_{n,k})$, and
\begin{align*}
(p_{n,k},q_{n,k}) \ &:= \ \bigg(\frac{\partial\varphi_{n,k}}{\partial t}(t_{n,k},s_{n,k},x_{n,k},y_{n,k}),D_x\varphi_{n,k}(t_{n,k},s_{n,k},x_{n,k},y_{n,k})\bigg), \\
(-p_{n,k}',-q_{n,k}') \ &:= \ \bigg(\frac{\partial\varphi_{n,k}}{\partial s}(t_{n,k},s_{n,k},x_{n,k},y_{n,k}),D_{y}\varphi_{n,k}(t_{n,k},s_{n,k},x_{n,k},y_{n,k})\bigg);
\end{align*}
\item[(iv)] The following inequalities hold (we denote by $I$ the $2d \times 2d$ identity matrix and by $D_{(x,y)}^2 \varphi_{n,k}$ the Hessian matrix of $\varphi_{n,k}$ with respect to $(x,y)$)
\begin{equation}
\label{E:N-N'}
-\frac{1}{\alpha} I \ \leq \
\begin{pmatrix}
N_{n,k} & 0 \\
0 & -N_{n,k}'
\end{pmatrix}
\ \leq \ D_{(x,y)}^2 \varphi_{n,k} (t_{n,k},s_{n,k},x_{n,k},y_{n,k}).
\end{equation}
\item[(v)] $\varphi_{n,k}$ converges uniformly in $\R^{2+2d}$ and in $C^{2}(B_{\hat r}(t_n,s_n,x_n,y_n))$ (where $B_{\hat r}(t_n,s_n,x_n,y_n)$ is the ball in $\R^{2+2d}$ of radius $\hat r$ and centered at $(t_n,s_n,x_n,y_n)$) towards $\psi_{n,\alpha} := R^\alpha[\varphi_n](\cdot,(p_n,p_n',q_n,q_n'))$, where, for any $\xi\in\R^{2+2d}$,
\[
R^\alpha[\varphi_n](z,\xi) := \sup_{|z'-z| \leq 1} \bigg\{\varphi_n(z') - \xi.(z'-z) - \frac{|z'-z|^2}{2\alpha}\bigg\}, \qquad \forall\,z\in\R^{2+2d}.
\] 
\end{enumerate}
Then, from Lemma \ref{L:viscosityJets} and the viscosity subsolution property to \reff{HJB2}-\reff{condterm2} of $\tilde u$, we have:
\begin{align*}
\rho \tilde u(t_{n,k}, x_{n,k}) - p_{n,k} -\sup_{a\in A} \bigg[ b(x_{n,k},a). q_{n,k} + \frac1{2} \text{tr}\big(\sigma\sigma\trans(x_{n,k},a)N_{n,k}\big) & \\
+ I_a^{1,\delta}(t_{n,k},x_{n,k},\varphi_{n,k}(\cdot,s_{n,k},\cdot,y_{n,k})) + I_a^{2,\delta}(t_{n,k},x_{n,k},q_{n,k},\tilde u) & \\
+ \tilde f\big(t_{n,k},x_{n,k},a\big)\bigg]& \ \leq \ 0.
\end{align*}
On the other hand, from the viscosity supersolution property to \reff{HJB2}-\reff{condterm2} of $\tilde w$, we have:
\begin{align*}
\rho \tilde w(s_{n,k},y_{n,k}) - p_{n,k}' - \sup_{a\in A} \bigg[ b(y_{n,k},a).q_{n,k}' + \frac1{2} \text{tr}\big(\sigma\sigma\trans(y_{n,k},a)N_{n,k}'\big) & \\
+ I_a^{1,\delta}(s_{n,k},y_{n,k},-\varphi_{n,k}(t_{n,k},\cdot,x_{n,k},\cdot)) + I_a^{2,\delta}(s_{n,k},y_{n,k},q_{n,k}',\tilde w) & \\
+ \tilde f\big(s_{n,k},y_{n,k},a\big)\bigg]& \ \geq \ 0.
\end{align*}
For every $k\in\N^*$, consider $a_k\in A$ such that
\begin{align}
\label{E:Ineq_a_k}
\rho \tilde u(t_{n,k}, x_{n,k}) - p_{n,k} - b(x_{n,k},a_k). q_{n,k} - \frac1{2} \text{tr}\big(\sigma\sigma\trans(x_{n,k},a_k)N_{n,k}\big) & \\
- I_{a_k}^{1,\delta}(t_{n,k},x_{n,k},\varphi_{n,k}(\cdot,s_{n,k},\cdot,y_{n,k})) - I_{a_k}^{2,\delta}(t_{n,k},x_{n,k},q_{n,k},\tilde u) & \notag \\
- \tilde f\big(t_{n,k},x_{n,k},a_k\big) & \ \leq \ \frac{1}{k}. \notag
\end{align}
From the compactness of $A$, we can suppose that $a_k\rightarrow a_\infty\in A$, up to a subsequence. Moreover, for every $a\in A$ we have
\begin{align}
\label{E:Ineq_a_k'}
\rho \tilde w(s_{n,k},y_{n,k}) - p_{n,k}' - b(y_{n,k},a).q_{n,k}' - \frac1{2} \text{tr}\big(\sigma\sigma\trans(y_{n,k},a)N_{n,k}'\big) & \\
- I_a^{1,\delta}(s_{n,k},y_{n,k},-\varphi_{n,k}(t_{n,k},\cdot,x_{n,k},\cdot)) - I_a^{2,\delta}(s_{n,k},y_{n,k},q_{n,k}',\tilde w) & \notag \\
- \tilde f\big(s_{n,k},y_{n,k},a\big) & \ \geq \ 0. \notag
\end{align}
Set $r^* := 2\sup_{(a,e)\in A\times(E\cap\{|e|\leq\delta\})} (|\beta(x^*,a,e)|\vee|\beta(y^*,a,e)|)$, where from \eqref{E:UniqProof4} we define $(x^*,y^*) := \lim_{n\rightarrow\infty}(x_n,y_n)$, and $\alpha^* := \hat\alpha(r^*)$. Notice that for all $n\in\N\backslash\{0\}$ we have $\sup_{(a,e)\in A\times(E\cap\{|e|\leq\delta\})} (|\beta(x_n,a,e)|\vee|\beta(y_n,a,e)|) < r^*$, up to a subsequence. Therefore, sending $k$ to infinity, we get $\varphi_{n,k}\rightarrow\psi_{n,\alpha}$, as $k$ tends to infinity, uniformly in $C^{2}(B_{r^*}(t_n,s_n,x_n,y_n))$ for any $0 < \alpha \leq \alpha^*$. Moreover, from assumption \textbf{(H$\lambda$)}(iii) we have
\begin{align*}
&\limsup_{k\rightarrow\infty} \int_{E\cap\{|e|\leq\delta\}} \big(\tilde u(t_{n,k},x_{n,k} + \beta(x_{n,k},a_k,e)) - \tilde u(t_{n,k},x_{n,k}) - \beta(x_{n,k},a_k,e).q_{n,k}\big) \lambda(a_{k},de) \\
&\leq \ \int_{E\cap\{|e|\leq\delta\}} \big(\tilde u(t_n,x_n + \beta(x_n,a_\infty,e)) - \tilde u(t_n,x_n) - \beta(x_n,a_\infty,e).q_n\big) \lambda(a_\infty,de).
\end{align*}
Therefore, from \eqref{E:Ineq_a_k} we obtain
\begin{align*}
\rho \tilde u(t_n, x_n) - p_n - b(x_n,a_\infty). q_n - \frac1{2} \text{tr}\big(\sigma\sigma\trans(x_n,a_\infty)M_{n,\alpha}\big) & \\
- I_{a_\infty}^{1,\delta}(t_n,x_n,\psi_{n,\alpha}(\cdot,s_n,\cdot,y_n)) - I_{a_\infty}^{2,\delta}(t_n,x_n,q_n,\tilde u) - \tilde f\big(t_n,x_n,a_\infty\big) & \ \leq \ 0. \notag
\end{align*}
A fortiori, if we take the supremum over $a\in A$ we conclude
\begin{align}
\label{tilde_u_a_n0}
\rho \tilde u(t_n, x_n) - p_n -\sup_{a\in A} \bigg[ b(x_n,a).q_n + \frac1{2} \text{tr}\big(\sigma\sigma\trans(x_n,a)M_{n,\alpha}\big) & \notag \\
+ I_a^{1,\delta}(t_n,x_n,\psi_{n,\alpha}(\cdot,s_n,\cdot,y_n)) + I_a^{2,\delta}(t_n,x_n,q_n,\tilde u) + \tilde f\big(t_n,x_n,a\big)\bigg]& \ \leq \ 0,
\end{align}
for any $0 < \alpha \leq \alpha^*$. On the other hand, letting $k$ to infinity in \eqref{E:Ineq_a_k'} for every fixed $a\in A$, and then taking the supremum, we end up with
\begin{align}
\label{tilde_w_a_n0}
\rho \tilde w(s_n,y_n) - p_n' - \sup_{a\in A} \bigg[ b(y_n,a).q_n' + \frac1{2} \text{tr}\big(\sigma\sigma\trans(y_n,a)M_{n,\alpha}'\big) & \notag \\
+ I_a^{1,\delta}(s_n,y_n,-\psi_{n,\alpha}(t_n,\cdot,x_n,\cdot)) + I_a^{2,\delta}(s_n,y_n,q_n',\tilde w) + \tilde f\big(s_n,y_n,a\big)\bigg]& \ \geq \ 0,
\end{align}
for any $0 < \alpha \leq \alpha^*$. Moreover, from \eqref{E:N-N'} we have
\begin{equation}
\label{E:M_nD^2alpha}
-\frac{1}{\alpha} I \ \leq \
\begin{pmatrix}
M_{n,\alpha} & 0 \\
0 & -M_{n,\alpha}'
\end{pmatrix}
\ \leq \ D_{(x,y)}^2 \psi_{n,\alpha} (t_n,s_n,x_n,y_n)
\end{equation}
and by direct calculation
\begin{equation}
\label{E:M_nD^2}
D_{(x,y)}^2 \psi_{n,\alpha} (t_n,s_n,x_n,y_n) \ = \ D_{(x,y)}^2 \varphi_n (t_n,s_n,x_n,y_n) + o(1), \qquad \text{as }\alpha\rightarrow0^+.
\end{equation}
\emph{Step 5.} From \eqref{tilde_u_a_n0}, for any $n$, consider $a_n\in A$ such that
\begin{align}
\label{tilde_u_a_n}
\rho \tilde u(t_n, x_n) - p_n - b(x_n,a_n). q_n - \frac1{2} \text{tr}\big(\sigma\sigma\trans(x_n,a_n)M_{n,\alpha}\big) & \notag \\
- I_{a_n}^{1,\delta}(t_n,x_n,\psi_{n,\alpha}(\cdot,s_n,\cdot,y_n)) 
- I_{a_n}^{2,\delta}(t_n,x_n,q_n,\tilde u) - \tilde f\big(t_n,x_n,a_n\big)& \ \leq \ \frac{1}{n}. 
\end{align}
On the other hand, from \eqref{tilde_w_a_n0} we deduce that
\begin{align}
\label{tilde_w_a_n}
\rho \tilde w(s_n,y_n) - p_n' -  b(y_n,a_n).q_n' - \frac1{2} \text{tr}\big(\sigma\sigma\trans(y_n,a_n)M_{n,\alpha}'\big) & \notag \\
- I_{a_n}^{1,\delta}(s_n,y_n,-\psi_{n,\alpha}(t_n,\cdot,x_n,\cdot)) - I_{a_n}^{2,\delta}(s_n,y_n,q_n',\tilde w) - \tilde f\big(s_n,y_n,a_n\big)& \ \geq \ 0. 
\end{align}
By subtracting \reff{tilde_w_a_n} to \reff{tilde_u_a_n}, we obtain:
\begin{align}
\label{ineq}
\rho (\tilde u(t_n,x_n) - \tilde w(s_n,y_n)) \ &\leq \ \frac{1}{n} +  p_n - p_n' + \Delta F_n + \Delta I_n^{1,\delta} + \Delta I_n^{2,\delta} \\
&\quad \ + b(x_n,a_n).q_n - b(y_n,a_n).q_n' \nonumber \\
&\quad \ + \frac1{2} \text{tr}\big(\sigma\sigma\trans(x_n,a_n)M_{n,\alpha} - \sigma\sigma\trans(y_n,a_n)M_{n,\alpha}'\big), \nonumber
\end{align}
where
\begin{align*}
\Delta F_n \ &= \ \tilde f\big(t_n,x_n,a_n\big) - \tilde f\big(s_n,y_n,a_n\big), \\
\Delta I_n^{1,\delta} \ &= \ I_{a_n}^{1,\delta}(t_n,x_n,\psi_{n,\alpha}(\cdot,s_n,\cdot,y_n)) - I_{a_n}^{1,\delta}(s_n,y_n,-\psi_{n,\alpha}(t_n,\cdot,x_n,\cdot)), \\
\Delta I_n^{2,\delta} \ &= \ I_{a_n}^{2,\delta}(t_n,x_n,q_n,\tilde u) - I_{a_n}^{2,\delta}(s_n,y_n,q_n',\tilde w).
\end{align*} 
We have
\[
p_n - p_n' \ = \ \frac{\partial\varphi_n}{\partial t}(t_n,s_n,x_n,y_n) + \frac{\partial\varphi_n}{\partial s}(t_n,s_n,x_n,y_n) \ = \ 0.
\]
By the uniform Lipschitz property of $b$ with respect to $x$, and \reff{E:UniqProof5}, we see that
\begin{align*}
&\lim_{n\rightarrow\infty} \big(b(x_n,a_n).q_n - b(y_n,a_n).q_n'\big) \\
&= \ \lim_{n\rightarrow\infty} \big(b(x_n,a_n).D_x \varphi_n (t_n,x_n,y_n) + b(y_n,a_n).D_y \varphi_n (t_n,x_n,y_n)\big) \ = \ 0.
\end{align*}
Regarding the trace term in \eqref{ineq}, by the uniform Lipschitz property of $\sigma$ with respect to $x$, \eqref{E:M_nD^2alpha}, \eqref{E:M_nD^2}, and \reff{E:UniqProof5}, we obtain
\begin{align*}
&\limsup_{n\rightarrow\infty} \limsup_{\alpha\rightarrow0^+} \text{tr}\big(\sigma\sigma\trans(x_n,a_n)M_{n,\alpha}-\sigma\sigma\trans(y_n,a_n)M_{n,\alpha}'\big) \ \leq \ 0.
\end{align*}
Moreover, from assumption \textbf{(HBC)} and \eqref{E:UniqProof5}-\eqref{E:UniqProof6}, we find
\beqs
\lim_{n\rightarrow\infty} |\Delta F_n| = 0.
\enqs
Concerning the integral term $\Delta I_n^{1,\delta}$, we have, for some $\vartheta',\vartheta''\in(0,1)$,
\begin{align*}
\Delta I_n^{1,\delta} \ &= \ \int_{E\cap\{|e|\leq\delta\}} \big[ D_x^2 \psi_{n,\alpha} (t_n,s_n,x_n + \vartheta'\beta(x_n,a_n,e),y_n)\beta(x_n,a_n,e).\beta(x_n,a_n,e) \\
&\quad \ + D_y^2 \psi_{n,\alpha} (t_n,s_n,x_n,y_n + \vartheta''\beta(y_n,a_n,e))\beta(x_n,a_n,e).\beta(x_n,a_n,e) \big] \lambda(a_n,de).
\end{align*}
Therefore, using \eqref{E:M_nD^2} we see that there exists a positive constant $C_n'$, depending only on $(x_n,y_n)$, the Lipschitz constant of $\beta$, and on $\sup_{\vartheta',\vartheta''\in[0,1]}|D_x^2\varphi_n(t_n,s_n,x_n+\vartheta'\beta(x_n,a_n,e),y_n)|\vee|D_y^2\varphi_n(t_n,s_n,x_n,y_n+\vartheta''\beta(y_n,a_n,e))|$, such that
\begin{equation}
\label{E:I1deltaIneq}
\limsup_{\alpha\rightarrow0^+} |\Delta I_n^{1,\delta}| \ \leq \ C_n' \int_{E\cap\{|e|\leq\delta\}} \big(1\wedge|e|^2\big) \lambda(a_n,de).
\end{equation}
Finally, it remains to consider the integral term $\Delta I_n^{2,\delta}$. Integrating inequality \eqref{E:IneqI2delta}, with $d = \beta(x_n,a_n,e)$ and $d' = \beta(y_n,a_n,e)$, we find
\begin{align*}
I_{a_n}^{2,\delta}(t_n,x_n,q_n,\tilde u) \ &\leq \ I_{a_n}^{2,\delta}(s_n,y_n,q_n',\tilde w) + n \int_{E\cap\{|e|>\delta\}} \frac{|\beta(x_n,a_n,e) - \beta(y_n,a_n,e)|^2}{2} \lambda(a_n,de) \\
&\quad \ + \gamma\int_{E\cap\{|e|>\delta\}} \big(|x_n + \beta(x_n,a_n,e)|^2 - |x_n|^2\big) \lambda(a_n,de) \\
&\quad \ + \gamma\int_{E\cap\{|e|>\delta\}} \big(|y_n + \beta(y_n,a_n,e)|^2 - |y_n|^2\big) \lambda(a_n,de).
\end{align*}
Then, it follows from assumption \textbf{(HFC)}(ii) that there exists a positive constant $C''$, depending only on the function $\beta$, such that (recalling that by Cauchy-Schwarz inequality we have $|a+b|^2 - |a|^2 \leq |b|^2 + 2|a||b|$, $a,b\in\R$)
\begin{align}
\label{E:I2deltaIneq}
I_{a_n}^{2,\delta}(t_n,x_n,q_n,\tilde u) \ &\leq \ I_{a_n}^{2,\delta}(s_n,y_n,q_n',\tilde w) + nC''\frac{|x_n-y_n|^2}{2} \int_E \big(1\wedge|e|^2\big) \lambda(a_n,de) \notag \\
&\quad \ + \gamma C''\big(1 + |x_n|^2 + |y_n|^2\big)\int_E \big(1\wedge|e|^2\big) \lambda(a_n,de).
\end{align}
From assumption \textbf{(HFC)}(iii) we see that $\sup_{a\in A}\int_E (1\wedge|e|^2)\lambda(a,de) < \infty$. Moreover, from \eqref{E:UniqProof3} we have that $|x_n|^2 + |y_n|^2$ is bounded by a constant, independent of $n$ and $\gamma$. So that, enlarging the constant $C''$ appearing in \eqref{E:I2deltaIneq} if necessary, we find
\begin{align}
\label{E:I2deltaIneq2}
I_{a_n}^{2,\delta}(t_n,x_n,q_n,\tilde u) \ &\leq \ I_{a_n}^{2,\delta}(s_n,y_n,q_n',\tilde w) + nC''\frac{|x_n-y_n|^2}{2} + \gamma C''.
\end{align}
In conclusion, plugging \eqref{E:I2deltaIneq2} into \eqref{ineq}, we obtain
\begin{align}
\label{ineq2}
\rho (\tilde u(t_n,x_n) - \tilde w(s_n,y_n)) \ &\leq \ \frac{1}{n} +  p_n - p_n' + \Delta F_n + \Delta I_n^{1,\delta} \\
&\quad \ + nC''\frac{|x_n-y_n|^2}{2} + \gamma C'' \notag \\
&\quad \ + b(x_n,a_n).q_n - b(y_n,a_n).q_n' \nonumber \\
&\quad \ + \frac1{2} \text{tr}\big(\sigma\sigma\trans(x_n,a_n)M_{n,\alpha} - \sigma\sigma\trans(y_n,a_n)M_{n,\alpha}'\big). \nonumber
\end{align}
Then, taking the $\limsup_{\alpha\rightarrow0^+}$ in both sides of \eqref{ineq2} and using \eqref{E:I1deltaIneq}, we get
\begin{align}
\label{ineq3}
\rho (\tilde u(t_n,x_n) - \tilde w(s_n,y_n)) \ &\leq \ \frac{1}{n} +  p_n - p_n' + \Delta F_n + C_n' \int_{E\cap\{|e|\leq\delta\}} \big(1\wedge|e|^2\big) \lambda(a_n,de) \\
&\quad \ + nC''\frac{|x_n-y_n|^2}{2} + \gamma C'' \notag \\
&\quad \ + b(x_n,a_n).q_n - b(y_n,a_n).q_n' \nonumber \\
&\quad \ + \frac1{2} \limsup_{\alpha\rightarrow0^+} \text{tr}\big(\sigma\sigma\trans(x_n,a_n)M_{n,\alpha} - \sigma\sigma\trans(y_n,a_n)M_{n,\alpha}'\big). \nonumber
\end{align}
Now, taking the $\limsup_{\delta\rightarrow0^+}$ in both sides of \eqref{ineq3}, we deduce
\begin{align}
\label{ineq4}
\rho (\tilde u(t_n,x_n) - \tilde w(s_n,y_n)) \ &\leq \ \frac{1}{n} +  p_n - p_n' + \Delta F_n + nC''\frac{|x_n-y_n|^2}{2} + \gamma C'' \\
&\quad \ + b(x_n,a_n).q_n - b(y_n,a_n).q_n' \nonumber \\
&\quad \ + \frac1{2} \limsup_{\alpha\rightarrow0^+} \text{tr}\big(\sigma\sigma\trans(x_n,a_n)M_{n,\alpha} - \sigma\sigma\trans(y_n,a_n)M_{n,\alpha}'\big). \nonumber
\end{align}
Recall from \eqref{E:tilde_u-tilde_w>0} that $\rho (\tilde u(t_n,x_n) - \tilde w(s_n,y_n)) \geq \rho(\tilde u(\hat t,\hat x) - \tilde w(\hat t,\hat x))/2$. Therefore, taking the $\limsup_{n\rightarrow\infty}$ in \eqref{ineq4}, we conclude
\[
0 \ < \ \frac{\tilde u(\hat t,\hat x) - \tilde w(\hat t,\hat x)}{2} \ \leq \ \gamma C'',
\]
which is a contradiction for $\gamma$ small enough.
\ep


\vspace{9mm}


\small

\begin{thebibliography}{}

\bibitem{bar94} Barles G. (1994): ``Solutions de viscosit\'{e} des \'{e}quations de Hamilton-Jacobi'', {\it Math\'{e}matiques et Applications}, Springer-Verlag.

\bibitem{barles97} Barles G., Buckdahn R. and E. Pardoux (1997): ``Backward stochastic differential equations and integral-partial differential equations'', {\it Stochastics and Stochastics Reports}, \textbf{60}, 57-83.

\bibitem{barlesimbert08} Barles G. and C. Imbert (2008): ``Second-order elliptic integro-differential equations: viscosity solutions' theory revisited'', {\it Annales de l'Institut Henri Poincar\'{e}}, \textbf{25}, 567-585.

\bibitem{bechschw05} Becherer D. and M. Schweizer (2005): ``Classical Solutions to Reaction-Diffusion Systems for Hedging Problems with Interacting It\^o and Point Processes'', {\it Annals of Applied Probability}, \textbf{15}, 1111-1144.


\bibitem{brezis10} Brezis H. (2010): {\it Functional Analysis, Sobolev Spaces and Partial Differential Equations}, Springer.

\bibitem{cinlar_jacod_protter_sharpe80} \c{C}inlar E., Jacod J., Protter P., and M. J. Sharpe (1980): {\it Semimartingales and Markov Processes}, Probability Theory and Related Fields, \textbf{54}, 161-219.

\bibitem{choukcossopham13} Choukroun S., Cosso A., and H. Pham (2013): ``Reflected BSDEs with nonpositive jumps, and controller-and-stopper games'', preprint arXiv:1308.5511.

\bibitem{conffuhrman12} Confortola F. and M. Fuhrman (2012): ``Backward stochastic differential equations and optimal control of marked point processes'', to appear in {\it SIAM Journal on Control and Optimization}.


\bibitem{crepey11} Cr\'{e}pey S. (2011): ``About the Pricing Equations in Finance'', in {\it Paris-Princeton Lectures in Mathematical Finance 2010}, Lecture Notes in Mathematics, Springer, 63-203.

\bibitem{crepeymat08} Cr\'{e}pey S. and A. Matoussi (2008): ``Reflected and Doubly Reflected BSDEs with Jumps: A Priori Estimates and Comparison'', {\it Annals of Applied Probability}, \textbf{18}, 2041-2069.

\bibitem{cuchfilmayteich11} Cuchiero C., Filipovi\'{c} D., Mayerhofer E., and J. Teichmann (2011): ``Affine Processes on Positive Semidefinite Matrices'', {\it Annals of Applied Probability}, \textbf{21}, 397-463.

\bibitem{elkaroui_peng_quenez} El Karoui N., Peng S., and M. C. Quenez (1997): ``Backward stochastic differential equations in finance'', {\it Mathematical Finance}, \textbf{7}, 1-71.

\bibitem{filovschm11} Filipovi\'{c} D., Overbeck L., and T. Schmidt (2011): ``Dynamic CDO Term Structure Modeling'', {\it Mathematical Finance}, \textbf{21}, 53-71.

\bibitem{flemso} Fleming W. H. and Soner, H. M. (2006) : ``Controlled Markov Processes and Viscosity Solutions'', {\it Stochastic Modelling and Applied Probability} \textbf{25}, Springer,
New York.

\bibitem{fuhrmanpham13} Fuhrman M. and H. Pham (2013): ``Dual and backward SDE representation for optimal control of non-Markovian SDEs'', preprint arXiv:1310.6943.

\bibitem{hu_peng09} Hu M. and S. Peng (2009): ``$G$-L\'evy Processes under Sublinear Expectations'', preprint arXiv:0911.3533.


\bibitem{jac75} Jacod J. (1975): ``Multivariate Point Processes: Predictable Projection, Radon-Nikodym Derivatives, Representation of Martingales'', Z. Wahrsch. verw. Geb., \textbf{31}, 235-253.

\bibitem{jac79} Jacod J. (1979): {\it Calcul Stochastique et Probl\`emes de Martingales}, Lecture Notes in Mathematics,  Springer, Berlin.

\bibitem{jacprotter82} Jacod J. and P. Protter (1982): ``Quelques remarques sur un nouveau type d'\'{e}quations diff\'{e}rentielles stochastiques'', {\it S\'{e}minaire de Probabilit\'{e}s XVI}, Lecture Notes in Mathematics, Springer, \textbf{920}, 447-458.

\bibitem{jacshiryaev03} Jacod J. and A. N. Shiryaev (2003): {\it Limit Theorems for Stochastic Processes}, second edition, Springer, Berlin.


\bibitem{kaziposszhou12I} Kazi-Tani N., Possamai D., and C. Zhou (2012): ``Second Order BSDEs with Jumps, Part I: Formulation and Uniqueness'', preprint arXiv:1208.0757.

\bibitem{kaziposszhou12II} Kazi-Tani N., Possamai D., and C. Zhou (2012): ``Second Order BSDEs with Jumps, Part II: Existence and Applications'', preprint arXiv:1208.0763.

\bibitem{khalanpha13a} Kharroubi I., Langren\'e N., and H. Pham (2013): ``Discrete time approximation of fully nonlinear HJB equations via BSDEs with nonpositive jumps'', preprint arXiv:1311.4505.

\bibitem{khalanpha13b} Kharroubi I., Langren\'e N., and H. Pham (2013): ``A numerical algorithm for fully nonlinear HJB equations: an approch by control randomization'', preprint arXiv:1311.4503.

\bibitem{khaetal10} Kharroubi I., Ma J., Pham H., and J. Zhang (2010): ``Backward SDEs with constrained jumps and quasi-variational inequalities'', {\it Annals of Probability}, {\bf 38}, 794-840.

\bibitem{khapha12} Kharroubi I. and H. Pham (2012): ``Feynman-Kac representation for Hamilton-Jacobi-Bellman IPDEs", to appear on \emph{Annals of Probability}.


\bibitem{neufeld_nutz14} Neufeld A. and M. Nutz (2014): ``Nonlinear L\'evy Processes and their Characteristics'', preprint arXiv:1401.7253.

\bibitem{peng00} Peng S. (2000): ``Monotonic limit theorem for BSDEs and non-linear Doob-Meyer decomposition'', {\it Probability Theory and Related Fields}, {\bf 16}, 225-234.

\bibitem{peng06} Peng S. (2006): ``G-expectation, G-Brownian motion and related stochastic calculus of It\^o type'', Proceedings of 2005, Abel symposium, Springer.


\bibitem{pham09} Pham H. (2009). ``Continuous-time Stochastic Control and Optimization with Financial Applications'', Vol. \textbf{61}, Springer.

\bibitem{revuzyor99} Revuz D. and M. Yor (1999): ``Continuous Martingales and Brownian Motion'', Springer-Verlag.

\end{thebibliography}
\end{document}